\documentstyle[twoside,12pt]{article}
\input psfig
\newif\ifpctex
\pctexfalse
\ifpctex\hoffset=-lin\voffset=-lin\fi
\ifpctex\def\bb{bf}\else\newfont{\bb}{msbm10 scaled\magstep1}\fi

\newcommand{\A}{{\cal A}}

\newcommand{\F}{{\cal F}}

\newcommand{\N}{{\cal N}}

\newcommand{\R}{{\cal R}}

\newcommand{\X}{{\cal X}}

\newcommand{\ph}{\varphi}

\def\nz{\mbox{\bb N}}
\def\rz{\mbox{\bb R}}
\def\zz{\mbox{\bb Z}}
\def\cz{\mbox{\bb C}}

\newcommand{\liminfn}{\liminf_{n\to\infty}}

\newcommand{\qed}{\hfill$\Box$\linebreak\medskip\par}

\newcommand{\pr}{{\em Proof:\quad}}

\newcommand{\be}{\begin{enumerate}}
\newcommand{\ee}{\end{enumerate}}
\newcommand{\bdm}{\begin{displaymath}}
\newcommand{\edm}{\end{displaymath}}
\newcommand{\beq}{\begin{equation}}
\newcommand{\eeq}{\end{equation}}
\newcommand{\beqa}{\begin{eqnarray}}
\newcommand{\eeqa}{\end{eqnarray}}
\newcommand{\beqas}{\begin{eqnarray*}}
\newcommand{\eeqas}{\end{eqnarray*}}

\catcode`@=11 \catcode`!=11

\expandafter\ifx\csname fiverm\endcsname\relax
  \let\fiverm\fivrm
\fi
  
\let\!latexendpicture=\endpicture 
\let\!latexframe=\frame
\let\!latexlinethickness=\linethickness
\let\!latexmultiput=\multiput
\let\!latexput=\put
 
\def\@picture(#1,#2)(#3,#4){%
  \@picht #2\unitlength
  \setbox\@picbox\hbox to #1\unitlength\bgroup 
  \let\endpicture=\!latexendpicture
  \let\frame=\!latexframe
  \let\linethickness=\!latexlinethickness
  \let\multiput=\!latexmultiput
  \let\put=\!latexput
  \hskip -#3\unitlength \lower #4\unitlength \hbox\bgroup}

\catcode`@=12 \catcode`!=12

\input pictex.tex

\catcode`@=11 \catcode`!=11
  
\let\!pictexendpicture=\endpicture 
\let\!pictexframe=\frame
\let\!pictexlinethickness=\linethickness
\let\!pictexmultiput=\multiput
\let\!pictexput=\put

\def\beginpicture{%
  \setbox\!picbox=\hbox\bgroup%
  \let\endpicture=\!pictexendpicture
  \let\frame=\!pictexframe
  \let\linethickness=\!pictexlinethickness
  \let\multiput=\!pictexmultiput
  \let\put=\!pictexput
  \let\input=\@@input   
  \!xleft=\maxdimen  
  \!xright=-\maxdimen
  \!ybot=\maxdimen
  \!ytop=-\maxdimen}

\let\frame=\!latexframe

\let\pictexframe=\!pictexframe

\let\linethickness=\!latexlinethickness
\let\pictexlinethickness=\!pictexlinethickness

\let\\=\@normalcr
\catcode`@=12 \catcode`!=12

\def\IMSmarkvadjust{0 pt}
\def\IMSmarkhadjust{0 pt}
\def\IMSmarkhpadding{0 pt}
\def\IMSpubltext{Published in modified form:}
\def\SBIMSMark#1#2#3{
 \font\SBF=cmss10 at 10 true pt
 \font\SBI=cmssi10 at 10 true pt
 \setbox0=\hbox{\SBF \hbox to \IMSmarkhpadding{\relax}
                Stony Brook IMS Preprint \##1}
 \setbox2=\hbox to \wd0{\hfil \SBI #2}
 \setbox4=\hbox to \wd0{\hfil \SBI #3}
 \setbox6=\hbox to \wd0{\hss
             \vbox{\hsize=\wd0 \parskip=0pt \baselineskip=10 true pt
                   \copy0 \break%
                   \copy2 \break%
                   \copy4 \break}}
 \dimen0=\ht6   \advance\dimen0 by \vsize \advance\dimen0 by 8 true pt
                \advance\dimen0 by -\pagetotal
	        \advance\dimen0 by \IMSmarkvadjust
 \dimen2=\hsize \advance\dimen2 by .25 true in
	        \advance\dimen2 by \IMSmarkhadjust

  \openin2=publishd.tex
  \ifeof2\setbox0=\hbox to 0pt{}
  \else 
     \setbox0=\hbox to 3.1 true in{
                \vbox to \ht6{\hsize=3 true in \parskip=0pt  \noindent  
                {\SBI \IMSpubltext}\hfil\break
                \input publishd.tex 
                \vfill}}
  \fi
  \closein2
  \ht0=0pt \dp0=0pt
 \ht6=0pt \dp6=0pt
 \setbox8=\vbox to \dimen0{\vfill \hbox to \dimen2{\copy0 \hss \copy6}}
 \ht8=0pt \dp8=0pt \wd8=0pt
 \copy8
 \message{*** Stony Brook IMS Preprint #1, #2. #3 ***}
}

\def\kiesps#1{\centerline{\psfig{figure=#1,height=2in}}} 
\def\kies#1{#1} 
\def\text#1{\mbox{#1}}
\def\Box{\sqcup\llap{$\sqcap$}}

\newtheorem {lemma}{Lemma}[section]
\newtheorem {theo}{Theorem}[section]
\newtheorem {bemerkung}{Remark}[section]
\newtheorem {remark}{Remark}[section]
\newtheorem {prop}{Proposition}[section]
\newtheorem {koro}{Corollary}[section]

\newtheorem {beispiel} {Example}[section]

\newcommand{\sect}[1]{\section{#1}\setcounter{equation}{0}}
\newcommand\Swiak{\'Swi\accent'30atek}
\newcommand\Swia{\'Swi\accent'30atek }
\renewcommand{\rho}{\varrho}

\newcommand\cSn[1]{c_{S_{#1}}}
\newcommand\cmSn[1]{z_n}
\newcommand\fSn[1]{f^{S_{#1}}}

\newcommand\st{\,\,;\,\,}

\begin{document}
\title
{Polynomial maps with a Julia set of positive Lebesgue measure:
Fibonacci maps\footnote{A mistake which was pointed out
to us by J.C. Yoccoz has been corrected in this version}}
\author{
Sebastian van Strien, University of Amsterdam, the Netherlands
\thanks{e-mail: strien at fwi.uva.nl.}
\\
Tomasz Nowicki, University of Warsaw, Poland
\thanks{supported by the NWO, KBN-GR91.}
}
\date{May 25, 1994}
\maketitle
\thispagestyle{empty}
\SBIMSMark{1994/3}{September 1994}{earlier received March, 1994}

%
\ifx\beginpic\undefined\else \fi

\chardef\oldatcat=\the\catcode`\@
\catcode`\@=11

\newskip\hsssglue \hsssglue=0pt plus 1fill minus 1fill \def\hsss{\hskip\hsssglue}

\newdimen\unitlength \newdimen\linethickness
\newdimen\@picheight \newdimen\@xdim \newdimen\@ydim \newdimen\@len \newdimen\@save
\newcount\@multicount \newcount\@xarg \newcount\@yarg
\newbox\@picbox \newbox\@mpbox

\font\tenln=line10     \font\tenlnw=linew10
\font\tencirc=lcircle10 \font\tencircw=lcirclew10
\font\smallln=linew10 scaled 483 

\def\thinlines{\let\linefont=\tenln \let\circlefont=\tencirc
  \linethickness=\fontdimen8\linefont}
\def\thicklines{\let\linefont=\tenlnw \let\circlefont=\tencircw
  \linethickness=\fontdimen8\linefont}
\thinlines

\def\beginpic(#1,#2)(#3,#4){\@picheight=#2\unitlength
  \setbox\@picbox=\hbox to#1\unitlength\bgroup\let\line=\@line
    \kern-#3\unitlength \lower#4\unitlength\hbox\bgroup\ignorespaces}
\def\endpic{\egroup\hss\egroup
  \ht\@picbox=\@picheight \dp\@picbox=\z@
  \leavevmode\box\@picbox}

\def\put(#1,#2)#3{\raise#2\unitlength\rlap{\kern#1\unitlength #3}\ignorespaces}

\def\multiput(#1,#2)(#3,#4)#5#6{\@multicount=#5
 \@xdim=#1\unitlength \@ydim=#2\unitlength \setbox\@mpbox=\hbox{#6}%
 \loop\ifnum\@multicount>0
   \raise\@ydim\rlap{\kern\@xdim \unhcopy\@mpbox}%
   \advance\@xdim#3\unitlength \advance\@ydim#4\unitlength
   \advance\@multicount\m@ne \repeat\ignorespaces}

\def\makebox(#1,#2)#3{\setbox\@picbox=\hbox to#1\unitlength{\hss#3\hss}%
  \@ydim=\ht\@picbox \advance\@ydim-\dp\@picbox
  \ht\@picbox=#2\unitlength \dp\@picbox=\z@
  \leavevmode\lower.5\@ydim\box\@picbox}

\newif\ifneg
\def\@line(#1,#2)#3{\@xarg=#1 \@yarg=#2 \@len=#3\unitlength \leavevmode
 \ifnum\@xarg<0 \reverseline \else \negfalse \@ydim=\z@\fi
 \ifnum\@xarg=0 \@vline
 \else\ifnum\@yarg=0 \@hline \else\@sline\fi\fi
 \ifneg\kern-\@len\else\@save=\@ydim\fi}
\def\reverseline{\negtrue \kern-\@len \@xarg=-\@xarg
 \@ydim=\@len \multiply\@ydim\@yarg \divide\@ydim\@xarg \@yarg=-\@yarg}

\def\@hline{\vrule height.5\linethickness depth.5\linethickness width\@len}
\def\@vline{\kern-.5\linethickness\vrule width\linethickness
  \ifnum\@yarg<0 height\z@ depth\else depth\z@ height\fi\@len
  \kern-.5\linethickness}

\def\@sline{\setbox\@picbox=\hbox{\linefont \count@=\@xarg \multiply\count@ 8
 \ifnum\@yarg>0 \advance\count@\@yarg \advance\count@-9
 \else \advance\count@-\@yarg \advance\count@ 55 \fi \char\count@}%
 \ifnum\@yarg<0 \@picheight=-\ht\@picbox \advance\@ydim\@picheight
 \else \@picheight=\ht\@picbox \fi
 \@xdim=\wd\@picbox \@save=\@ydim
 \loop\ifdim\@xdim<\@len \raise\@ydim\copy\@picbox
  \advance\@xdim\wd\@picbox \advance\@ydim\@picheight \repeat
 \advance\@xdim-\@len \kern-\@xdim
 \multiply\@xdim\@yarg \divide\@xdim\@xarg \advance\@ydim-\@xdim
 \raise\@ydim\box\@picbox}

\def\vector(#1,#2)#3{\@line(#1,#2){#3}%
 \ifnum\@xarg=0 \@vvector \else\ifnum\@yarg=0 \@hvector \else\@svector\fi\fi}
\def\@hvector{\ifneg\rlap{\linefont\char27}\else
 \smash{\llap{\linefont\char45}}\fi} 
\def\@vvector{\ifnum\@yarg<0 \raise-\@len\rlap{\linefont\char63}%
 \else\setbox\@picbox=\rlap{\linefont\char54}\advance\@len-\ht\@picbox
 \raise\@len\box\@picbox\fi}

\def\@svector{\setbox\@picbox=\hbox to\z@{\linefont
 \ifnum\@yarg<0 \count@=55 \@yarg=-\@yarg \else\count@=-9 \fi
 \ifneg\multiply\@xarg16 \multiply\@yarg2
 \else\hss 
  \ifnum\@xarg>2 \multiply\@xarg9 \multiply\@yarg2 \advance\count@29
  \else\ifnum\@yarg>2 \multiply\@xarg16 \multiply\@yarg9 \advance\count@-20
   \else\multiply\@xarg24 \multiply\@yarg3 \fi\fi\fi
  \advance\count@\@xarg \advance\count@\@yarg \char\count@
  \ifneg\hss\fi}
 \raise\@save\box\@picbox}

\def\disk#1{\@len=#1\unitlength \count@='160 \@diskcirc}
\def\circle#1{\@len=#1\unitlength \count@='140 \@diskcirc}
\def\@diskcirc{\setbox\@picbox=\hbox{\circlefont\char\count@}\@xdim=\wd\@picbox
 \leavevmode \ifdim\@len>15.499\@xdim \@bigdc \else \@smalldc\fi}
\def\@bigdc{\ifnum\count@<'160 \@bigcirc
 \else \@len=15\@xdim \@diskcirc\fi}
\def\@smalldc{{\advance\@len-.5\@xdim
 \loop\ifdim\@xdim<\@len \advance\count@\@ne \advance\@xdim\wd\@picbox\repeat
 \hbox{\circlefont\char\count@}}}
\def\@bigcirc{{\circlefont\count@=15
 \setbox\@picbox=\hbox{\char\count@}\@xdim=\wd\@picbox
 \ifdim\@len>2.5\@xdim \@len=2.5\@xdim\fi
 \advance\@len-.125\wd\@picbox
 \loop\ifdim\@xdim<\@len \advance\count@ 4 \advance\@xdim.25\wd\@picbox\repeat
 \@ydim=.5\@xdim \advance\@ydim.5\linethickness
 \setbox\@picbox=\vbox{\hbox{\char\count@\advance\count@-3\char\count@}%
  \nointerlineskip
  \hbox{\advance\count@\m@ne\char\count@\advance\count@\m@ne\char\count@}}%
 \kern-\@ydim\lower\@ydim\box\@picbox}}

\newif\ifovaltl \newif\ifovaltr \newif\ifovalbl \newif\ifovalbr
\ovaltltrue \ovaltrtrue \ovalbltrue \ovalbrtrue
\def\oval(#1,#2){\@xdim=#1\unitlength \@ydim=#2\unitlength
 {\circlefont \setbox\@picbox=\hbox{\char0}
 \ifdim\@xdim<\wd\@picbox \@xdim=\wd\@picbox\fi
 \ifdim\@ydim<\wd\@picbox \@ydim=\wd\@picbox\fi
 \@save=\@xdim \ifdim\@ydim<\@save \@save=\@ydim \fi
 \count@=39
 \loop \setbox\@picbox=\hbox{\char\count@}\ifdim\@save<\wd\@picbox
  \advance\count@-4 \repeat
 \setbox\strutbox=\hbox{\vrule height\ht\@picbox depth\dp\@picbox width\z@
   \kern\wd\@picbox}%
 \@save=.5\wd\@picbox \advance\@save-.5\linethickness
 \setbox0=\hbox to\@xdim{\ifovaltl\char\count@\else\strut\fi
  \kern-\@save\leaders\hrule height\ifovaltl\linethickness\else\z@\fi\hfil
  \leaders\hrule height\ifovaltr\linethickness\else\z@\fi\hfil\kern\@save
  \ifovaltr\advance\count@-3\char\count@\else\strut\fi\kern-\wd\@picbox}%
  \advance\count@\m@ne
 \setbox2=\hbox to\@xdim{\ifovalbl\char\count@\else\strut\fi
  \kern-\@save\leaders\hrule height\ifovalbl\linethickness\else\z@\fi\hfil
  \leaders\hrule height\ifovalbr\linethickness\else\z@\fi\hfil\kern\@save
  \ifovalbr\advance\count@\m@ne\char\count@\else\strut\fi\kern-\wd\@picbox}%
 \@save=\@ydim \advance\@save-\wd\@picbox \divide\@save 2
 \setbox\@picbox=\vbox{\box0\nointerlineskip
  \hbox to\@xdim{\vrule height\@save width\ifovaltl\linethickness\else\z@\fi
    \hfil\ifovaltr\vrule width\linethickness\kern-\linethickness\fi}%
  \nointerlineskip
  \hbox to\@xdim{\vrule height\@save width\ifovalbl\linethickness\else\z@\fi
    \hfil\ifovalbr\vrule width\linethickness\kern-\linethickness\fi}%
  \nointerlineskip\box2}%
  \@save=.5\@ydim \advance\@save.5\linethickness \leavevmode
  \kern-.5\@xdim \kern-.5\linethickness \lower\@save\box\@picbox}}

\def\cpic#1\endcpic{\vcenter{\hbox{\beginpic#1\endpic}}}


\newdimen\@xi \newdimen\@xii \newdimen\@xiii \newdimen\@xiv
\newdimen\@xpt \newdimen\@xoldpt
\newdimen\@yi \newdimen\@yii \newdimen\@yiii \newdimen\@yiv
\newdimen\@ypt \newdimen\@yoldpt
\def\squine(#1,#2,#3,#4,#5,#6){\setbox\@picbox\hbox{\tencirc q}%
 \global\@xoldpt=#1\unitlength \global\@yoldpt=#4\unitlength \kern\@xoldpt
 \@xi=\@xoldpt \@xii=#2\unitlength \@xiii=#3\unitlength
 \@yi=\@yoldpt \@yii=#5\unitlength \@yiii=#6\unitlength
 \squinerec
 \@xpt=#3\unitlength \@ypt=#6\unitlength \@addpoint
 \raise\@ypt\copy\@picbox}
\newif\iffar
\def\squinerec{\farfalse \testnear\@xi\@xiii \testnear\@yi\@yiii
 \iffar \decast \fi}
\def\testnear#1#2{\@save=#1\advance\@save-#2%
 \ifdim\@save<\z@ \@save=-\@save\fi \ifdim\@save>\p@ \fartrue \fi}
\def\decast{\@xpt=\@xi \advance\@xpt\@xii \divide\@xpt2
 \advance\@xii\@xiii \divide\@xii2
 \@xiv=\@xpt \advance\@xiv\@xii \divide\@xiv2
 \@ypt=\@yi \advance\@ypt\@yii \divide\@ypt2
 \advance\@yii\@yiii \divide\@yii2
 \@yiv=\@ypt \advance\@yiv\@yii \divide\@yiv2
 \begingroup\@xii=\@xpt \@xiii=\@xiv
  \@yii=\@ypt \@yiii=\@yiv \squinerec\endgroup
 \@xpt=\@xiv \@ypt=\@yiv \@addpoint
 \@xi=\@xiv \@yi=\@yiv \squinerec}
\def\@addpoint{
 \global\advance\@xoldpt-\@xpt \wd\@picbox=-\@xoldpt
 \raise\@yoldpt\copy\@picbox \global\@xoldpt=\@xpt \global\@yoldpt=\@ypt}

\catcode`\@=\oldatcat


\begin{abstract}
In this paper we shall show that there exists $\ell_0$
such that for each even integer $\ell\ge \ell_0$
there exists $c_1\in \rz$ for which
the Julia set of
$z\mapsto z^\ell+c_1$ has positive Lebesgue measure.
This solves an old problem.

{\bf Editor's note:} In 1997, it was shown by Xavier Buff that there was a
serious flaw 
leaving a gap in the proof. Currently (1999), the question of positive
measure Julia sets remains open. 
\end{abstract}

\tableofcontents

\sect{Introduction and statement of results}

Since the work of Julia and Fatou from the 1920's
there has been a continuous interest in the dynamics of rational maps.
One of the main objects of study
is the Julia set. This is the closure of the set of repelling
periodic points or -- equivalently -- the complement of the set of points
which have neighbourhoods on which the iterates of the
map form a normal family. It was shown back in the 20's that the Julia set
of a polynomial map is nowhere dense
and that its Julia set is the boundary of the set of points
whose iterates do not tend to infinity.
Thus it was natural to conjecture
that the Julia set of such maps have Lebesgue measure zero.
In this paper we shall show that this conjecture is false.
Given a map $f\colon \cz\to \cz$, let $\omega(z)$
be the set of accumulation points of the sequence
$z,f(z),f^2(z),\dots$.
\vskip 0.3cm

\noindent
{\bf Main Theorem} \\
{\em For each sufficiently large even integer $\ell$
there exists $c_1\in \rz$ such that the map
$f(z)=z^\ell+c_1$ has the following properties:
\begin{itemize}
\item the set $\omega(0)$ is a Cantor set with zero Lebesgue measure;
\item the set of points $z\in \cz$ for which
$\omega(z)$ is contained in $\omega(0)$ has positive Lebesgue measure;
\item the set of points whose forward iterates remain bounded
has no interior.
\end{itemize}
In particular, the Julia set of
$z\mapsto z^\ell+c_1$ has positive Lebesgue measure.
This map has the Fibonacci dynamics (to be defined in the next section).}

\bigskip
In other words, the Julia set from our example is `thin but heavy'.
In the picture below we have drawn the Julia set
of the unimodal polynomial $z\mapsto z^\ell+c_1$
when $\ell=16$ and $c_1=-1.04710851003600355\ldots\in \rz$
is chosen so that this map has the Fibonacci map dynamics.
To determine $c_1$, we have used a program of Tangerman which determines,
given $\ell$, the corresponding coefficient $c_1$ to any required precision,
\cite{Tan}. The Julia set is drawn, using a well-known program of Milnor.

\begin{figure}[htp]\hfil
\kiesps{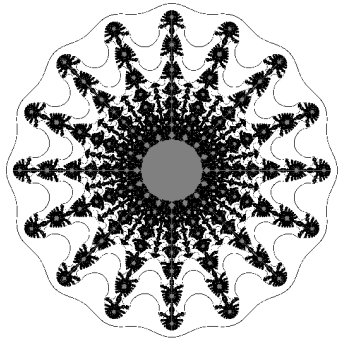} 
\caption{\protect{\small The Julia set of the Fibonacci map
$z\mapsto z^\ell+c_1$ where  $\ell=16$ and
$c_1=-1.04710851003600355\ldots$ (with some equipotentials
drawn in as well). We do not know whether
this value of $\ell=16$ is large enough for our theorem
to hold.}}
\end{figure}

We should point out that -- as far as we know --
this is the only example of a rational mapping with
a `heavy but thin' 
Julia set, but that several numerical
studies and mathematical results
already indicated that some rational
maps should have a nowhere dense Julia set
with positive Lebesgue measure, see \cite{Dou}
and \cite{Je}.
We did not make precise estimates on how
large we have to take $\ell$
but one should think of $\ell$ as being pretty large.

For entire functions `thin but heavy' Julia sets where constructed before
by McMullen, \cite{McM0}, see also Eremenko and Lyubich, \cite{EL}.
Shishikura, \cite{Sh1}, has shown that there exist (non-real) quadratic maps
whose Julia set has Hausdorff dimension two. From this he is able to conclude
that the boundary of the Mandelbrot set has also Hausdorff dimension two.
We expect that our methods may be helpful in improving this result
by showing that the Lebesgue measure of the Julia set of
such a map is positive. We also believe that this should imply that either
the boundary of the Mandelbrot set has positive Lebesgue measure
or that there are queer domains in the Mandelbrot set
(i.e., open sets in the parameter plane of conjugate non-hyperbolic maps).
\bigskip  

The inspiration for the Main Theorem
came from the analogous result on interval mappings
which was proved by the authors jointly with Gerhard Keller
and Henk Bruin:

\vskip0.3cm

\noindent
{\bf Theorem }\quad (The real case \cite{BKNS})\\
{\em For each sufficiently large $\ell\in \rz$
there exists $c_1\in \rz$ such that for the map
$f(x)=|x|^\ell+c_1$ the
set $\{z\st \omega(z)\subset \omega(c)\}$ is a
set of positive Lebesgue measure in $\rz$.}

\bigskip
\medskip
The maps we consider have Fibonacci-type dynamics. Such maps
were introduced by Hofbauer and Keller in the real context
as examples with very slow recurrence, see \cite{HK}. In the complex
setting these maps came up in \cite{BH}. In that paper
Branner and Hubbard study cubic maps with one critical point
escaping to infinity. They associate a tableau to each map
and the Fibonacci map has again the worst possible behaviour
w.r.t. the tableau rules. We believe that a related
cubic-like map also has a Julia set with positive Lebesgue measure:
consider a polynomial map of degree $\ell+1$ which preserves
the real line and with two critical points.
One of these critical points has order two and escapes to infinity and the
other has degree $\ell$. One can take this map so that its kneading
sequence is the same as that of the real
cubic map Branner and Hubbard considered.
Note that the Julia sets of these maps
are Cantor sets. A preliminary investigation suggests that for $\ell$ large
enough, these  Julia sets have positive Lebesgue measure also.

Let us now give a short outline of the ideas
needed for the proof of the Main Theorem.
Certain real estimates form one of the main ingredients for the proof
of the Main Theorem. Several of these
real estimates were proved already in \cite{BKNS},
and extend the estimates made in \cite{fibo}.
They follow from cross-ratio distortion results for
interval maps which were developed in the mid 1980's --
for an extensive overview of these tools, see the monograph
\cite{MS}. We would like to emphasize that some of the estimates
in this paper are really much stronger than
those from \cite{BKNS}. These cannot be derived
by applying Koebe and follow from two ideas.
Firstly, some estimates 
which show that maps, which are not like Moebius transformations,
satisfy improved Koebe estimates, see for example Proposition~\ref{rootl}.
Secondly, estimates which show that if one has 
a converging sequence of maps with an `almost neutral point' 
then -- up to a map with very small distortion --
one can compare their composition with the solution
of a particular differential equation,
see Theorems~\ref{miracle} and \ref{miracle2}.
For an abstract statement of this type of result, see
Theorem~\ref{comdv}. In this way we get
an asymptotic expression for some high
iterate of $f$ even though this limit
is extremely non-linear and one has extremely little
Koebe space. (Presumably similar estimates
should also work in a more general context.)

The second type of ingredients come from complex analysis:
using the Koebe Lemma and the Schwarz Lemma we are able to
show that the real estimates imply that
some return maps are polynomial-like in the sense
defined by Lyubich and Milnor \cite{LM}. This is 
an extremely important step because with this,
and because of the renormalization theory of Sullivan~\cite{S2}
and of McMullen~\cite{McM}, we can improve the real estimates.
The idea to apply renormalization theory also in
this case is due to Lyubich, see for example \cite{L4} and
\cite{L5}. We should emphasize, that many of the
estimates on this paper rely on this idea.

The third ingredient is that of a certain
induced map: this induced map is a very
natural consequence of the real analogue of Yoccoz's puzzle construction,
see Martens \cite{Mar}. This induced map is applied to the
so-called Fibonacci map which was `invented' in the real
context by Hofbauer and Keller \cite{HK} and
in the complex context by Branner and Hubbard \cite{BH}.
(For results on this maps, see \cite{LM}, \cite{fibo} and \cite{BKNS}.)
Interestingly, rather than using the Yoccoz partition of the Julia
set we found it advantageous in our proof to define a dynamical partition
which is based on two
circular curves (rather than by curves formed by equipotentials and
rays). In fact, this was inspired by  a kind of procedure
which was originated by Misha Lyubich and was first published in 
\cite{LM}. Something  a little similar is done in \cite{Sw2} and
\cite{GJ}, see also \cite{GS}. Because of this choice of domains, 
the proof that we get a polynomial-like map
requires some careful real estimates.
To get good estimates on the shape of these
curves we combine the real bounds (including the improved
estimate referred to above) with complex tools such as the
Koebe and Schwarz Lemma combined with
more delicate estimates which use the renormalization results.
In some sense this is the key part of this paper.
To get good estimates on the distortion of the induced map
we also use results which are better than those following
from Koebe, see Section~\ref{asytex}. More precisely,
we will decompose a high iterate of $f$ into a composition of
maps $\phi_i$. These maps $\phi_i$ send some preimage $z_i$ of $c$
into another preimage $z_{i-1}$.
By translating $z_i$ and $z_{i-1}$ to the origin and rescaling
small intervals to unit size, we can pretend that such a map
has a fixed point which is almost neutral
and in fact is close to a map of the type $z\mapsto z-z^3$. Then,
using a method which is reminiscent of Ecalle cylinders,
see \cite{Sh1}, we obtain a very good estimate for the composition
of the maps $\phi_i$. Putting all this together will give that
for $n$ and $\ell$ sufficiently large, certain iterates $f^{S_n-1}$
of $f$ map a neighbourhood of the critical value $c_1$ to a neighbourhood
of the critical point $c=0$ approximately as
$$c_1+z \mapsto \sqrt{M_\ell(z^2)},$$
where $M_\ell$ is a Moebius transformation which becomes more and
more degenerate as $\ell$ tends to infinity. The precise non-linearity
obtained from the composition of the Moebius transformation
with the quadratic map will give us our improved Koebe bounds.

The fourth and last ingredient in our proof is that of
a probabilistic (random walk) analysis of the behaviour
of typical points: this tool uses the language of martingales.
In fact, the result we use to apply these ideas was
proved by Gerard Keller, is stated in
Section~\ref{ransec}, and is
also one of the essential ingredients in \cite{BKNS}.
Many people have though of the idea to use
such random walk arguments. For example,
Guckenheimer and Johnson, use this terminology
in \cite{GJ}. Martens and van Strien discussed this
idea extensively in 1989, and this approach became
the motivation for the main result in \cite{Mar}.
Also, Luybich, Sutherland and Tangerman
performed computer experiments several years ago
to check the likelihood of our Main Theorem using
a random walk on the Yoccoz puzzle, \cite{LST}.
The first papers in which this idea was successfully
applied were \cite{fibo} and \cite{BKNS}.

Because this `random walk'
approach is not so usual in this subject,
we would like to explain these ideas by stating a simple
version of the result we use. Consider $r\in (0,1)$
and the map $F\colon [0,1)\to [0,1)$ which for $n\ge 1$
sends the interval $[r^{n+1},r^n)$ in an affine way to
$[0,r^{n-1})$ and which on $[r,1)$ is equal to the identity map.
Then it is not hard to show that
there exists $r_0>0$ such that for each $r\in (r_0,1)$
there exists a set $D$ of positive Lebesgue measure such that
$x\in D$ implies that $F^n(x)\to 0$ as $n\to \infty$.
This is not surprising: when $r_0$ is close to one,
a point in $[r_n,r_{n+1})$ moves with probability $1-r$
to the right and with probability $r$ to the left.
More precisely, the chance to move $i\in \{-1,0,1,\dots\}$
states from a given interval $[r^{n+1},r^n)$
is equal to $(1-r)r^{i+1}$ and so the expected drift is
$$\sum_{i\ge -1} i \cdot (1-r) \cdot r^{i+1}$$
which tends to infinity as $r\uparrow 1$.
With a random walk argument this implies that points typically move
to the states with higher index, i.e., to the origin.

In this paper, we will have a similar random walk
model. Here $F$ will be some iterate of $f(z)=z^\ell+c_1$
and the role of the intervals $[r^n,r^{n+1})$ will be replaced
by some nested sequence of annuli in the complex plane
surrounding $0\in \cz$.
If we are able to show that there is a set $D\subset \cz$
of positive Lebesgue measure such that $z\in D$ implies that
$F^n(z)\to 0$ as $n\to \infty$, then it follows that
points in $D$ are not in the basin of $\infty$. Since the
map $f$ will be chosen in such a way that it has no
periodic attractors, we obtain that $D\subset J(f)$.
Hence $J(f)$ has positive Lebesgue measure!
As we shall see, however, the random walk model is
in this case considerably more subtle then
in the previous one-dimensional model. One important difference
-- which also explains the difficulty to
get conclusive evidence from the numerical experiments
-- is that the probability to go `further away from zero'
is certainly not small. However, as we shall
see at the last section of this paper, the `probability'
to move $i$ `states' closer to $0$ (in one step
of the induced map $F$) is roughly of the order $\frac{1}{i^2} e^{-i/\ell}$,
where $\ell$ is the order of the critical point.
This implies that the expected drift is equal to
$$\sum_i i \cdot \frac{1}{i^2} e^{-i/\ell}.$$
This sum is of the order $\log(\ell)$ and therefore
grows relatively slowly with $\ell$.
Therefore one might have to take an extremely large $\ell$
to offset all types of constants and get a proper
drift towards $0$.
We shall elaborate on this issue in the last section of this paper.
\bigskip

One of the main reasons why one is interested in the Lebesgue
measure of the Julia set is the Measurable Riemann Mapping Theorem.
Indeed, as became apparent through Sullivan's work,
one way to solve the well-known stability conjecture that
the set of polynomials which are structurally stable
form a dense set is through the Measurable Riemann Mapping
Theorem. The real quadratic version of this conjecture
was solved by \Swia (\cite{Sw2}) using this idea of Sullivan.
\Swia shows that any two real quadratic polynomials $P$ and $Q$
which are conjugate are quasi-symmetrically conjugate on the
real line. Using Sullivan's pullback method this implies that they
are quasiconformally conjugate on the Riemann sphere.
By considering the Beltrami coefficient of the conjugacy
and by using the Measurable Riemann Mapping Theorem
one obtains a path of polynomial maps $[0,1]\ni t\mapsto P_t$ with $P_0=P$
and $P_1=Q$ and where $t$ is defined on a neighbourhood of
$[0,1]$ in $\cz$. It is easy to show that this is only possible
if $P$ is structurally stable.

Now if the Julia set of these quadratic maps would have
zero Lebesgue measure, then one could substantially
simplify this proof: in this case
the conjugacy would not need to be quasi-symmetric on the
real line in order to obtain a quasiconformal extension.
More precisely, using the $\lambda$-Lemma, see \cite{MSS}
and also \cite{McM}[Theorem 4.7],
there would be a quasiconformal conjugacy
which is conformal outside the Julia set
(provided $J(f)$ does not disconnect the plane).
This was -- of course -- one of the motivations for
Lyubich and Shishikura's result that non-renormalizable quadratic maps
have a Julia set with zero Lebesgue measure, see \cite{L2} and \cite{Sh2}.
Lyubich's method is based on the combinatorial
pattern of the Yoccoz puzzle \cite{Y} and
of the moduli of annuli argument of \cite{BH}. This last argument
states that the modulus of the preimage by a quadratic map
of an annulus (perhaps of higher genus) is either
equal to or otherwise half the size of the modulus of the original annulus.

This last method breaks down entirely if we consider
polynomials with critical points of higher order.
As follows from the next result, in our specific example,
other methods can be used to show that the Fibonacci maps
from the Main Theorem
do not form counter examples to the stability conjecture:

\vskip 0.5cm

\noindent
{\bf Theorem A}\quad
{\em For each even $\ell\ge 4$ one has the following properties.
\begin{itemize}
\item
For each $\ell$ there exists a unique parameter $c_1\in \rz$
such that $f(z)=z^\ell+c_1$ has Fibonacci dynamics.
\item There are no measurable invariant linefields on $J(f)$.
\item There exists a nested sequence of discs $D_n$
centered at the origin and disjoint topological discs $D_n^0, D_n^1$
which are compactly contained in $D_n$
such that the maps
$$R_n\colon (D_n^0\cup D_n^1)\to D_n$$
defined for $z\in D_n^0\cup D_n^1$
by
$$R_n(z)=\{f^i(z)\st i>0\mbox{ is minimal with }f^i(z)\in D_n\}$$
converge -- up to scaling -- as $n\in 4\nz$ tends to infinity.
The discs $D_n$ are chosen as in Section~\ref{sectqcr}.
\end{itemize}}

\vskip 0.5cm

Let us clarify the last part of this theorem.
Fix $i\in \{1,2,3,4\}$. Then the sequence
$R_n\colon (D_n^0\cup D_n^1)\to D_n$
converges as $n\in 4\nz+i$ tends to infinity.
To say that such a sequence of maps converges is perhaps
unclear because the domains of the maps vary. However,
$D_n$ is a Euclidean disc and -- as we shall see in
Section~\ref{sectqcr}, $R_n|D^0_n\to D_n$
is a branched covering onto (with a single $\ell$ fold
branching point) and $R_n|D^1_n\to D_n$ is a diffeomorphism.
So if we take $\Lambda_n$ the scaling map from $D_n$ to the
unit disc, the two
inverses $\Lambda_n\circ R_n^{-1}\circ \Lambda_n^{-1}$
become maps defined on the unit disc
(the inverse of the first map is $\ell$-valued).
If these two inverses converge then we say that the
above sequence of maps $R_n$ converge.

\medskip

In Theorem A we use the pullback method of Sullivan \cite{S2}
and McMullen's results on renormalization, see \cite{McM}.
We would like to thank Misha Lyubich for suggesting
this result to us. The first two statements of Theorem A
are standard, see Section~\ref{sectqcr}. 
We would like to thank Jacek Graczyk for
some very useful discussions on the third part of this
result.

\bigskip\medskip

For simplicity we shall denote the Lebesgue measure
of a set $A$ in $\rz^n$ by $|A|$. Moreover, if $a,b\in \rz$
then $(a,b)$ denotes the interval connecting $a$ and $b$.
For $a\in \rz$ we define
$$\hat a=-a$$
and $a^\#$ will denote either $a$ or $\hat a$ depending on the
context (for example depending on the parity of
some integer $n$). Finally, given two sequences of
positive real numbers $u_n$ and $v_n$ (depending also
on some parameter $\ell$), we write
$u_n\le C v_n$ if for each sufficiently large
$\ell$ there exists $n_0(\ell)$ so that
this inequality holds for  $n\ge n_0(\ell)$.
The same letter $C$ will be used throughout this paper
for several such universal constants.

The authors would like to thank Gerard Keller for allowing
us to use his result on random walks from
Section~\ref{ransec} which was also one of
the essential ingredients in \cite{BKNS}.
We also would like to thank Misha Lyubich for some very helpful
discussions on the renormalization theory of polynomial-like maps.
It is a pleasure to thank Adrien Douady
and Jean-Christoff Yoccoz for helpful discussions.
During a discussion with Yoccoz
a mistake was found in the last part of a previous version of this paper.
To fix this, we had to develop the improved
Koebe estimates from Section~\ref{asytex}. Folkert Tangerman's
notes on a method of McMullen's to get renormalization results
were very useful. Discussions with Jacek Graczyk on this
and other aspects of the paper were very much appreciated.
We have included some computer generated pictures of the
Julia set and the Yoccoz puzzle.
These were made using programs written
by Folkert Tangerman, Scott Sutherland and Misha Lyubich.


\sect{Combinatorial properties}
\label{seccomb}
As is well-known, see \cite{HK} or \cite{LM} and also
\cite{MS},
the Fibonacci map is a non-renormalizable
unimodal interval map for which the closure of the
forward orbit of the critical point $c$ is a
minimal Cantor set $\omega(c)$. In this section we
want to construct `by hand' this Fibonacci map. The main reason for
doing this, is that  it also gives a nice sequence
of induced maps, and a good covering of $\omega(c)$.
This covering will be used in Section~\ref{secqsr}
to show that $\omega(c)$ is a Cantor set
of `bounded geometry', provided $\ell$ is large.

\subsection{Construction of the Fibonacci map}
Rather than giving the kneading invariant of the map,
or its kneading map, and check that it satisfies some admissibility
conditions, we shall construct by hand a topological version of the map.
This will be done inductively, and at the same time
we shall construct a partition of the interval
which is the real analogue of the Yoccoz puzzle,
see \cite{Y} and also see \cite{LM} in this context.
This partition was also used by Martens \cite{Mar}
and Keller and Nowicki \cite{fibo}. Something similar was done
in \cite{GJ}. A complex extension of
the return maps which we construct in this section
are crucial in the remainder of this paper.

Let us first introduce some notation.
If $I\subset \rz$ is a bounded closed interval
then we say that $f\colon I \to I$ is a unimodal map if $f$
is continuous, $f$ has a unique extremal point
$c\in I$ (which is a minimum)
and $f(\partial I)\subset \partial I$.
If $I$ is unbounded then we require that $I=\rz$ and
replace the last condition by $f(x)\to \infty$ as $|x|\to \infty$.
Hence for each $x\ne c$ there exists
$\hat x\ne x$ such that $f(x)=f(\hat x)$.
We define $x_n$ to be $f^n(x)$ and for simplicity we shall
often assume that $c$ is equal to $0$.

Let us start by taking a unimodal map $f_0\colon I\to I$,
and assume that $f_0$ has an orientation reversing fixed point $q$
and a minimum at $c\in \text{int}(I)$ (and so $c_1<c$). Assume that
$c_2\in (\hat q,1)$ and $c_3\in (c,\hat q)$.
We shall modify $f_0$ in $(q,\hat q)$ repeatedly to suit our needs.

\kies{
\begin{figure}[htp]
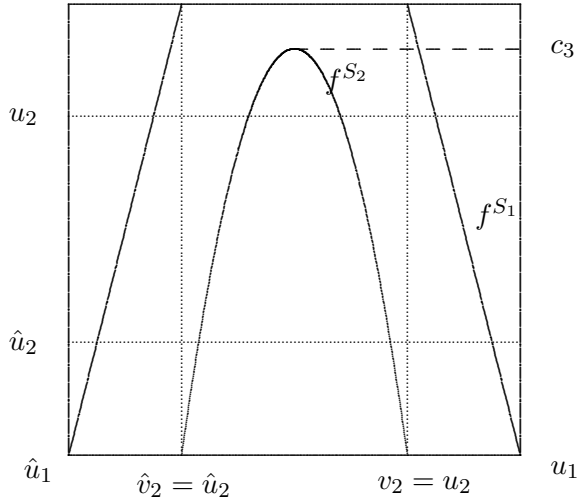
 \hfil
\beginpicture
\dimen0= 3cm \setcoordinatesystem units <\dimen0,\dimen0>
\setplotarea x from -1 to 1, y from -1 to 1
\setlinear \plot -1 -1 -1 1 1 1 1 -1 -1 -1 /
\setdots <1.5pt>
\plot -0.5 -1 -0.5 1 /
\plot 0.5 -1 0.5 1 /
\plot -1 0.5 1 0.5 /
\plot -1 -0.5 1 -0.5 /
\setdashes
\plot 0 0.8 1 0.8 /
\setsolid
\plot -1 -1 -0.5 1 /
\plot 0.5 1 1 -1 /
\setquadratic \plot -0.5 -1 0 0.8 0.5 -1 /
\put {\small $c_3$} [l]   <4mm,0mm> at  1   0.8
\put {\small $\hat u_1$} [l]        <-6mm,-2mm> at -1   -1
\put {\small $\hat v_2=\hat u_2$} [l]   <-6mm,-4mm> at -0.5 -1
\put {\small $v_2=u_2$} [l]   <-4mm,-4mm> at 0.5  -1
\put {\small $\hat u_2$} [l]   <-8mm,0mm> at -1 -0.5
\put {\small $u_2$} [l]        <-8mm,-0mm> at -1  0.5
\put {\small $u_1$} [l]   <4mm,-2mm> at  1   -1
\put {\small $f^{S_2}$} [l]    <-2mm,2mm> at  0.2  0.6
\put {\small $f^{S_1}$} [l]    <-3mm,-4mm> at  0.9   0.2
\endpicture
\caption[ ]{{\small The return map $R_2\colon U_2\to U_2$.
Here $u_0=q$, $u_1=\hat q$, $S_1=2$ and $S_2=3$.}}
\end{figure}

}

Define $S_0=1$, $S_1=2$ and $S_k=S_{k-1}+S_{k-2}$ for
$k\ge 2$. Let us define
$$u_0=q\text{ and }u_1=\hat q$$
and consider the first return map $R_2$ of $f_0$ to
$$U_2=(q,\hat q)=(\hat u_1,u_1).$$
Since $c_3\in U_2$ this first return map
consists of three branches:
two diffeomorphic ones
$U_2^1,U_2^2$ where
the return time is equal to $S_1=2$
(these intervals are symmetric)
and one which is defined on a `central' interval $U_2^0$
containing $c$ on which the map has a fold and on
which the return time is equal to $S_2=3$.

Now we modify $f_0$ on the central interval
$U_2^0$ (i.e., we keep $f_0$ the same outside
this interval) and call the new map $f_1$.
We do this so that $R_2(U_2^0)$ strictly contains the closure of
$U_2^0$. The reason this can be done, is because
there exists a neighbourhood around $c_1$ which is
disjoint from $U_2^0$ and which is mapped homeomorphically
onto $U_2^0$ by $f^{S_2-1}$.
Because $R_2$ is a first return map,
this modification does not affect $R_2|(U_2\setminus U_2^0)$.
This implies that $f_1^{S_2}(U_2^0)$ contains
one component of $U_2\setminus U_2^0$ and that
$f_1^{S_2}(c)$ is contained in the other component
(which we will call $U_2^1$). We should emphasize
that we have complete freedom where inside $U_2^1$ to choose $f_1^{S_2}(c)$.
Now we let
$$U_2^1=(u_2,x_2)\text{ and }U_2^0=(\hat v_2,v_2)$$
where we make the choices so that $u_2$ is the endpoint of $U_2^1$ which is
closer to $c$ and so that $u_2$ and $v_2$ are equal.
Note that $R_2(c)\in U_2^1$ and that $R_2(u_2)=u_1$.

\kies{
\begin{figure}[htp]
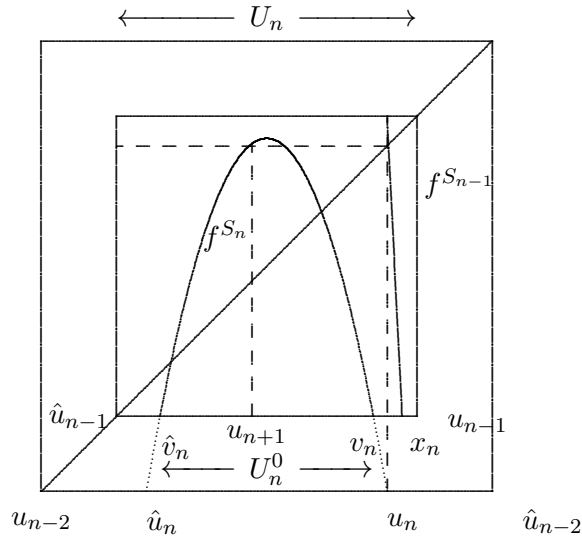
 \hfil \beginpicture
\dimen0=2cm \setcoordinatesystem units <\dimen0,\dimen0>
\setplotarea x from -1.5 to 1.5, y from -1.5 to 1.5
\setlinear \plot -1 -1 -1 1 1 1 1 -1 -1 -1 /
\plot -1.5 -1.5 -1.5 1.5 1.5 1.5 1.5 -1.5 -1.5 -1.5 /
\plot -1.5 -1.5 1.5 1.5 /
\plot 0.8 1 0.9 -1 /
\setdashes
\plot -0.1 -1 -0.1 0.8 /
\plot 0.8 -1.5 0.8 1 /
\plot 0.8 0.8 -1 0.8  /
\setsolid
\setquadratic \plot -0.709 -1 0 0.85 0.709 -1 /
\setdots <1pt>
\setquadratic \plot -0.8 -1.5 0 0.85 0.8 -1.5 /
\put {\small $v_n$} [l]        <-3mm,-4mm> at 0.7  -1
\put {\small $\hat v_n$} [l]   <0mm,-4mm> at -0.7  -1
\put {\small $u_n$} [l]        <0mm,-4mm> at 0.8  -1.5
\put {\small $\hat u_n$} [l]        <0mm,-4mm> at -0.8  -1.5
\put {\small $\hat u_{n-1}$} [l]        <-9mm,0mm> at -1   -1
\put {\small $u_{n+1}$} [l]        <-4mm,-2.5mm> at -0.07  -1
\put {\small $x_n$} [l]        <1mm,-4mm> at 0.9  -1
\put {\small $u_{n-1}$} [l]   <4mm,-1mm>  at  1   -1
\put {\small $u_{n-2}$} [l]        <-4mm,-4mm> at -1.5 -1.5
\put {\small $\hat u_{n-2}$} [l]   <4mm,-4mm> at 1.5 -1.5
\put {\small $f^{S_n}$} [l]    <-2.5mm,4mm> at  -0.3   0
\put {\small $f^{S_{n-1}}$} [l]    <1.3mm,7mm> at  0.98   0.2
\betweenarrows {\small $U_n^0$} <0pt,0mm> from -0.709 -1.35 to 0.709 -1.35
\betweenarrows {\small $U_n$}  <0pt,3mm> from -1 1.5 to 1 1.5
\endpicture
\vskip 2mm
\caption[ ]{\protect{\label{nsf_top2}}
{\small $R_n$ on the central branch $U_n^0=(v_n,\hat v_n)$
can be extended to $U_{n+1}=(u_n,\hat u_n)$ as shown.
$R_n$ is a surjection from $U_n^1=(v_n,x_n)$ onto $U_n$
and can also be extended to a monotone map onto $U_{n-1}$.
Moreover, $u_{n+1}$ is defined to be the point in $U_n^0$
for which $f^{S_n}(u_{n+1})=u_n$ and which is on the same side of
$c$ as $c_{S_{n+1}}$.}}
\end{figure}

}

We continue with the construction inductively.
So assume that $f_n$, $U_n=(u_{n-1},\hat u_{n-1})$,
$U_n^1=(x_n,u_n)$, $U_n^0=(v_n,\hat v_n)$,
$R_n\colon (U_n^0\cup U_n^1)\to U_n$
are already constructed for $n < N$.
Here we label these points so that
$u_n$ and $v_n$ are on the same side of $c$ and so that $u_n \in (x_n,c)$.
Also assume that
\begin{description}
\item[a)]  $R_n\colon U_n^1\to U_n$ is a homeomorphism,
$R_n\colon U_n^0\to U_n$ has one extremum
and $R_n(\partial U_n^0)\subset \partial U_n$;
\item[b)] $R_n(U_n^0) \supset U_n^0$;
\item[c)]  Moreover, $R_n(c)$ is contained in the other set $U_n^1$
and we have complete freedom where to place $R_n(c)$ in $U_n^1$
by changing $f_n$ inside $U_n^0$; but, we shall choose it so that
\item[d)]  $R_n \circ R_n(c) \in U_n^0$.
\end{description}
In Figure~\ref{nsf_top2} the graph of $R_n\colon U_n\to U_n$ is drawn
over the intervals $U_n^0$ and $U_n^1$ (and the extension of the central
interval is also depicted). Now we can proceed the construction
inductively as follows by taking
$U_N=(u_{N-1},\hat u_{N-1})$ and letting
$R_N$ to be the first return map to $U_N$.
Then define $U_N^{1,2} \subset U_{N-1}^0$ to be the two intervals
on which $R_N$ coincides with $R_{N-1}$
and which are mapped diffeomorphically by $R_N$ onto $U_N$.
The map $R_N$ has a `unimodal' branch over a
central interval $U_N^0$ (containing $c$),
and we modify $f_{N-1}$ in such a way that $R_N(U_N^0) \supset U_N^0$.
This proves properties a) and b) for $n=N$.
Of course, b) implies that $R_N(U_N^0)$ contains one of the sets $U_N^i$;
so let us call this set $U_N^2$. By the inductive hypothesis
c) we can even modify $f_{N-1}$ on $U_N^0$ such that
$R_N(c)$ is contained in the other set $U_N^1$
and we have complete freedom where to place $R_N(c)$ in $U_N^1$.
However, we shall choose it so that $R_N \circ R_N(c) \in U_N^0$.
This proves statements c) and d) for $n=N$.
Let us call the modified function $f_N$
and write
$$U_N^1=(x_N,u_N)\text{ and }U_N^0=(v_N,\hat v_N),$$
where
$u_N$ and $v_N$ lie on the same side of $c$ and $u_N \in (x_N,c)$.
We have $U_N^1 \cup U_N^0 \subset U^0_{N-1} \subset U_N$.

In this way we get sequences of points, intervals return maps and modified
functions $f_n$. Without loss of generality,
we may assume that $|U_n| \to 0$ as
$n \to \infty$, and as $f_n$ will only be modified on $U_k$ for $k \geq n$,
there exists a unimodal limit function $f$. In this way we have
shown how to construct a topological version of
a {\it Fibonacci interval map}.

\kies{
\begin{figure}[htp]
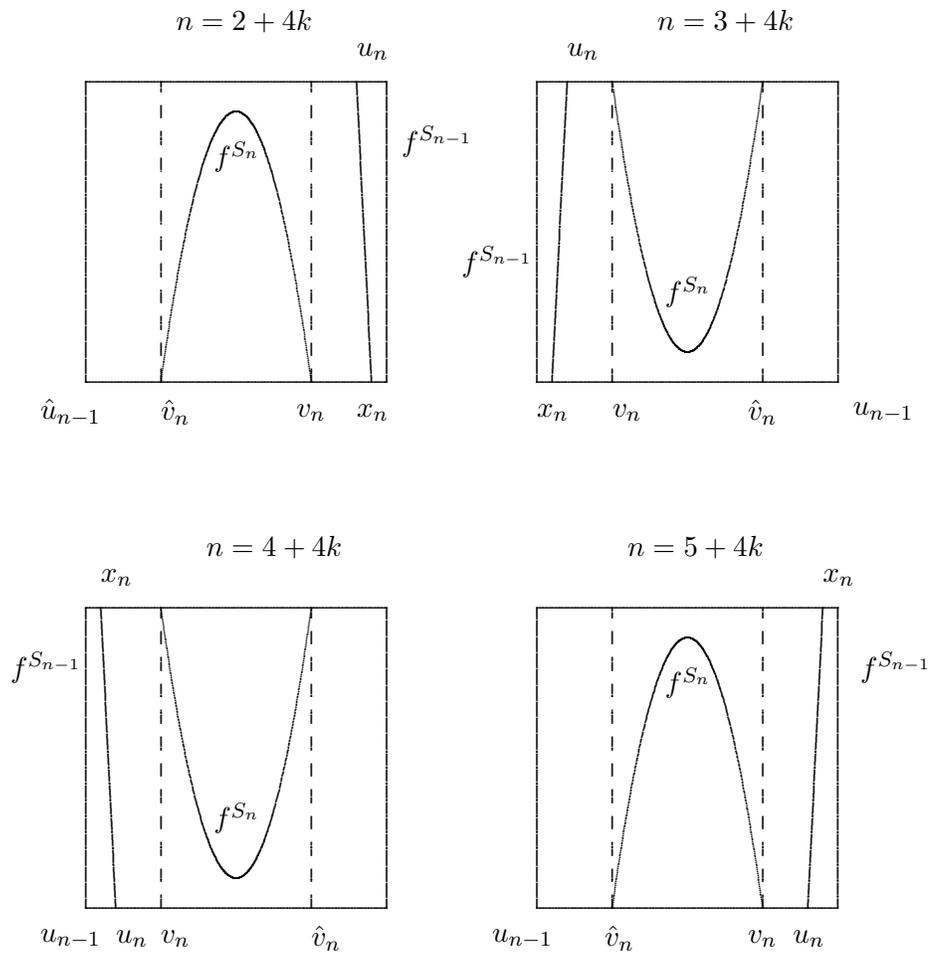
 \hfil \beginpicture
\dimen0=2cm
\setcoordinatesystem units <\dimen0,\dimen0> point at 0 3.5
\setplotarea x from -1 to 1, y from -1 to 1.3
\setlinear \plot -1 -1  -1 1  1 1  1 -1  -1 -1 /
\setlinear
\plot -0.9 1 -0.8 -1 /
\setdashes
\plot -0.5 -1 -0.5 1 /
\plot 0.5 -1 0.5 1 /
\setsolid
\setquadratic \plot -0.5 1 0 -0.8 0.5 1 /
\put {\small $u_{n-1}$} [l]        <-6mm,-4mm> at -1   -1
\put {\small $x_n$} [l]        <0mm,4mm> at -0.9 1
\put {\small $u_n$} [l]   <0mm,-4mm> at -0.8 -1
\put {\small $v_n$} [l]        <0mm,-4mm> at -0.5  -1
\put {\small $\hat v_n$} [l]   <0mm,-4mm>  at  0.5  -1
\put {\small $f^{S_n}$} [l]   <-3mm,4mm> at  0   -0.6
\put {\small $f^{S_{n-1}}$} [l]    <-12mm,4mm> at  -0.9   0.4
\put {\small $n=4+4k$} [l] <-4mm,8mm> at 0 1
\setcoordinatesystem units <\dimen0,\dimen0> point at -3 3.5
\setplotarea x from -1 to 1, y from -1 to 1.3
\setlinear \plot -1 -1  -1 1  1 1  1 -1  -1 -1 /
\setlinear
\plot 0.8 -1 0.9 1 /
\setdashes
\plot -0.5 -1 -0.5 1 /
\plot 0.5 -1 0.5 1 /
\setsolid
\setquadratic \plot -0.5 -1 0 0.8 0.5 -1 /
\put {\small $u_{n-1}$} [l]        <-6mm,-4mm> at -1   -1
\put {\small $x_n$} [l]        <0mm,4mm> at 0.9 1
\put {\small $u_n$} [l]        <-2mm,-4mm> at 0.8 -1
\put {\small $v_n$} [l]        <-2mm,-4mm> at 0.5  -1
\put {\small $\hat v_n$} [l]   <-1mm,-4mm>  at  -0.5  -1
\put {\small $f^{S_n}$} [l]   <-3mm,0mm> at  0   0.5
\put {\small $f^{S_{n-1}}$} [l]    <5mm,4mm> at  0.9   0.4
\put {\small $n=5+4k$} [l] <-8mm,8mm> at 0 1
\setcoordinatesystem units <\dimen0,\dimen0> point at 0 0
\setplotarea x from -1 to 1, y from -1 to 1.3
\setlinear \plot -1 -1  -1 1  1 1  1 -1  -1 -1 /
\setlinear
\plot 0.8 1 0.9 -1 /
\setdashes
\plot -0.5 -1 -0.5 1 /
\plot 0.5 -1 0.5 1 /
\setsolid
\setquadratic \plot -0.5 -1 0 0.8 0.5 -1 /
\put {\small $\hat u_{n-1}$} [l]    <-6mm,-4mm> at -1   -1
\put {\small $x_n$} [l]        <-2mm,-4mm> at 0.9 -1
\put {\small $u_n$} [l]        <0mm,+4mm> at 0.8 1
\put {\small $v_n$} [l]        <-2mm,-4mm> at 0.5  -1
\put {\small $\hat v_n$} [l]   <0mm,-4mm>  at  -0.5  -1
\put {\small $f^{S_n}$} [l]   <-3mm,0mm> at  0   0.5
\put {\small $f^{S_{n-1}}$} [l]    <4mm,4mm> at  0.9   0.4
\put {\small $n=2+4k$} [l] <-8mm,8mm> at 0 1
\setcoordinatesystem units <\dimen0,\dimen0> point at -3 0
\setplotarea x from -1 to 1, y from -1 to 1.3
\setlinear \plot -1 -1  -1 1  1 1  1 -1  -1 -1 /
\setlinear
\plot -0.9 -1 -0.8 1 /
\setdashes
\plot -0.5 -1 -0.5 1 /
\plot 0.5 -1 0.5 1 /
\setsolid
\setquadratic \plot -0.5 1 0 -0.8 0.5 1 /
\put {\small $u_{n-1}$} [l]        <2mm,-4mm> at 1   -1
\put {\small $x_n$} [l]        <-2mm,-4mm> at -0.9 -1
\put {\small $u_n$} [l]   <0mm,4mm> at -0.8 1
\put {\small $v_n$} [l]        <0mm,-4mm> at -0.5  -1
\put {\small $\hat v_n$} [l]   <-2mm,-4mm>  at  0.5  -1
\put {\small $f^{S_n}$} [l]   <-3mm,4mm> at  0   -0.6
\put {\small $f^{S_{n-1}}$} [l]    <-12mm,4mm> at  -0.9   -0.4
\put {\small $n=3+4k$} [l] <-4mm,8mm> at 0 1
\endpicture
\caption[ ]{\protect{\label{returnmaps}}
{\small The successive first return maps $R_n\colon U_n\to U_n$
for $n\ge 2$ (i.e., for $k\ge 0$; note however that $v_2=u_2$).}}
\end{figure}

}

Of course, this construction merely gives a continuous map.
However, using a general fullness result from the theory
of interval maps, in any reasonable family of unimodal maps
one can find maps with the same combinatorial properties:

\begin{lemma}\label{realizability}

Consider a family of $C^1$ unimodal maps $g_t
\colon [-1,1]\to [-1,1]$ such that $g_0$ has no periodic points
of period $>1$ and $g_1$ is surjective.
Then there exists a parameter $t'$ such that $g=g_{t'}$
is a Fibonacci map.

In particular, there exists for any $\ell\in 2\nz$
a Fibonacci map in the family type $z\mapsto z^\ell+t$, $t\in \rz$.
\end{lemma}
\pr This follows immediately from the fullness of such families,
see for example Section II.4 in \cite{MS}.
Indeed, this  fullness result implies
that there exists such a parameter $t'$ such that $g_{t'}$
has the same kneading invariant as the Fibonacci map constructed above.
\qed

\subsection{Topological properties of the Fibonacci map}
Let $\{S_k\}$ be the Fibonacci numbers, i.e. $S_0=1$, $S_1=2$ and
$S_k=S_{k-1}+S_{k-2}$.
We prove the following  properties of a Fibonacci map.
Define $U_n=(u_{n-1},\hat u_{n-1}),$
$U_n^0=(v_n,\hat v_n)$, $U_n^1=(u_n,x_n)$
and $R_n\colon U_n^0,U_n^1\to U_n$ be the
two branches of the return map as above. In the next lemma we show
how these return maps are related to $f$ and
what the orbit of the intervals $U_n^i$ look like.

\begin{lemma}
\label{indi}
Let $f$ be a Fibonacci map and take $n\in \nz$. Then
one has the following properties.
\begin{enumerate}
\item  $R_n|U_n^0$ coincides with $f^{S_n}|U_n^0$, and
$R_n|U_n^{1,2}$ coincides with $f^{S_{n-1}}|U_n^{1,2}$.
\item  $R_n(u_n)=f^{S_{n-1}}(u_n)=u_{n-1}$.
\item  $c_{S_n} \in (c_{S_{n-1}},\hat c_{S_{n-1}})$ and
$c_k \notin (c_{S_{n-1}},\hat c_{S_{n-1}})$ for $0<k<S_n$.
\item  Let $c_{-k}$ and $\hat c_{-k}$ be the points in $f^{-k}(c)$ which are
closest to $c$, then $c_{-S_n} \in (c_{-S_{n-1}},\hat c_{-S_{n-1}})$ and
$c_{-k} \notin (c_{-S_{n-1}},\hat c_{-S_{n-1}})$ for $0<k<S_n$.
\item  $f^{S_{n-1}}(U_n^0) \subset U_{n-1}^1$.
\item  For every $n \ge k \geq 2$, $f^i(U_n^0) \cap U_k \subset
U_k^0 \cup U_k^1$ for each $i=0,1,\ldots,S_n$
and $f^i(U_n^1) \cap U_k \subset U_k^0 \cup U_k^1$
for $i=0,1,\dots,S_{n-1}$, except if $n=k$
in which case $f^{S_n}(U_n^0) \subset U_n$
$f^{S_{n-1}}(U_n^1)\subset U_n$.
\item  $\omega(c) \cap U_n \subset U_n^0 \cup U_n^1$ for every $n \geq 2$.
\end{enumerate}
\end{lemma}
\pr
By the above construction these properties hold for $n=2$.
Since, by property b),\,\, $R_n|U_n^{1,2} = R_{n-1}|U_{n-1}^0$ and
$R_n|U_n^0 = R_{n-1}|U_{n-1}^1 \circ R_{n-1}|U_{n-1}^0$,
the first statement follows immediately from induction.
Statement 2) follows from the choice of $u_n$ and surjectivity of $R_n|U_n^1$.
Since $R_n$ is a first return map, $u_n$ is a preimage of $q$
and since $R_n|U_n^0 = f^{S_n}|U_n^0$, it follows that
$c_{S_n}$ are the successive closest returns to $c$.
So let us prove 4). The intervals $U_n^{1,2}$
contain precritical points in $f^{-S_{n-1}}(c)$, because
$R_n|U_n^{1,2} = f^{S_{n-1}}|U_n^{1,2}$ is surjection onto
$U_n$. Both branches are part of the central branch of the first return map
$R_{n-1}$. As $R_{n-1}(u_{n+1},\hat u_{n+1})$ does not contain $c$,
the interval
$U_n=(u_{n+1},\hat u_{n+1})$
contains no point in $\cup_{i=1}^{S_{n-1}} f^{-i}(c)$.
So $c_{-S_{n-1}}\cap U_n^{1,2}$ is a closest precritical point.
This proves 4).

In order to prove 5), observe that
since $U_n^0 \subset U_{n-1}^0$, property c)
implies $f^{S_{n-1}}(U_n^0) \cap U_{n-1}^1 \neq \emptyset$.
Therefore $f^{S_{n-1}}(U_n^0) \subset U_{n-1}^1$
since $f^{S_{n-2}}$ maps a neighbourhood of $U_{n-1}^1$
homeomorphically onto a neighbourhood of $U_{n-1}$
and so if $f^{S_{n-1}}(U_n^0)$ is not completely contained
in $ U_{n-1}^1$ then
$f^{S_n}(U_n^0)=f^{S_{n-2}}\circ f^{S_{n-1}}(U_n^0) \not\subset U_{n-1}$,
contradicting that $f^{S_n}|U_n^0$ is a branch of the
first return map to $U_n$. This proves 5).
Let us now prove statement 6) for $2\le k\le n< N$
by induction on $N$. For $N=3$ this statement is obvious.
Assume statement 6) holds for $2\le k\le n< N$
and let us show it also holds for $2\le k\le n\le N$.
Because $U_N^0 \subset U_{N-1}^0$, the induction assumption implies
$f^i(U_N^0) \cap U_k \subset U_k^0 \cup U_k^1$ for $0 \leq i \leq S_{N-1}$.
However $f^{S_{N-1}}(U_N^0) \subset U_{N-1}^1$, hence
$f^i(U_N^0) \cap U_k \subset
f^{i-S_{N-1}}(U_{N-1}^1) \subset U_k^0 \cup U_k^1$
for $S_{N-1} < i < S_N$ (because $S_N-S_{N-1}=S_{N-2}$ and using the
the second part of the
induction hypothesis fot the last inclustion). The other part of
statement 6) is proved similarly.

Since $U_n^0$ and $U_n^1$ are contained in the
interior of $U_{n-1}^0$, statement 7) follows immediately.
\qed

Now we will discuss the ordering of some
crucial dynamically defined points.
Firstly, we let $T_n\ni c_1$ be the maximal interval
for which $F^{S_n-1}|T_n$ is a diffeomorphism
and define
$$
y_n=\fSn{n}(\cSn{n+2})\ ,\quad y_n^f=f(y_n).
$$
Also define $w_n^f$, $r_n^f$ to be the points in $T_n$
to the left of
$c_1$ so that
$$f^{S_n-1}(w_n^f)=\hat u_{n-1}\mbox{ and }
f^{S_n-1}(r_n^f)=\hat u_{n-2}$$
(note that $w_n^f$ is not the image of a point $w_n\in [-1,1]$
so the notation is only to suggest that $w_n^f$ lies near $c_1$).
As before, let $x_n^f$ be the point in the interval $T_{n-1}$
for which $f^{S_{n-1}-1}(x_n^f)=\hat u_{n-1}$.
For simplicity we also write
$$d_n=c_{S_n}\text{ and }d_n^f=f(d_n).$$
Moreover, we shall write $z_n$ for one of the
two points in $f^{-S_n}(c)$ closest to $c$.

\begin{prop}\label{top2}
The points $u_n^f$, $d_n^f$, $x_n^f$, $y_n^f$, $w_n^f$ and $z_n^f$
are ordered as in the picture below
(we state the ordering near $c_1$ rather than near $c$
so that we do not need to be careful about on which side of $c$
these points lie).
\end{prop}

\begin{figure}[htp]
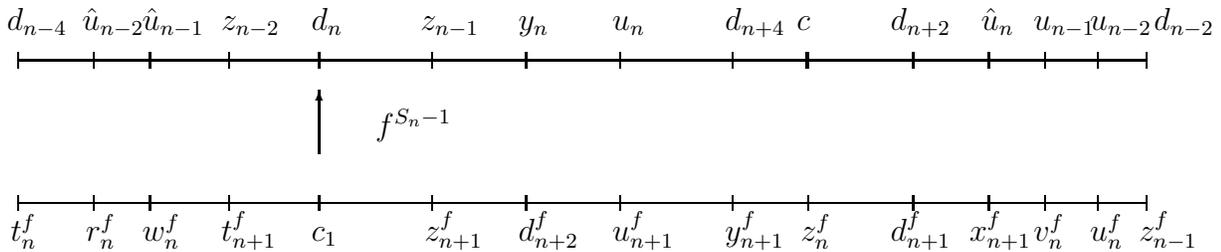

\vskip 0.2cm
\hbox to \hsize{\hss\unitlength=5mm
\beginpic(20,4)(-20,0) \let\ts\textstyle
\put(4,4){\line(-1,0){30}}
\put(4,3.8){\line(0,1){0.4}} \put(4.2,4.8){$d_{n-2}$}
\put(2.7,3.8){\line(0,1){0.4}} \put(2.5,4.8){$u_{n-2}$}
\put(1.3,3.8){\line(0,1){0.4}} \put(1,4.8){$u_{n-1}$}
\put(-0.2,3.8){\line(0,1){0.4}} \put(-0.4,4.8){$\hat u_{n}$}
\put(-2.2,3.8){\line(0,1){0.4}} \put(-2.8,4.8){$d_{n+2}$}
\put(-5,3.8){\line(0,1){0.4}}
\put(-5.05,3.8){\line(0,1){0.4}}
 \put(-5.3,4.8){$c$}
\put(-7,3.8){\line(0,1){0.4}} \put(-7.2,4.8){$d_{n+4}$}
\put(-10,3.8){\line(0,1){0.4}} \put(-10.2,4.8){$u_{n}$}
\put(-12.5,3.8){\line(0,1){0.4}} \put(-12.7,4.8){$y_{n}$}
\put(-15,3.8){\line(0,1){0.4}} \put(-15.3,4.8){$z_{n-1}$}
\put(-18,3.8){\line(0,1){0.4}} \put(-18.2,4.8){$d_{n}$}
\put(-20.4,3.8){\line(0,1){0.4}} \put(-20.6,4.8){$z_{n-2}$}
\put(-22.5,3.8){\line(0,1){0.4}} \put(-22.7,4.8){$\hat u_{n-1}$}
\put(-24,3.8){\line(0,1){0.4}} \put(-24.3,4.8){$\hat u_{n-2}$}
\put(-26,3.8){\line(0,1){0.4}} \put(-26.3,4.8){$d_{n-4}$}

\put(4,0.2){\line(-1,0){30}}
\put(4,0){\line(0,1){0.4}} \put(3.8,-0.8){$z_{n-1}^f$}
\put(2.7,0){\line(0,1){0.4}} \put(2.5,-0.8){$u_{n}^f$}
\put(1.3,0){\line(0,1){0.4}} \put(1,-0.8){$v_n^f$}
\put(-0.2,0){\line(0,1){0.4}} \put(-0.7,-0.8){$x_{n+1}^f$}
\put(-2.2,0){\line(0,1){0.4}} \put(-2.8,-0.8){$d_{n+1}^f$}
\put(-5,0){\line(0,1){0.4}} \put(-5.2,-0.8){$z_{n}^f$}
\put(-7,0){\line(0,1){0.4}} \put(-7.2,-0.8){$y_{n+1}^f$}
\put(-10,0){\line(0,1){0.4}} \put(-10.2,-0.8){$u_{n+1}^f$}
\put(-12.5,0){\line(0,1){0.4}} \put(-12.7,-0.8){$d_{n+2}^f$}
\put(-15,0){\line(0,1){0.4}} \put(-15.2,-0.8){$z_{n+1}^f$}
\put(-18,0){\line(0,1){0.4}} \put(-18.2,-0.8){$c_1$}
\put(-20.4,0){\line(0,1){0.4}} \put(-20.6,-0.8){$t_{n+1}^f$}
\put(-22.5,0){\line(0,1){0.4}} \put(-22.7,-0.8){$w_{n}^f$}
\put(-24,0){\line(0,1){0.4}} \put(-24.2,-0.8){$r_{n}^f$}
\put(-26,0){\line(0,1){0.4}} \put(-26.2,-0.8){$t_{n}^f$}
\put(-18,1.5){\vector(0,1){1.7}}
\put(-16.5,2){$f^{S_{n}-1}$ }
\endpic\hss}
\vskip 4mm
\caption[ ]{\label{nsf_top4}
{\small Points and their images under $f^{S_{n}-1}$.
Note that $c_1$ is the mimumum of $f\colon \rz\to \rz$.
The points $u_n,w_n,x_n$ are in the full orbit of the fixed point
$u_0$ whereas $d_n=f^{S_n}(c)$ and $y_n=f^{S_n}(d_{n+2})$
are forward iterates of $c$. The point $z_n$ is a point
in $f^{-S_n}(c)$ nearest to $c$. We should note that the position
of $u_{n-1}$ and $\hat u_{n-1}$ should be interchanged for $n$ even
(in that case $f^{S_n}(v_n)=\hat u_{n-1}$).}}
\end{figure}


\pr
The proof of these statements can be found in \cite{fibo}.
\qed

Next we shall show that the set
$\omega(c)$ of accumulation points $c,f(c),f^2(c),\dots$
is a minimal Cantor set. This means that
each point $f^k(c)$ is the limit of some sequence
$f^{n(k)}(c)$ with $n(k)\to \infty$. Moreover, we shall show that
this Cantor set can be covered in a very natural way.
In the next section this covering shall be used
to show that this Cantor set has `bounded geometry' provided
the critical point of $f$ has order $\ell>2$.

\begin{lemma} The union
\beq
\label{cover}
\bigcup_{i=0}^{S_n-1} f^i(U_n^0) \cup
    \bigcup_{i=0}^{S_{n-1}-1} f^i(U_n^1)
\eeq
is a cover of $\omega(c)$ with mutually disjoint intervals
(the closures of the intervals are disjoint if $n\ge 3$).
Moreover, $\omega(c)$ is a minimal Cantor set.

\end{lemma}
\pr
For $n=2$ statement (\ref{cover}) is easily verified.
Because of 1) in the previous lemma, $f^{S_n}(U_n^0) \subset U_n$ and
$f^{S_{n-1}}(U_n^1) \subset U_n$. But due to 7) and since $\omega(c)$
is forward invariant,
$f^{S_n}(U_n^0 \cap \omega(c))$ and $f^{S_{n-1}}(U_n^1 \cap \omega(c))$
are both contained in $U_n \cap \omega(c) \subset U_n^0 \cup U_n^1$.
This proves the covering property. To show that
the covering consists of disjoint intervals,
mark that $f^i(U_n) \cap U_{n-1} = \emptyset$
for $0 < i < S_{n-1}$.
This is easily verified by similar arguments as in the
previous lemma. In fact, $f^{S_{n-2}}(U_n)$ is adjacent to $U_{n-1}$,
and $f^{S_{n-1}}(U_n)=(d_{n-1},u_{n-1}) \supset U_{n-1}$.
It follows that
$$U_n,f(U_n),\ldots,f^{S_{n-1}-1}(U_n)$$
are mutually disjoint.
The interval $U_n^0$ is symmetric,
so $f(U_n^0) \cap f(U_n^1) = \emptyset$.
Hence 
\[
U_n^0, f(U_n^0),\ldots,f^{S_{n-1}-1}(U_n^0)\text{ and }
U_n^1, f(U_n^1),\ldots,f^{S_{n-1}-1}(U_n^1)
\] 
are all mutually disjoint.
$f^{S_{n-1}}(U_n^0) \subset U_{n-1}^1$ and using induction,
$f^{S_{n-1}+i}(U_n^0) \subset f^i(U_{n-1}^1)$ is disjoint from
$f^i(U_n^0) \supset f^i(U_n^0 \cup U_n^1)$ for $0 \leq i < S_{n-2}$.
This proves that also the intervals $f^i(U_n^0)$, $S_{n-1} \leq i < S_n$,
are mutually disjoint and disjoint from the other intervals.

The fact that $\omega(c)$ is covered by this union implies
that $orb(x) \cap U_n \neq \emptyset$ for every $x \in \omega(c)$
and every $n \geq 2$. So $\omega(c)$ is a minimal Cantor set:
for each $x\in \omega(c)$ one has $\omega(x)\ni c$.
\qed

\bigskip
Next we define a sequence of nested sets $F_n$, each consisting of $2^{n-1}$
intervals, which generates a Cantor set $\cap_n F_n$
such that $\cap_n F_n \supset \omega(c)$.
Let $W$ be the interval containing $c_1$ such that
$f$ maps $W$ diffeomorphically onto $U_2$.
Take
$$F_2^1 = \{U_2^0, U_2^1\}$$
and
$$F_2=F_2^1 \vee ((f|W)^{-1}(F_2^1)). $$
By the previous lemma $F_2^1$ is a covering of $\omega(c) \cap U_2$
and since $\omega(c)\subset U_2\cup W$ this
implies $F_2$ is a covering of $\omega(c)$.
$F_3$ is defined in three steps:
$$F_3^1 = \{U_3^0,U_3^1\},$$
$$F_3^2 = F_3^1 \vee (f^{S_1}|U_2^1)^{-1}(F_2) $$
and
$$F_3=F_3^2\vee ((f|W)^{-1}(F_3^2)).$$
In general, we define $F_n=F_n^{n-1}\vee ((f|W)^{-1}(F_n^{n-1}))$
where
\beqas
    F_n^1 &=& \{U_n^0,U_n^1\}, \\
    F_n^2 &=& F_n^1 \vee (f^{S_{n-2}}|U_{n-1}^1)^{-1}(F_n^1), \\
    ..    &=& ..\\
    F_n^i &=& F_n^{i-1} \vee (f^{S_{n-i}}|U_{n-i+1}^1)^{-1}(F_n^{i-1}), \\
    ..    &=& ..\\
    F_n^{n-1} &=& F_n^{n-2} \vee (f^{S_1}|U_{2}^1)^{-1}(F_n^{n-2}).
\eeqas
Clearly all intervals in $F_n$ are disjoint.

\begin{lemma}\label{covering}
$F_n \supset \omega(c)$ for every $n$. Moreover, $F_n$
consists of $2^n$ components and each component of $F_n$ contains
exactly two components of $F_{n+1}$.
\end{lemma}
\pr
Because $\omega(c) \cap U_{n} \subset U_{n}^0 \cup U_{n}^1$,
in order to prove that $F_n^2$ is a covering of $\omega(c)$
it suffices to prove that
\beqa
\label{impli}
x \in \omega(c) \cap U_{n-1}\text{ implies }f^{S_{n-2}}(x) \in U_{n}.
\eeqa
In fact, (\ref{impli}) also implies inductively that
$F_n^i$ covers $\omega(c)\cap U_{n-i+1}$ (by replacing
in (\ref{impli}) $n$ by $n-i$ it follows that $F_n^i$ covers
$U_{n-i+1}$ if the previous collection $F_n^{i-1}$ already covers
$U_{n-i+2}$.) To prove (\ref{impli}),
note that Lemma~\ref{indi} implies that
$f^{S_{n-1}}(U_n^0)$ is the first return to
$U_{n-1}$, and $f^{S_{n-1}}(U_n^0) \subset U_{n-1}^1$.
Similarly, $f^{S_{n-2}}(U_{n-1}^1)$ is the first return of $U_{n-1}^1$
to $U_{n-1}$. In particular,
$f^{S_{n-2}}(f^{S_{n-1}}(U_n^0))$ is the first return of
$f^{S_{n-1}}(U_n^0)$ to $U_{n-1}$.
But 
\[
f^{S_{n-2}}(f^{S_{n-1}}(U_n^0)) = f^{S_n}(U_n^0) \subset U_n,
\]
because $R_n|U_n^0 = f^{S_n}|U_n^0$.
For the same reason $f^{S_{n-1}}(U_n^1)$ is the first return of $U_n^1$ to
$U_{n-1}$, and $f^{S_{n-1}}(U_n^1) \subset U_n$.
It follows from (\ref{cover}) that $x \in U_{n-1}^1 \cap \omega(c)$
implies $x = f^{S_{n-1}}(y)$ for some $y \in U_n^0 \cap \omega(c)$, and
therefore that $f^{S_{n-2}}(x)=f^{S_n}(y) \subset U_n$.
\qed

\sect{Real bounds for smooth Fibonacci maps}

In this section we shall state and prove
some results on the metric properties of a
smooth Fibonacci map. First we shall quickly state
the main tool needed for these estimates.

\subsection{The cross-ratio tool and the Koebe Principle}
Let $j\subset t$ be intervals and let $l,r$ be the components
of $t\setminus j$. Then
the cross-ratio of this pair of intervals is defined
as
$$C(t,j):=\frac{|t|}{|l|}\frac{|j|}{|r|}.$$
Let $f$ be a smooth function  mapping
$t,l,j,r$ onto $T,L,J,R$ diffeomorphically.
Define 
\[
B(f,t,j)=\frac{|T|\, |J|}{|t|\, |j|}\,\frac{|l|\, |r|}{|L|\, |R|}\
=\frac{C(T,J)}{C(t,j)}
.
\]
It is well known that if the {\it Schwarzian derivative} of $f$, i.e.,
$Sf=f'''/f'-3(f''/f')^2/2$, is negative then
$B(f,t,j)\ge 1$. It is easy to check that our map
$f(z)=z^\ell+c_1$ satisfies $Sf(x)<0$ for $x\in \rz$.

We say that a set $t\subset \rz^k$
contains a {\it $\tau$-scaled} neighbourhood of a disc
$j\subset \rz^k$ with midpoint $x$ and radius $r$
if $t$ contains the ball around $x$ with radius $(1+\tau)r$.

\bigskip
\begin{prop} [Real Koebe Principle]
Let $Sf<0$. Then for any intervals $j\subset t$ and any $n$
for which $f^n|t$ is a diffeomorphism one has
the following.
If $f^n(t)$ contains a $\tau$-scaled
neighbourhood of $f^n(j)$
then
\beq
\label{koebee}
\frac{|Df^n(x)|}{|Df^n(y)|}\le
\left[\frac{1+\tau}{\tau}\right]^2
\eeq
for each $x,y\in j$.
Moreover, there exists a universal function $K(\tau)>0$
which does not depend on $f$, $n$ and $t$
such that
$$|l|,|r|\ge K(\tau)\cdot |j|.$$
\end{prop}

\subsection{The bounds}

Bounds on the relative position of the points $u_n$
and $d_n=c_{S_n}$ are essential in this paper.
They are given in the following theorem.
(All the results in this section
also hold if $f$ is a $C^2$ Fibonacci map
using the disjointness statements as in
\cite{BKNS}.)
\bigskip

\begin{figure}[htp]
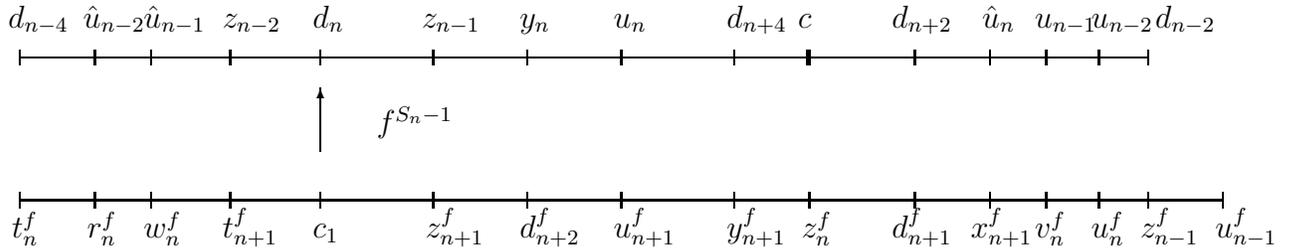

\vskip 0.7cm
\hbox to \hsize{\hss\unitlength=5mm
\beginpic(20,4)(-20,0) \let\ts\textstyle
\put(4,4){\line(-1,0){30}}
\put(4,3.8){\line(0,1){0.4}} \put(4.2,4.8){$d_{n-2}$}
\put(2.7,3.8){\line(0,1){0.4}} \put(2.5,4.8){$u_{n-2}$}
\put(1.3,3.8){\line(0,1){0.4}} \put(1,4.8){$u_{n-1}$}
\put(-0.2,3.8){\line(0,1){0.4}} \put(-0.4,4.8){$\hat u_{n}$}
\put(-2.2,3.8){\line(0,1){0.4}} \put(-2.8,4.8){$d_{n+2}$}
\put(-5,3.8){\line(0,1){0.4}}
\put(-5.05,3.8){\line(0,1){0.4}}
 \put(-5.3,4.8){$c$}
\put(-7,3.8){\line(0,1){0.4}} \put(-7.2,4.8){$d_{n+4}$}
\put(-10,3.8){\line(0,1){0.4}} \put(-10.2,4.8){$u_{n}$}
\put(-12.5,3.8){\line(0,1){0.4}} \put(-12.7,4.8){$y_{n}$}
\put(-15,3.8){\line(0,1){0.4}} \put(-15.3,4.8){$z_{n-1}$}
\put(-18,3.8){\line(0,1){0.4}} \put(-18.2,4.8){$d_{n}$}
\put(-20.4,3.8){\line(0,1){0.4}} \put(-20.6,4.8){$z_{n-2}$}
\put(-22.5,3.8){\line(0,1){0.4}} \put(-22.7,4.8){$\hat u_{n-1}$}
\put(-24,3.8){\line(0,1){0.4}} \put(-24.3,4.8){$\hat u_{n-2}$}
\put(-26,3.8){\line(0,1){0.4}} \put(-26.3,4.8){$d_{n-4}$}

\put(6,0.2){\line(-1,0){32}}
\put(6,0){\line(0,1){0.4}} \put(5.8,-0.8){$u_{n-1}^f$}
\put(4,0){\line(0,1){0.4}} \put(3.8,-0.8){$z_{n-1}^f$}
\put(2.7,0){\line(0,1){0.4}} \put(2.5,-0.8){$u_{n}^f$}
\put(1.3,0){\line(0,1){0.4}} \put(1,-0.8){$v_n^f$}
\put(-0.2,0){\line(0,1){0.4}} \put(-0.7,-0.8){$x_{n+1}^f$}
\put(-2.2,0){\line(0,1){0.4}} \put(-2.8,-0.8){$d_{n+1}^f$}
\put(-5,0){\line(0,1){0.4}} \put(-5.2,-0.8){$z_{n}^f$}
\put(-7,0){\line(0,1){0.4}} \put(-7.2,-0.8){$y_{n+1}^f$}
\put(-10,0){\line(0,1){0.4}} \put(-10.2,-0.8){$u_{n+1}^f$}
\put(-12.5,0){\line(0,1){0.4}} \put(-12.7,-0.8){$d_{n+2}^f$}
\put(-15,0){\line(0,1){0.4}} \put(-15.2,-0.8){$z_{n+1}^f$}
\put(-18,0){\line(0,1){0.4}} \put(-18.2,-0.8){$c_1$}
\put(-20.4,0){\line(0,1){0.4}} \put(-20.6,-0.8){$t_{n+1}^f$}
\put(-22.5,0){\line(0,1){0.4}} \put(-22.7,-0.8){$w_{n}^f$}
\put(-24,0){\line(0,1){0.4}} \put(-24.2,-0.8){$r_{n}^f$}
\put(-26,0){\line(0,1){0.4}} \put(-26.2,-0.8){$t_{n}^f$}
\put(-18,1.5){\vector(0,1){1.7}}
\put(-16.5,2){$f^{S_{n}-1}$ }
\endpic\hss}
\label{mappp}
\vskip 3mm
\caption[ ]{{\small
In the top figure the actual scaling
is completely different for large $\ell$: $|d_{n+2}-c|/|d_{n+4}-c|$ is of order
$1-C\frac{1}{\ell}$ whereas the mutual distance of all points
in the top figure on one  component of $\rz\setminus \{c\}$
is of order $(C/\ell)|d_{n+2}-c|$.
All the distances between the marked points in the bottom figure
(which shows the situation near $c_1$) are of the same order.}}
\end{figure}

\begin{theo}[The real bounds]
\label{realbounds}
There exists $\ell_0\ge 4$ such that
if $f$ is a real unimodal Fibonacci map with a critical
point of order $\ell\ge \ell_0$ with $Sf<0$ then one there exist
universal constants
$0<\lambda<\mu\in (0,1)$
such that the ratio between two consecutive terms
$$|d_{n+1}^f-c_1|<|u_n^f-c_1|< |z_{n-1}^f-c_1|<|d_n^f-c_1|$$
is between $\lambda$ and $\mu$ for all $n$ sufficiently large.
In fact, all the distances in the bottom
part of Figure~\ref{mappp} are of the same order. From this
it follows that the distances near $c$ as stated in the
caption of this figure.
Moreover,
$$\frac{|d_{n-2}^f-c_1|}{|d_n^f-c_1|}\ge 3.85 $$
and therefore
$$\frac{|d_{n-4}^f-c_1|}{|d_n^f-c_1|}\ge 14 $$
for all $n$ sufficiently large.
\end{theo}
\pr The last two inequalities can be found in \cite{fibo}
and also in Lemma 3.3 in \cite{BKNS}.
In Theorem 3.1 of \cite{BKNS} it is shown that
$$\frac{|d^f_n-c_1|}{|u_n^f-c_1|},
\frac{|d^f_n-c_1|}{|d_{n+1}^f-c_1|}\text{ and }
\frac{|u^f_n-c_1|}{|u_{n+1}^f-c_1|}
$$
are bounded and bounded away from one.
Hence there exists uniform constants $C_1,C_2$ such that
$$\frac{C_1}{\ell}
\le
\frac{|d_n-u_n|}{|u_n-c|},
\frac{|d_n-c|-|d_{n+1}-c|}{|d_n-c|},
\frac{|u_n-c|-|u_{n+1}-c|}{|u_n-c|}
\le \frac{C_2}{\ell}$$
for all $n$ large.
From this, by considering the
map drawn in Figure~\ref{mappp} and by
the Koebe Principle one obtains that all distances
are comparable is size. For example,
these inequalities imply that
$[d_{n-2},c]$ is a uniformly scaled neighbourhood of $[u_{n-2},d_{n+2}]$
and by Koebe it follows that
$[z_{n-1}^f,z_n^f]$ is also a scaled neighbourhood of
$[u_n^f,d_{n+1}^f]$. Hence
$$\frac{|z_{n-1}^f-c_1|}{|u_n^f-c_1|}\text{ and }
\frac{|d_{n+1}^f-c_1|}{|z_n^f-c_1|}.$$
are both bounded away from one. Continuing in this way
the proposition follows.
\qed

\bigskip
In fact, we should remark that the last theorem holds for $\ell_0 =4$.
We shall not need this however, and
since the necessary real bounds are only proved in
\cite{BKNS} for $\ell_0$ sufficiently large we only claim the
existence of such an integer $\ell_0$.

We should point out that the previous theorem is false
if $\ell=2$. In that case, $|u_n^f-c_1|/|u_{n+1}^f-c_1|$
goes exponentially fast to infinity, see
\cite{LM} and \cite{fibo}.

Let $T_n=(z_{n-1}^f,t_{n-1}^f)$ be the maximal interval
containing $c_1$ on which $f^{S_n-1}$ is a diffeomorphism
and let $w_n^f\in T_n$ be so that $f^{S_n}(w_n^f)=u_{n-1}^f$.
Then we have the following estimate,
see also Figure~\ref{nsf_qcr2}. This estimate will be needed in
Section~\ref{sectqcr}.

\begin{prop} [Bounds near $c_1$]
\label{43ineqprop}
There exists $\ell_0\ge 4$ such that
if $f$ is a real unimodal Fibonacci map with a critical
point of order $\ell\ge \ell_0$ and $Sf<0$ then
$$\frac{|u_{n-1}^f-c_1|}{|w_n^f-c_1|}\ge \frac{4}{3}$$
for all $n$ sufficiently large.
\end{prop}
\pr
To prove this proposition we use the following lemma.

\begin{lemma}
Let $J'\subset J\subset T$ be intervals
on which $f$ is a diffeomorphism and assume that $Sf<0$.
Then
\beq
\label{growcrr}
B(f,T,J)\ge B(f,T,J').
\eeq
Furthermore, if $f(x)=x^\ell$,
$T=[0,\gamma]$ and $J=[\alpha,\beta]\subset T$
then
$$B(f,T,J)\ge \ell (1-\frac{\alpha}{\gamma}).$$
\end{lemma}
\pr
We may assume that one boundary of $J'$ coincides with one boundary of
$J$ (by applying the lemma twice in this situation
we get the lemma also for general intervals $J$).
Let $L'$ and $R'$ be the components of $T\setminus J'$ which
are labeled so that $R'$ and $R$ both lie on the right hand side of $J$
and $J'$. In order to be definite, assume that the left endpoints of
$J'$ and $J$ coincide. This means that $L'=L$.
It follows that (\ref{growcrr}) is equivalent
to
$$\frac{|f(J)||f(R')|}{|f(R)||f(J)|} \ge
\frac{|J||R'|}{|R||J|}.$$
If we define $\hat T=J\cup R$, $\hat L=J'$, $\hat J=J\setminus J'$
and $\hat R=R'$ then this last inequality becomes
$$
\frac{|f(\hat L\cup \hat J)||f(\hat J \cup \hat R)|}
{|f(\hat L)||f(\hat R)|}
\ge
\frac{|\hat L\cup \hat J||\hat J \cup \hat R|}
{|\hat L||\hat R|}
$$
which is equivalent to the usual cross-ratio expansion:
$$\frac{|f(\hat J)||f(\hat T)|}
{|f(\hat L)||f(\hat R)|}
\ge
\frac{|\hat J||\hat T|}
{|\hat L||\hat R|}.$$ This completes the proof of the first
part of the lemma.

It follows from the first part that we may assume that
$\beta=\alpha$.
Since $f(x)=x^\ell$,
$$B(f,(0,\gamma),\{\alpha\})=
\frac{\gamma^\ell}{\gamma}\cdot \ell\alpha^{\ell-1}
\cdot \frac{\alpha}{\alpha^\ell}\cdot
\frac{\gamma-\alpha}{\gamma^\ell-\alpha^\ell}=
\ell(1-\frac{\alpha}{\gamma})\cdot
\frac{\gamma^\ell}{\gamma^\ell-\alpha^\ell}
\ge \ell(1-\frac{\alpha}{\gamma}).$$
This completes the proof of this lemma.
\qed

\noindent
{\em Proof of Proposition \ref{43ineqprop}:}
Now we can prove the previous proposition.
\beqas
&&
B\left(f^{S_n}, (t_n^f,z_n^f),(c_1,w_n^f)\right)\\
&&\quad =
B\left(f^{S_n-1},(t_n^f,z_n^f),(c_1,w_n^f)\right)
\cdot B\left(f,(d_{n-4},c),(d_n,\hat u_{n-1})\right)
\\
&&\quad \ge 1\cdot  \ell (1-(\frac{|d_n^f-c_1|}{|d_{n-4}^f-c_1|})^{1/\ell})
\ge
\ell ( 1- (\frac{1}{14})^{1/\ell})
\ge 4 ( 1- (\frac{1}{14})^{1/4})> 1.9
\eeqas
where we have used the previous lemma, the inequality
from Theorem~\ref{realbounds} and $\ell\ge 4$.
Now $f^{S_n}(t_n)=d_{n-4}$, $f^{S_n}(z_n)=c_1$,
$f^{S_n}(c_1)=d_n^f$, $f^{S_n}(w_n^f)=u_{n-1}^f$.
Rewriting this last inequality and using the order structure
of the points on the real line,
gives
\beqas
\frac{|u_{n-1}^f-c_1|}{|w_n^f-c_1|}
&\ge & 1.9 \cdot
\frac{|d_{n-4}^f-u_{n-1}^f|}{|d_{n-4}^f-c_1|}
\cdot
\frac{|u_{n-1}^f-c_1|}{|u_{n-1}^f-d_n^f|}
\cdot
\frac{|d_n^f-c_1|}{|z_n^f-c_1|}
\cdot
\frac{|t_n^f-z_n^f|}{|t_n^f-w_n^f|}
\\
&\ge & 1.9 \cdot
\frac{|d_{n-4}^f-d_{n-2}^f|}{|d_{n-4}^f-c_1|}
\cdot
1 \cdot 1 \cdot 1
\ge 1.9 \cdot \left(1-\frac{1}{3.85}\right)\ge \frac{4}{3}.
\eeqas
\qed

The next bounds require that we already know the map
satisfies some renormalization properties,
and is used in Section~\ref{nested} to prove
that certain discs really lie nested.
\medskip

\begin{prop} [Improved bounds near $c_1$ 
if renormalization holds]
\label{boundifreno}
If $\ell\ge \ell_0$, $f$ is as above and
\beq
\lim_{n\to \infty } \,\,\, \frac{|d_n-c|/|d_{n-2}-c|}
{|d_{n-2}-c|/|d_{n-4}-c|}\to 1,
\label{renolimits}
\eeq
then we have the following property.
If $z_{n-1}^f<l_n^f<c_1<s_n^f<t_n^f$ are so that
$$|d_n-c|<|f^{S_n-1}(l_n^f)-c|=|f^{S_n-1}(s_n^f)-c|$$
then
$$\liminf_{n\to \infty}\frac{|l_n^f-c_1|}{|s_n^f-c_1|}\ge 1$$
Moreover, (\ref{renolimits}) implies that
if we take $l_n^f=u_n^f$ and
$r_n^f\in (c_1,t_n^f)\subset T_n$ so that
$f^{S_n}(r_n^f)=\hat u_{n-2}$, then
$|f^{S_n}(r_n^f)-c|=|u_{n-2}-c|=|f^{S_n-1}(u_n^f)-c|$ and
$$\liminf_{n\to \infty} \frac{|u_n^f-c_1|}{|r_n^f-c_1|}>1.$$
\end{prop}
\pr
Consider $f^{S_n-1}$ on $t=(l_n^f,t_n^f)$
and let $j=(c_1,s_n^f)$, $l=(l_n^f,c_1)$
and $r=(s_n^f,t_n^f)$.

\begin{figure}[htp]
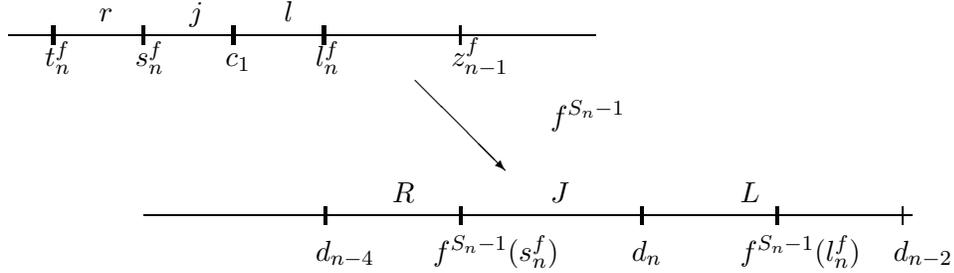

\hbox to \hsize{\hss\unitlength=6mm
\beginpic(20,8)(-20,1) \let\ts\textstyle
\put(-7,7){\line(-1,0){13}}
\put(-10,6.8){\line(0,1){0.4}} \put(-10.2,6.3){\small $z_{n-1}^f$}
\put(-13,6.8){\line(0,1){0.4}} \put(-13.2,6.3){\small $l_n^f$}
\put(-15,6.8){\line(0,1){0.4}} \put(-15.2,6.3){\small $c_1$}
\put(-17,6.8){\line(0,1){0.4}} \put(-17.2,6.3){\small $s_n^f$}
\put(-19,6.8){\line(0,1){0.4}} \put(-19.2,6.3){\small $t_n^f$}
\put(-19.05,6.8){\line(0,1){0.4}}
\put(-17.05,6.8){\line(0,1){0.4}}
\put(-15.05,6.8){\line(0,1){0.4}}
\put(-13.05,6.8){\line(0,1){0.4}}
\put(-13.9,7.3){\small $l$} \put(-16,7.3){\small $j$} \put(-18,7.3){\small $r$}

\put(0,3){\line(-1,0){17}}
\put(-0.2,2.8){\line(0,1){0.4}} \put(-0.4,2){\small $d_{n-2}$}
\put(-3,2.8){\line(0,1){0.4}} \put(-3.8,2){\small $f^{S_n-1}(l_n^f)$}
\put(-6,2.8){\line(0,1){0.4}} \put(-6.2,2){\small $d_n$}
\put(-10,2.8){\line(0,1){0.4}} \put(-10.6,2){\small $f^{S_n-1}(s_n^f)$}
\put(-13,2.8){\line(0,1){0.4}} \put(-13.2,2){\small $d_{n-4}$}
\put(-2.958,2.8){\line(0,1){0.4}}
\put(-5.958,2.8){\line(0,1){0.4}}
\put(-9.958,2.8){\line(0,1){0.4}}
\put(-12.958,2.8){\line(0,1){0.4}}
\put(-3.8,3.3){\small $L$} \put(-8,3.3){\small $J$} \put(-11.5,3.3){\small $R$}
\put(-11,6){\vector(1,-1){2}} \put(-8,5){\small $f^{S_n-1}$}
\endpic\hss}
\caption[ ]{{\small The proof of Proposition~\ref{boundifreno}}}
\end{figure}

\bigskip

Write $a=|f^{S_n-1}(l_n^f)-c|=|f^{S_n-1}(s_n^f)-c|$.
Then $|T|=|d_{n-4}-c|+a$, $|L|=a$, $|J|=a$
and $|R|=|d_{n-4}-c|-a$.
Using the cross-ratio inequality gives
\beqas
\frac{|l_n^f-c_1|}{|r_n^f-c_1|}
&=&
\frac{|j|}{|l|}
\ge
\frac{|L|}{|T|}\frac{|R|}{|J|}
=
\left(\frac{|d_n-c|+a}{|d_{n-4}-c|+a}\right)
\left(\frac{|d_{n-4}-c|-a}{a-|d_n-c|}\right)\\
&\ge&
\left(\frac{|d_n-c|+|d_{n-2}-c|}{|d_{n-4}-c|+|d_{n-2}-c|}\right)
\left(\frac{|d_{n-4}-c|-|d_{n-2}-c|}{|d_{n-2}-c|-|d_n-c|}\right)
\to 1\text{ as }n\to \infty .
\eeqas
Here we have used that the fourth expression is decreasing
in $a\in (0,|d_{n-2}-c|)$ and in the last limit that
(\ref{renolimits}) holds.
To prove the last assertion of the proposition,
note that because of Proposition~\ref{realbounds},
$\limsup_{n\to \infty}\frac{|u_{n-2}-c|}{|d_{n-2}-c|}<1$.
Hence in the second inequality above one has in fact a gain
by a factor which is uniformly strictly larger than one.
\qed

\sect{Background in complex analysis and hyperbolic geometry}

\subsection{Applications of the Schwarz Lemma}
\label{sectionsch}
First we shall review some results from hyperbolic geometry.
Define
$$\cz_J=(\cz\setminus \rz)\cup J$$
where $J\subset \rz$ is an interval.
This set is the complex plane slitted in two infinite rays on the
real line. It is easy to show $\cz_J$ is conformally equivalent
to the upperhalf plane and that
$$D_k(J)=\{z; \text{ the hyperbolic distance to $J$ is at most $k$}\}$$
consists of the intersection with the upper and lower half plane
of two Euclidean discs which are symmetric to each
other with respect to the real line and whose boundaries
intersect the boundary points of $J$, see \cite[pages 485-486]{MS}.
Moreover, $k$ is determined by the external angle $\alpha$ at which the discs
intersect the real line. We also denote this set by
$$D(J;\alpha).$$
 For later use, we define
$D_*(J)$ to be the disc symmetric w.r.t. the real-line and which
intersects the real line in $\partial J$
with angle $\pi/2$.
\vskip0.2cm

\kies{
\begin{figure}[htp]
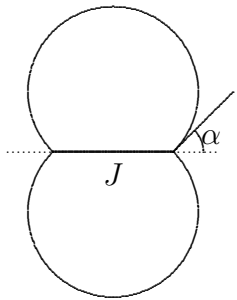
 \hfil \beginpicture
\dimen0=.2cm
\setcoordinatesystem units <\dimen0,\dimen0>
\setplotarea x from -6 to 6, y from -9 to 8
\circulararc 270 degrees from 4 0  center at  0 4
\circulararc 270 degrees from -4 0  center at  0 -4
\setlinear
\plot -4 0 4 0 /
\plot -4 0.05 4 0.05 /
\setdots <2pt>
\plot  -7 0 -4 0 /
\plot   4 0  7 0 /
\setsolid \plot 4 0 8 4 /
\circulararc 45 degrees from 6 0  center at  4 0
\put {$J$} <0mm,-3mm> at 0 0
\put {$\alpha$} <0mm,0mm> at 6.5 1
\endpicture
\caption[ ]{{\small A Poincar\'e neighbourhood of $J$
in $\cz_J$.}}
\end{figure}

}

\begin{lemma} [Schwarz Lemma]
\label{schwarz}

Let $I,J\subset \rz$ be two intervals.
If $G\colon \cz_J\to \cz_I$ is a univalent map
which maps $I$ diffeomorphically onto $J$ then
$G(D_*(J))\subset D_*(I)$.

In particular, let $F\colon \cz\to \cz$ be a real polynomial map
whose critical points are on the real line and
which maps $I$ diffeomorphically onto $J$ then there exists
a set $D\subset D_*(I)$ with $D\cap \rz=I$
which is mapped diffeomorphically by $F$ onto $D_*(J)$.
\end{lemma}
\pr The first statement follows immediately from
the Lemma of Schwarz, which states that any univalent map between
 hyperbolic Riemann surfaces strictly
contracts the Poincar\'e metric.

Since $F$ is a real polynomial and $F$ has no critical values in
$\cz_J$, the inverse $G=F^{-1}\colon \cz_J\to \cz_I$ is a well defined
univalent map. So let $D$ be the inverse of $D_*(J)$ under $G$
and apply the first part of this lemma.
\qed

\subsection{The Koebe Lemma}

As before, we say that a set $t\subset \rz^k$
contains a {\it $\tau$-scaled} neighbourhood of a disc
$j\subset \rz^k$ with midpoint $x$ and radius $r$
if $t$ contains the ball around $x$ with radius $(1+\tau)r$.
(Here we take the standard metric on $\rz^k$.) Then we get
the following classical annalogue of the
real Koebe Principle:

\begin{lemma} [Koebe Lemma]
\label{koebe}

Suppose that $D'\subset \cz$ contains
a $\tau$-scaled neighbourhood of the
disc $D\subset \cz$.
Then for any univalent function $f\colon D'\to \cz$
one has
$$\frac{|f'(x)|}{|f'(y)|}\le \left[\frac{1+\tau}{\tau}\right]^2
\text{ for all }z,y\in D.$$
\end{lemma}
\pr
This result is well known and can be found in for example \cite{ahl1}
and \cite{ahl2} or in \cite{bie}.
\qed

When $J$ is a real interval then take $D(J;\alpha)$
as in the beginning of this section.

\begin{prop}\label{corkoebe}
Assume that $D\subset D'$ and $f\colon D'\to \cz$ are as
in the previous lemma and assume that $f$
maps the real line to the real line.
For each $\alpha\in (\pi/2,\pi)$
there exists $\alpha'\in (\alpha,\pi)$
such that
if $J$ is a real interval in $D$
then $$f(D(J;\alpha))\supset D(f(J);\alpha').$$
(Note that $D(J;\alpha)$
is convex since $\alpha\in (\pi/2,\pi)$.)
\end{prop}
\pr Follows quite easily from the Koebe Lemma.
\qed

\sect{Quasisymmetric rigidity on the real line}
\label{secqsr}
As a preparation for the
proof of Theorem A we shall prove in this section
the following theorem. (This theorem also holds for
$C^2$ maps if we use the disjointess and distortion results
of \cite{BKNS}.)

\medskip
\begin{theo}\label{qsrig}
There exists an integer $\ell_0\ge 4$ with the following property.
Let $f$ be a real unimodal Fibonacci map with $Sf<0$
and with a critical point of order $\ell\ge \ell_0$. Then
$\omega(c)$ has bounded geometry (for the definition see below).
Moreover, there exists
$K<\infty$ and $n_0$ such that for each $n,m\ge n_0$
with $n-m\in 2\zz$,
there exists a quasiconformal homeomorphism $h\colon \cz\to \cz$
which conjugates the first return map of
$R_n\colon \omega(c)\cap U_n\to \omega(c)\cap U_n$
to the first return map of
$R_m\colon \omega(c)\cap U_m\to \omega(c)\cap U_m$.
This map $h$ is symmetric w.r.t. the real line
and $h$ preserves the orientation on the
real line iff $n-m\in 4\zz$.
\end{theo}

\bigskip
Again, using the statement below Theorem~\ref{realbounds},
the above theorem also holds for $\ell_0=4$.
\bigskip
We shall prove this result by constructing a suitable covering
of $\omega(c)\cap U_n$. Firstly, we define
a {\it  presentation}  of a Cantor set $C$ to be a
decreasing collection $F_n \supset F_{n+1}$ of closed sets
such that
\begin{itemize}
\item each $F_n$ is a finite union of closed intervals whose boundary
points are in $C$;
\item each connected component of $F_n$ contains the same
number $a_n$ of connected components of $F_{n+1}$ and
\item $\cap _{n=0}^\infty
F_n = C$.
\end{itemize}
Each component of $F_n$ is called an {\it interval}
of generation  $n$ and each component of $F_n\setminus  F_{n+1}$
is called a {\it  gap} of generation
$n+1$. Of course, there are many presentation of a Cantor set.
\bigskip

 We say that the presentation $\{F_n; n=0,1,2,\dots \} $ of $C$
 has {\it bounded geometry} by $\mu\in (0,1)$
 such that for any interval or gap $I$ of generation
 $n$ and any interval or gap $J \subset I$
 of generation  $n+1$,
       $$ 0< (1-\mu) <\frac {|J|}{|I|} < \mu < 1 .$$
 It follows from the above definition that if the presentation $\{F_n;
 n=0,1,\dots\} $
 of $C$ has bounded geometry then it  has {\it  bounded combinatorics},
  namely, the
 number $a_n$ of components of $F_{n+1}$ in each component of $F_n$ is
 bounded independently of $n$.

 \bigskip
 
 We need the following result.
 \medskip

 \begin{lemma} For each $\mu\in (0,1)$ there exists $K<\infty$
 with the following properties. Let $\{F_n^{(j)}; n=0,1,2,\dots \}$
 be presentations with  geometry bounded by $\mu<1$
 of the Cantor sets $C^{(i)}\subset \rz$, $j=1,2$. Suppose that these
 presentations have the same combinatorics, i.e., the number of
 components
 of $F_{n+1}^{(j)}$ in each component of $F_n^{(j)}$ does not depend on
 $j$.  Then there exists a $K$-quasiconformal homeomorphism $h\colon \cz \to
 \cz $ which is symmetric w.r.t. the real line
 and maps $F_n^{(1)}$ onto $F_n^{(2)}$
 (and therefore $C^{(1)}$ onto $C^{(2)}$.
\end{lemma}
\pr See \cite{MS}[Section VI.3].
\qed

\bigskip

\noindent
{\em Proof of Theorem~\ref{qsrig}:}
From the real bounds in Theorem~\ref{realbounds}
it follows that the size of the
intervals $U_k^0, U_k^1$ and also of the components
of $U_{k-1}\setminus (U_k^0\cup U_k^1)$ are the same up to a
multiplicative constant. Now take $0\le i<k$.
Since the map
$f^{S_{k-i}}\colon U_{k-i+1}^1\to U_{k-i}$ extends
to a diffeomorphism onto $(d_{k-i-2},d_{k-i-4})$;
since (again by Theorem~\ref{realbounds})\,
$(d_{k-i-2},d_{k-i-4})$ contains a uniformly scaled neighbourhood
of $U_{k-i}$, the Koebe Principle
implies that the map $f^{S_{k-i}}\colon U_{k-i+1}^1 \to U_{k-1}$
has uniformly bounded distortion.
Hence the size of the components of $F_k$ and of $F_{k-1}\setminus
F_k$ are all of the same order as the component of $F_{k-1}$
which contains them.
Finally, each component of $F_k$ contains
exactly two components of $F_{k-1}$ and $F_{i+k}\cap U_k$
consists of $2^i$ components.

Therefore and because of the previous lemma, it follows that
there exists a quasiconformal homeomorphism $h$
which is symmetric w.r.t. the real line
which for each $i\ge 0$ sends the $j$-th component of $F_{i+n}\cap U_n$
(say from the left)
to the $j$-th component of $F_{i+m}\cap U_m$ (from the left).
If $n-m\in 4\zz$ then $R_n\colon U_n^0\cup U_n^1\to U_n$
is conjugate to $R_n\colon U_m^0\cup U_m^1\to U_m$
in an orientation preserving way
and so this conjugacy also sends the $j$-th component
of $F_{i+n}$ to the $j$-th component of $F_{i+m}$.
It follows that the homeomorphism $h$ is a conjuagacy from
$U_n\cap \omega(c)$ to $U_m\cap \omega(c)$.
If $n-m\in 2\zz\setminus 4\zz$ then
$R_n\colon U_n^0\cup U_n^1\to U_n$
and $R_n\colon U_m^0\cup U_m^1\to U_m$
are still conjugate but the orientation is reversed;
so this conjugacy sends the $j$-th component
of $F_{i+n}$ to the ($2^i-j$)-th component of $F_{i+m}$.
However, we can also impose that the homeomorphism $h$
from the above lemma reverses orientation (because
$a_n=2$) and then with this choice $h$ again
becomes a conjugacy from
$U_n\cap \omega(c)$ to $U_m\cap \omega(c)$.
\qed

\sect{Quasiconformal rigidity of the return maps; renormalization
and the proof of Theorem A}
\label{sectqcr}

In this section we want to prove Theorem A.
Moreover, we shall prove a renormalization result:
up to rescaling the return maps $R_n\colon U_n^0\cup U_n^1\to U_n$
has at most four limits. This last property will be needed
in the proof of the Main Theorem. In fact, it is needed
in order to apply Proposition~\ref{boundifreno}.
As we mentioned before, we believe that the proof of
the Main Theorem should be independent of this
renormalization result (and of Theorem A); but so far
we have not been able to prove an analogue of
Proposition~\ref{boundifreno} by real methods which is
sufficient for our purposes.

Let $f\colon [0,1]\to [0,1]$ be a
unimodal Fibonacci map.
Let $U_n=[u_{n-1},\hat u_{n-1}]$, $U_n^0=[v_n,\hat v_n]$
and $U_n^1=[u_n,x_n]$ as before.
Moreover, let $R_n\colon U_n^0\cup U_n^1\to U_n$ be
the restriction of the first return map.
In this section we want to show that these first return maps
converge.

\medskip

\begin{theo}[Renormalization result for Fibonacci maps]
\label{reno}
There exists an integer $\ell_0\ge 4$ with the following property.
Fix $\ell\in 2\nz$ with $\ell\ge \ell_0$.
Let $f\colon [0,1]\to [0,1]$ be a
polynomial unimodal Fibonacci map with critical point $c$ of order $\ell$,
let $R_n\colon U_n^0\cup U_n^1\to U_n$ be as above
and let $\phi_n\colon U_n\to [0,1]$ be the affine rescaling map.
Then for each $i\in \{0,1,2,3\}$,
the sequence $\left\{\phi\circ R_{i+4k}\circ \phi^{-1}\right\}_{k\ge 0}$
converges in the $C^1$ topology.
\end{theo}

\bigskip
Again, using the statement below Theorem~\ref{realbounds},
the above theorem also holds for $\ell_0=4$.
Moreover, using the distortion results of \cite{BKNS}
and a shuffling lemma as in Lemma 6.4 of \cite{LM}
the above results also holds for $C^2$ maps.
We should note that in \cite{LM} a similar result is proved
for the case that $\ell=2$ but then $|u_n-c|/|u_{n-1}-c|$
and therefore $|U_n^0|/|U_n|$ tends to zero.
It follows that $|U_{n-1}|/|U_n|\to \infty$ and that
the renormalization map $R|U^0\to U_n$ tends
to a unimodal quadratic map.
In our case (when $\ell>2$)
the situation is more subtle, but on the other hand
in our case of the Cantor set
$\omega(c)$ has bounded geometry (as was shown in the previous section).
This bounded geometry (which does not hold if
$\ell=2$) will help us also a great deal (compare this section with
\cite{L5}).

As we shall also explain in this section,
the above renormalization result
is related to Sullivan's \cite{S2} and McMullen's \cite{McM}
result. These results will also imply
\bigskip

\noindent
{\bf Theorem A}\quad
There exists $\ell_0\ge 4$ with the following property.
{\em For each even $\ell\ge \ell_0$ one has the following properties.
\begin{itemize}
\item
For each $\ell$ there exists a unique parameter $c_1\in \rz$
such that $F(z)=z^\ell+c_1$ has Fibonacci dynamics.
\item There are no measurable invariant linefields on $J(F)$.
\item There exists a sequence of discs $D_n$ and relatively compact
topological discs $D_n^0, D_n^1$ in $D_n$
defined in Proposition~\ref{pl} below,
such that the maps
$$R_n\colon (D_n^0\cup D_n^1)\to D_n$$
defined for $z\in D_n^0\cup D_n^1$
by
$$R_n(z)=\{f^k(z)\st k>0\mbox{ is minimal with }f^k(z)\in D_n\}$$
converge -- up to scaling -- as $n\in 4\nz+i$ tends to infinity,
where $i\in \{0,1,2,3\}$.
\end{itemize}}

\bigskip
\noindent
{\bf Remarks}
\begin{enumerate}
\item In fact, the limits for $i=0$ and $i=2$
in Theorem~\ref{reno} are the same up to orientation and, similarly,
the limits for $i=1$  and $i=3$ are also equal
up to orientation. Similarly, in the last statement of Theorem A.
\item
We do not make any claims about the rate of convergence
in Theorem~\ref{reno}.
\item For the proofs of Theorem~\ref{reno} it would be sufficient to
assume that $f$ is $C^2$. In this case, we proceed
as in \cite[Lemma 6.4]{LM} or as in
\cite[Theorem VI.2.3]{MS} and show that
any limit of $C^2$ maps is an Epstein map.
We shall not discuss this here.
\end{enumerate}
\bigskip

The proof of Theorem~\ref{reno} is very similar to the proof
of D. Sullivan of the convergence of the renormalizations of Feigenbaum
and more general infinitely renormalizable maps, see
\cite{S1} and \cite{S2}. We shall refer to the exposition of these
results given in \cite{MS}. We should note that McMullen has given
an alternative proof of a substantial part of
Sullivan's results, see \cite{McM}. We shall use McMullen's
approach to this result, in order not to have to
develop Sullivan's theory of Riemann surface laminations for
Fibonacci-like maps.
The difference between  Sullivan's case
(of renormalizable maps) and ours (of maps
which are not renormalizable in the classical sense) is that
in the renormalizable
case the return maps have connected Julia sets whereas
in our case the relevant return maps have Julia sets which are
totally disconnected.
\bigskip

To start with the proof of Theorem~\ref{reno}
we first state

\bigskip
\begin{prop}\label{compactness}
There exists $\ell_0\ge 4$ with the following property.
Let $f\colon [0,1]\to [0,1]$ be a
$C^2$ unimodal Fibonacci map with critical point of order $\ell\ge \ell_0$.
Let $R_n\colon U_n\to U_n$ be the
sequence of return maps and $\phi_n\colon U_n\to [0,1]$
be the affine rescaling maps.
Then the closure of $\{\phi\circ R_n\circ \phi^{-1}\}_{n\ge 0}$
forms a compact
family in the $C^1$ topology.
\end{prop}
\pr This follows from the following considerations.
\begin{enumerate}
\item
From the bounds in the previous section, the relative
size of $U_n^0,U_n^1$ and of the components of $U_n\setminus U_n^0,U_n^1$
as subsets of $U_n$ are bounded from above and below. (Note that
the bounds are only claimed to be uniform in $n$ for each fixed $\ell$).
\item The diffeomorphism $R_n|U_n^1\colon U_n^1\to U_n$ can be extended
to a diffeomorphism onto $(d_{n-3},d_{n-5})\supset
(d_{n-2},d_{n-4})$; moreover,
$(d_{n-2},d_{n-4})$ contains  a $\tau(\ell)$-scaled neighbourhood of
$U_n$; in particular, by the Koebe Principle the distortion
of $R_n|U_n^1$ is uniformly bounded and because of 1) the
derivative of $R_n|U_n^1$ is bounded from above and below.
\item the unimodal map $R_n|U_n^0\colon U_n^0 \to U_n$
can be written as a composition of $f$ and a map
from a neighbourhood of $f(U_n^0)$ onto $(d_{n-2},d_{n-4})$,
Therefore, $R_n^0$ is a composition of $f\colon U_n^0\to f(U_n^0)$
and a map whose derivative is bounded from above and below.
\end{enumerate}
All this together implies the proposition.
\qed

Now we will show that the maps $R_n\colon U_n^0\cup U_n^1
\to U_n$ have a polynomial-like extension.
This notion is due to Douady and Hubbard \cite{DH},
see also \cite{MS}, which was extended to be suitable for
the present situation in \cite{LM}. We shall not give the general
definition, but just that of the case we will need.
\bigskip

\noindent
{\bf Definition.}
Let $D^0,D^1,D$ be topological discs bounded
by smooth curves and such that
the closures $D^i$ are disjoint and contained in the interior
of $D$. Then
$$R\colon (D^0\cup D^1) \to D$$
is {\it $\ell$-polynomial-like} if $R|D^1$ is a
univalent map onto $D$ and
if $R|D^0\colon D^0\to D$ is a $\ell$-fold
covering map, i.e., $R|D^0\to D$ is surjective
and the composition of a map of the type $z\mapsto z^\ell$
(up to translation) and a conformal map onto $I$.
The map $R$ is called {\it unbranched}
(using the terminology of McMullen \cite{McM})
if all $R$-iterates of the critical point of $R|D^0$
are contained in $D^0\cup D^1$.

\kies{
\begin{figure}[htp]
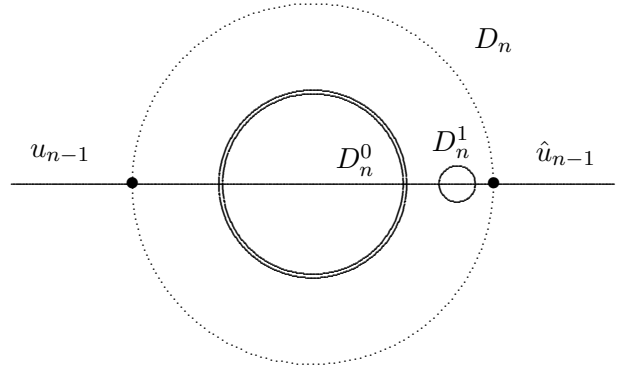
 \hfil
\beginpicture
\dimen0=0.8cm
\setcoordinatesystem units <\dimen0,\dimen0>
\setplotarea x from -5 to 5, y from -3 to 3
\setlinear
\plot -5 0 5 0 /
\multiput {\small $\bullet$} at -3 0 3 0 /
\put {\small $u_{n-1}$} at -4.2 0.5
\put {\small $\hat u_{n-1}$} at 4.2 0.5
\put {\small $D_n^0$} at 0.7 0.4
\put {\small $D_n^1$} at 2.3 0.7
\setdots <2pt>
\circulararc 360 degrees from 3 0  center at  0 0
\setsolid
\circulararc 360 degrees from 1.5 0  center at  0 0
\circulararc 360 degrees from 1.56 0  center at  0 0
\circulararc 360 degrees from  2.1 0  center at  2.4 0
\put{\small $D_n$} <0pt,0pt> at 3 2.4
\endpicture
\caption[ ]{{\small The complex extension of the first return map
$R_n\colon U_n^0\cup U_n^1\to U_n$
is a polynomial-like map
$R_n\colon D_n^0\cup D_n^1\to D_n$ if $\ell=4$.\,\,
The solid disc $D_n^0$ is mapped in a $\ell$-fold way
by $f^{S_n}$ onto the disc $D_n$ with
dotted boundary. The smaller disc
$D_n^1$ is mapped by $f^{S_{n-1}}$ univalently onto $D_n$.
Furthermore, $D_n\cap \rz=U_n^0=(u_{n-1},\hat u_{n-1})$,
$D_n^0\cap \rz=U_n^0=(v_n,\hat v_n)$ and
$D_n^1\cap \rz=U_n^1=(x_n,u_n)$. The domain of the
first return map to $D_n^0$ has infinitely many components
but all iterates of the critical point under the return map
are contained in $U_n^0\cup U_n^1\subset D_n^0\cup D_n^1$.}}
\end{figure}

}

\begin{prop} 
\label{pl}
There exists $\ell_0\ge 4$ with the following property.
Let $f$ be a polynomial Fibonacci map
with a critical point of even order $\ell\ge \ell_0$
and let $R_n\colon U_n^0\cup U_n^1\to U_n$
be the corresponding return maps for $n\ge 2$.
Let $R_n$ also denote the complex extension of this
return map to the disc $D_n=D_*(U_n)$.
Then $R_n$ is polynomial-like for $n$ sufficiently large:
there exists topological discs $D^i_n\subset D_n$
which are symmetric w.r.t. the real line, with
$D^i_n\cap \rz=U_n^i$ such that there exists a complex
extension
$$R_n\colon (D^0_n\cup D_n^1)\to D_n$$
of $R_n\colon U_n^0\cup U_n^1\to U_n$
which is $\ell$-polynomial-like and unbranched.
Moreover, the modulus of the disc $D_n\setminus (D^0_n\cup D^1_n)$
with two holes is bounded from above and below:
for each $\ell\ge \ell_0$ there exist universal constants
$C_1(\ell),C_2(\ell)$ such that for all $n$ sufficiently large,
$$C_1\text{diam}(D_n)
\le
\text{dist}(\partial D_n,\partial D_n^0),
\text{dist}(\partial D_n,\partial D_n^1),
\text{dist}(\partial D_n^0,\partial D_n^1),
\le
C_2\text{diam}(D_n).$$
Moreover, $R_n|D_n^0=f^{S_n}$, $R_n|D_n^1=f^{S_{n-1}}$
and
the distortion
of
$$f^{S_n-1}|D_n^0\text{ and of }f^{S_{n-1}}|D_n^1$$
is uniformly bounded. So the boundary of $D_n^i$ is smooth
and its shape is not too far from `round'.
\end{prop}
\pr Let $W_n$ be the interval containing $c_1$
such that $f^{S_n-1}\colon W_n\to U_n$
is a diffeomorphism. 
Since $R_n$ is a polynomial, the inverse
of
$f^{S_n-1}\colon W_n\to U_n$
extends to an analytic univalent map
$$f^{-(S_n-1)}\colon \cz_{U_n}\to \cz_{W_n}\ni c_1.$$
By the Lemma of Schwarz this univalent map contracts the Poincar\'e metrics
on these spaces, and therefore
$$D_n^{0,f}:=f^{-(S_n-1)}(D_*(U_n))\subset D_*(W_n)\ni c_1.$$
Similarly, let
$W_n'=f(U_n^1)$; this means that $f^{S_{n-1}-1}\colon W_n'\to U_n$
is also a diffeomorphism.
By the same argument, the inverse of this map
has a univalent holomorphic extension and therefore
\beq
\label{incl}
D_n^{1,f}:=f^{-(S_{n-1}-1)}(D_*(U_n))\subset D_*(W_n').
\eeq
So the inverses of the ball $D_n$ are both inside Euclidean
balls $D_*(W_n)$ and $D_*(W_n')$.
Now we have that $W_n=(v_n^f,w_n^f)\ni c_1$ and
$W_n'=(x_n^f,u_n^f)$.
Moreover, by Theorem~\ref{realbounds} these intervals
and the gap between them are of the same order.
Furthermore, it was shown in Proposition~\ref{43ineqprop}, that
\beq
\label{shorter}
|w_n^f-c_1|\le \frac{3}{4} |u_{n-1}^f-c_1|
\eeq
for all $n$ sufficiently large provided $\ell_0$ is sufficiently large.
Since $f(D_*(U_n))$ is equal to a Euclidean disc
centered at $c_1$ and with radius $|c_1-u_{n-1}^f|$
it follows that the closures of
both $D_*(W_n)$ and $D_*(W_n')$ are contained
in the interior of $f(D_*(U_n))$.
In fact, the modulus of the difference set
 -- a disc with two discs taken out -- is bounded and
bounded away from zero (in fact, these bounds
can be taken to be independent of $n$ and $\ell$ because of
the real bounds from Theorem~\ref{realbounds}).
Now take
$$D^0_n:=f^{-1}(D_n^{0,f})=R_n^{-1}(D_*(U_n))\subset D_*(U_n^0)$$
and let $D^1_n$ be the component of
$$f^{-1}(D_n^{1,f})=R_n^{-1}(D_*(U_n))\subset D_*(U_n^1)$$
which contains $(u_n,x_n)$.
So $R_n$ maps $D_n^1$ univalently onto $D_n$
and $D_n^0$ as an $\ell$-cover onto $D_n$.
Then the modulus of the set
$$D_n\setminus (D^0_n\cup D^1_n)$$
is bounded from below and above (in the sense mentioned above).
Since the inverses $f^{-(S_n-1)}$ and $f^{-(S_{n-1}-1)}$
even extend univalently to $\cz_{(d_{n-2},\hat d_{n-2})}$
it follows that the maps $f^{S_n-1}|f(D_n^0)$
and $f^{S_{n-1}-1}|f(D_n^1)$
have uniformly bounded distortion as in the previous proposition.
Since $f\colon D_n^1\to f(D_n^1)$ has bounded distortion
(since $|u_n^f-c_1|/|x_n^f-c_1|$ is bounded)
this implies the last sentence of the proposition.
\qed

\kies{
\begin{figure}[htp]
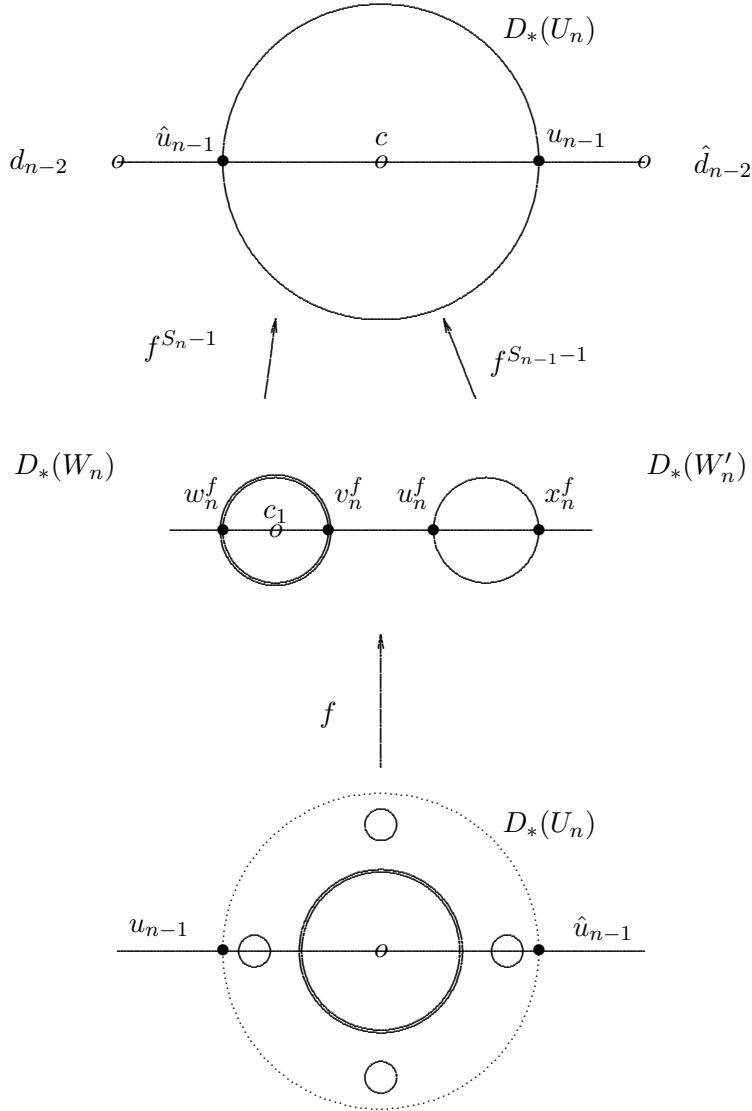
 \hfil
\beginpicture
\dimen0=0.7cm
\setcoordinatesystem units <\dimen0,\dimen0> point at 0 0 
\setplotarea x from -5 to 5, y from -2 to 2
\setlinear
\plot -5 0 5 0 /
\setsolid
\circulararc 360 degrees from 3 0  center at  0 0
\put{\small $c$} <0mm,3mm> at 0 0
\multiput {\small $\bullet$} at -3 0 3 0 /
\multiput {\small $o$} at -5 0 0 0 5 0 /
\put{\small $d_{n-2}$} at -6.5 0
\put{\small $\hat d_{n-2}$} at 6.5 0
\put{\small $\hat u_{n-1}$} <-5mm,3mm> at -3 0
\put{\small $u_{n-1}$} <5mm,3mm> at 3 0
\arrow <5pt> [0.2,0.4] from -2.2 -4.5 to -2 -3
\put{\small $f^{S_n-1}$} at -3.8 -3.5
\arrow <5pt> [0.2,0.4] from 1.8 -4.5 to 1.2 -3
\put{\small $f^{S_{n-1}-1}$} at 3 -3.8
\put{\small $D_*(U_n)$}  <0pt,4pt> at 3.2 2.3
\setcoordinatesystem units <\dimen0,\dimen0> point at 0 7
\setplotarea x from -4 to 4, y from -2 to 2
\setlinear
\plot -4 0 4 0 /
\setsolid
\circulararc 360 degrees from 1 0  center at  2 0
\circulararc 360 degrees from -3 0  center at  -2 0
\circulararc 360 degrees from -3.05 0  center at  -2 0
\multiput {\small $\bullet$} at -3 0 -1 0 1 0 3 0 /
\put {\small $o$} at -2 0
\put {\small $x_n^f$} <8pt,14pt> at 3 0
\put{\small $u_n^f$}  <-8pt,14pt> at 1 0
\put{\small $v_n^f$}  <8pt,14pt> at -1 0
\put{\small $w_n^f$}   <-8pt,14pt> at -3 0
\put{\small $c_1$} <0pt,6pt> at -2 0
\put{\small $D_*(W_n')$}  <0pt,4pt> at 6 1
\put{\small $D_*(W_n)$} <0pt,4pt> at -6 1
\put{\small $f$} at -1 -3.5
\arrow <5pt> [0.2,0.4] from 0 -4.5 to 0 -2
\setcoordinatesystem units <\dimen0,\dimen0> point at 0 15
\setplotarea x from -5 to 5, y from -3 to 3
\put {\small $o$} at 0 0
\setlinear
\plot -5 0 5 0 /
\multiput {\small $\bullet$} at -3 0 3 0 /
\put{\small $u_{n-1}$} <-5mm,3mm> at -3.5 0
\put{\small $\hat u_{n-1}$} <5mm,3mm> at 3.5 0
\setdots <2pt>
\circulararc 360 degrees from 3 0  center at  0 0
\setsolid
\circulararc 360 degrees from 1.5 0  center at  0 0
\circulararc 360 degrees from 1.55 0  center at  0 0
\circulararc 360 degrees from  2.1 0  center at  2.4 0
\circulararc 360 degrees from -2.1 0  center at  -2.4 0
\circulararc 360 degrees from 0 2.1 center at  0 2.4
\circulararc 360 degrees from 0 -2.1   center at  0 -2.4
\put{\small $D_*(U_n)$} <0pt,0pt> at 3.2 2.4
\endpicture
\caption[ ]{\protect{\label{nsf_qcr2}}
{\small The polynomial-like map if $\ell=4$.
The preimages of $D_*(U_n)$ under $f^{S_n-1}$
and $f^{S_{n-1}-1}$ are contained in the discs
$D_*(W_n)$ and $D_*(W_n')$.
To see that the points pull-back on the real line are as shown,
we refer to Figure~\ref{nsf_top4} (and the corresponding figure
if we replace there $n$ by $n-1$).
In the bottom picture the fat circle represents the
component of the preimage of $D_*(U_n)$ under $f^{S_n}$
containing $c$.
The inverse of $D_*(W_n')$ consists of $\ell$ (topological) discs.
Even though the $f$-inverse of $D_*(W_n)$ need not be convex,
it is contained in the disc $D_*(U_n)$ because $|w_n^f-c_1|\le
\frac{3}{4}|u_{n-1}^f-c_1|$ by Proposition~\ref{43ineqprop}.}}
\end{figure}

}

Now we want to show that all the return maps
associated to a polynomial Fibonacci map are quasiconformally conjugate.
For this we shall use the pullback argument of Sullivan,
see \cite{S2} and also Chapter VI of \cite{MS}.
\bigskip

\begin{theo}\label{qce}
There exists $\ell_0\ge 4$ such that for any
unimodal polynomial Fibonacci map $f$
with a critical point of order $\ell\ge \ell_0$,
there exists a constant $K(\ell)<\infty$
with the following properties.
Assume that $R_n\colon D^0_n\cup D^1_n\to D_n$
are the polynomial-like mappings from the previous proposition.
Then for any $n,m$ larger than some sufficiently large
$n_0(\ell)$, these maps $R_n,R_m$ are
$K$-quasiconformally conjugate.
\end{theo}
\pr
For simplicity, let us denote $R_n$, $D_n$, $D_n^0$, $D_n^1$
by $R$, $D$, $D^0$, $D^1$ and  similarly $R_m$, $D_m$, $D_m^0$, $D_m^1$
by $\tilde R$, $\tilde D$, $\tilde D^0$, $\tilde D^1$.
First a warning that we should be careful.
Indeed, as we will show below the filled Julia set of
$R\colon D^0\cup D^1\to D$,
$$K(R)=\{z; R^i(z)\in D^0\cup D^1\text{ for all }i\ge 0\}$$
has positive Lebesgue measure.
(In fact, the filled Julia set is equal to the Julia
set because the critical point of the map is recurrent.)
In particular, the moduli of the set $A_N=D_n\setminus K_N(R_n)$, where
$K_N(R)$ is the filled Julia set, i.e.,
$$K_N(R)=\{z; R^i(z)\in D^0\cup D^1\text{ for }i=0,1,2\dots,N\}$$
will not tend to infinity. So we cannot use the method of
\cite{BH} or rather that of Kahn, see also \cite{L2},
to extend the quasiconformal conjugacy across $K(R)$.

So instead we shall use the pullback argument from
\cite{S1}  following the exposition in \cite{MS}.
The idea of this, is that the quasiconformal conjugacy
between $R_n\colon U_n\cap \omega(c)\to U_n\cap \omega(c)$
to $R_m\colon U_m\cap \omega(c)\to U_m\cap \omega(c)$
from the previous section can be pulled back (stepwise for each $j$)
to a quasiconformal conjugacies between
$R_n\colon U_n\cap f^{-j}(\omega(c))\to U_n\cap f^{-j}(\omega(c))$
to $R_m\colon U_m\cap f^{-j}(\omega(c))\to U_m\cap f^{-j}(\omega(c))$.
Indeed, from Theorem~\ref{qsrig}, provided $\ell\in 2\nz$ is at least
$\ell_0$
the set $\omega(c)\cap U_n$ has a geometry which
is bounded uniformly for $n\ge n_0$.
In particular, there exists by this theorem for each
$\ell\ge \ell_0$
a uniform constant $K<\infty$ such that
for $n,m$ larger than some sufficiently larger $n_0$ there exists
a $K$-quasiconformal homeomorphism $h\colon \cz\to \cz$
which conjugates
$R=R_n\colon \omega(c)\cap D_n\to \omega(c)\cap D_n$
to
$\tilde R=R_m\colon \omega(c)\cap D_m\to \omega(c)\cap D_m$
and which is symmetric w.r.t. the real line
provided $n-m$ is even. If $n-m\in 2\zz\setminus 4\zz$
then this homeomorphism reverses the orientation
on the real line. Now let us change $h$ to a map $h_0$ so that
\begin{itemize}
\item it maps $\partial D$, $\partial D^0$,
$\partial D^0$ diffeomorphically to respectively
$\partial \tilde D$, $\partial \tilde D^0$, $\partial \tilde D^1$;
\item such that it conjugates
$R\colon (\partial D^0 \cup \partial D^1) \to \partial D$
to
$R\colon (\partial \tilde D^0 \cup \partial \tilde D^1)
\to \partial \tilde D$;
\item $h_0$ is symmetric: $h_0(\overline z)=\overline{h_0(z)}$.
\end{itemize}
This means that $h$ is a conjugacy between
the {\it fundamental domains}
$D\setminus (D^0\cup D^1)$ and
$\tilde D\setminus (\tilde D^0\cup \tilde D^1)$).
Since the boundaries of these sets are smooth curves, we can choose
$h$ such that it is $K'$-quasiconformal. The number $K'$ is finite
but can be much larger than $K$ (depending on the shape of
the fundamental domains).
In fact, because of the last sentence in Proposition~\ref{pl},
the shape of the boundaries of $D_n^i$ is `bounded'
and therefore the number $K'$ can be chosen independently of
$n$ and $m$ provided they are sufficiently large.

Now we can define inductively
a sequence of $K'$-quasiconformal $h_i$ such that
\beq
\label{comm}
\tilde R\circ h_{i+1}=h_i\circ R,
\eeq
$$h_{i+1}=h_i\text{ on }\{x; R^j_n(x)\notin D^0\cup D^1\text{ for some
}j=0,1,\dots,i\},$$
$$h_i\text{ is symmetric w.r.t. the real axis},$$
$$h_i\text{ conjugates }
R\text{ and }\tilde R\text{ along the critical orbits}.$$
So assume by induction that
we have already constructed $h_i$.
Since $h_i$ maps the critical value $v$ of $R$
to the critical value $\tilde v$ of $\tilde R$, there exists a unique
lift of $h_i$ to
a map $h_{i+1}\colon \cz\to \cz$ (i.e., such that (\ref{comm}) holds)
which maps $D^i$ onto $\tilde D^i$, which is symmetric w.r.t.
the real axis and for which $h_{i+1}|\rz$ has the same
orientation as $h_i|\rz$. Since $h_i$ is
quasiconformal and the other maps are conformal it follows that $h_{i+1}$ is
quasiconformal with the same conformal distortion as $h_i$. On the other
hand, since $h_i$ coincides with $h_0$ on $D\setminus (D^0\cup D^1)$
and $h_0$
conjugates $R$ with $\tilde R$ on the boundaries of $D^0\cup D^1$ we see that
$h_{i+1}$ coincides with $h_0$ on the boundary of $D^0\cup D^1$.
Hence it can be
extended continuously to $D$ by setting it equal to $h_0$ on $D\setminus
(D^0\cup D^1)$. This extension is quasiconformal and has
the same quasiconformal distortion
as $h_i$ because the boundary of $D^0\cup D^1$ is smooth
(hence has zero Lebesgue measure). Now
we claim that $h_{i+1}$ is a conjugacy from
the critical  orbit of $R$ to
the critical orbit of $\tilde R$. This follows because, by induction,
$h_i$ has this property and because
$R$ has the same combinatorial type as $\tilde R$.
This last statement holds because $n-m\in 2\zz$ and therefore
the $R=R_n, \tilde R=R_m$ are conjugate where the conjugacy is
orientation preserving precisely if $n-m\in 4\zz$.
Therefore (\ref{comm}) and the choice which was
made for the orientation of $h_{i+1}|\rz$ implies that $h_{i+1}$
conjugates $R$ and $\tilde R$ along the critical orbits.

We claim that the sequence $h_i$ converges uniformly to a quasiconformal
homeomorphism $h$ which is a conjugacy between $R$ and $\tilde R$. Indeed,
let $K$ be the quasiconformal distortion of $h_0$. Since all maps $h_i$ are
$K$-quasiconformal, and the set of $K$-quasiconformal
homeomorphisms is compact, we see that there are subsequences that 
converges uniformly. On the other hand, since $h_{i+1}$ is equal to
$h_i$ outside of $R^{-i}(D)$ we see that any two limits of convergent
subsequences must coincide in the complement of the filled
Julia set of $R$. The
claim follows because the filled Julia set of $R$ has empty interior.
(The interior components of the filled Julia set
are bounded components of the Fatou set, and by Sullivan's classification
theorem on wandering domains these are eventually periodic.
The periodic components of the Fatou set
contain iterates of the critical point on their
boundary. This is impossible by the minimality of the
orbit of the critical point.)
\qed
\bigskip

We wish to thank Misha Lyubich for pointing out that the
first two statements of Theorem A can be derived easily from
the Measurable Riemann Mapping Theorem, as
in Section VI.4 of \cite{MS}. In fact, it was he who convinced
us that Theorem A could be useful in this context.
\medskip

\noindent
{\em Proof of the first two statements of Theorem A: }
This is proved exactly as in Theorem 4.2a and Theorem 4.2b
of Chapter VI in \cite{MS}. The reason we can apply
this argument is that because we have a quasiconformal
conjugacy on the critical orbits, see Theorem~\ref{qsrig}.
For the details we refer to Section VI.4 of \cite{MS}, but let us
sketch the idea here. Firstly, the kneading invariant of $z^\ell+c_1$,
depends monotonically on $c_1\in \rz$. So if $[c_1,\tilde c_1]$
is the maximal interval of parameters with this kneading invariant,
then, by the Measurable Riemann Mapping Theorem,
there exists a family $H_u$, $u\in [c_1,\tilde c_1]$
of quasiconformal homeomorphisms
such that $H_u\circ f\circ H_u^{-1}$
is of the form $f(z)=z^\ell+w(u)$ with
$w(c_1)=c_1$ and $w(c_2)=c_2$. Moreover, by the theorem
of Ahlfors and Bers, $u\mapsto w(u)$ is analytic on a neighbourhood of
$[c_1,c_2]\subset \cz$ and so the image of $w$ contains $c_1$
and $c_2$ in its interior. Hence for each $t\in \rz$
near $[c_1,c_2]$, the map $z\mapsto z^\ell+t$ is
also conjugate to $z^\ell+c_1$. This contradicts the maximality of
the interval $[c_1,\tilde c_1]$.

The fact that one has no measurable linefield on the Julia
set of $f$ follows as in the Corollary on page 472 of
\cite{MS}. (The idea is that if an $f$-invariant measurable linefield
on the Julia set then we would obtain a
family of quasiconformal homeomorphisms
$H_u\colon \cz\to \cz$. From the invariance of the line-field,
we get that $H_u\circ f\circ H_u^{-1}$ is also holomorphic
for each $u$.
In fact, one gets that $H_u\circ f\circ H_u^{-1}(z)=z^\ell+w(u)$
where $w$ is a {\it non-constant} analytic function of $u$
(defined on a neighourhood of $c_1$.
This would show that for each $t\in \cz$
near $c_1$, the map $z\mapsto z^\ell+t$ is conjugate
to $z\mapsto z^\ell+c_1$, contradicting the first part
of Theorem A.
\bigskip

To prove the last part of Theorem A we
use an idea of McMullen which is discussed
in \cite{McM1} and \cite{McM2}.  We became aware of these
ideas through informal notes written by Folkert Tangerman.
We would like to thank him and Jacek Graczyk for some
useful discussions on this methods.
\bigskip

\noindent
{\em Proof of the last statement of Theorem of Theorem~\ref{reno}.}
Let us start by emphasizing that we shall fix $\ell$
in the proof of Theorem A. So we do not claim
that the constants (in the proof of this theorem) are independent
of $\ell$.

To clarify the strategy we shall
first define a sequence of renormalizations
of our map related to the Yoccoz puzzle,
and explain why we cannot use this sequence itself.
So let us first explain that one
can associate to a Fibonacci map
a sequence of sets $A_n^0,A_n^1\subset A_n$, $n\ge n_0$,
where $A_n^0,A_n^1$ are two disjoint closed topological discs
which are compactly contained in the interior of $A_n$.
This can be done so that
the return map $R_n\colon (A_n^0\cup A_n^1)\to A_n$ of $f$
$$R_n(x)=\{f^i(x)\st i\mbox{ minimal with }f^i(x)\in A_n\},$$
is equal to $f^{S_n}$ on $A_n^0$
and equal to $f^{S_{n-1}}$ on $A_n^1$.
Indeed, by Theorem~\ref{pl} there exists $n_0$
such that for $n\ge n_0$ and for
$D_n=D^*(u_{n-1},\hat u_{n-1})$, there exists two
discs $D_n^i\subset D_n$ and a polynomial-like map
$$R_n\colon (D^0_n\cup D^1_n)\to D_n$$
such that $R_n|D^0_n=f^{S_n}$ is a covering map
with branch point $c$ and with $R_n(c)\in D^1_n$
and such that $R_n|D^1_n=f^{S_{n-1}}$
is a diffeomorphism with $R_n^2(c)\in D^0_n$.
So fix $n\ge n_0$ and
write $A_0^0=D_n^0$, $A_0^1=D_n^0$, $A_0=D_n$ and $R=R_n$.
Starting with such a polynomial map
\beq
R\colon (A^0_0\cup A^1_1)\to A_0
\label{plms}
\eeq
we can define its renormalization $\R(R)$
as the
polynomial-like map
\beq
\R(R)\colon (A^0_1\cup A^1_1)\to A_1:=A_0^0
\eeq
where $A^1_1,A^2_1\subset A_1$ are defined as follows. Take
$A^1_1$ as the inverse of $A_1\subset D$
under the
$\ell$-fold covering map $(R|D^0)\colon D^0\to D$
and $A^2_1$ as the component containing $R^2(c)$
of the inverse of $A_1\subset D$
under $(R|D^1)\circ (R|D^0)$.

\kies{
\begin{figure}[htp]
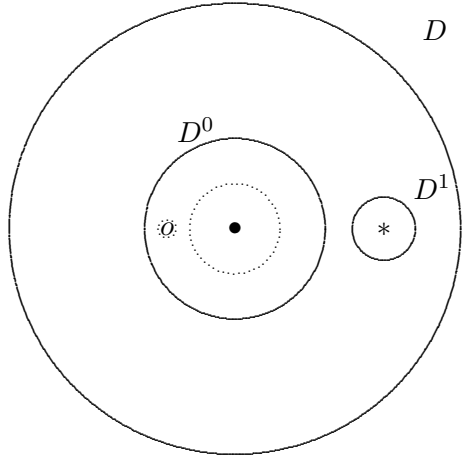
 \hfil
\beginpicture
\dimen0=0.6cm
\setcoordinatesystem units <\dimen0,\dimen0> point at 0 0 
\setplotarea x from -5 to 5, y from -5 to 5
\circulararc 360 degrees from 5 0  center at  0 0
\circulararc 360 degrees from 2 0  center at  0 0
\circulararc 360 degrees from 4 0  center at  3.3 0
\setdots <2pt>
\circulararc 360 degrees from 1 0  center at  0 0
\circulararc 360 degrees from -1.7 0  center at  -1.5 0
\put{\small $D$} <5mm,3mm> at 3.6 3.9
\put{\small $D^0$} <5mm,3mm> at -1.7 1.7
\put{\small $D^1$} <4mm,2mm> at 3.7 0.6
\put{\small $o$} at -1.5 0
\put{\small $*$} at 3.3 0
\put{\small $\bullet$} at 0 0
\endpicture
\caption[ ]{{\small The renormalization of
a polynomial-like mapping. The discs $A\subset D^0$ contains the
critical point $c$ which is denoted by the symbol
$\bullet$. The disc $D^1$ contains
its image $R(c)$ (denoted by $*$). The dics $B\subset D^1$ contains
$R^2(c)$ (denoted by $o$).
The map $R$ sends $A$ first into $D^1$ and
its second iterate is equal to $D^0$. There are $\ell$
topological discs
inside $D^0$ which $R$ sends diffeomorphically onto $D^0$,
but $B$ is the one which contains $R^2(c)$.}}
\end{figure}

}

Since $f$ is the Fibonacci map, this procedure can be
continued infinitely often. So let
$$\R^i(R)\colon (A_i^0\cup A_i^1)\to A_i,$$
$i=0,1,2,\dots$ be the sequence of renormalizations.
Note that in this construction we have
\beq
A_i=A_{i-1}^0.
\label{yopl}
\eeq
Then $\R^i(R)$ is again the first return map
from $A_i^0\cup A_i^1\subset A_i$ onto $A_i$.
This gives a kind of Yoccoz puzzle, associated to
the starting polynomial-like mapping (\ref{plms}).
The trouble, however, is that we are not able
to give lower bounds for the moduli of the
annuli (with two holes) $A_i\setminus (A_i^0\cup A_i^1)$
and also are not sure that the family $\R^i(R)$
is compact. Indeed, even though we have the lower bounds for the moduli
because of Proposition~\ref{pl} for $i=0$,
we do not even know for example that
$A_1=A_0^0$ lies inside the disc
$D_{{n_0}+1}=D^*(u_{n_0},\hat u_{n_0})\subset A_0=D_{n_0}$.
If we knew this, then at least we could apply Proposition~\ref{pl}
again and get inductively that $A_i\subset D_{{n_0}+i}$.
Unfortunately, Proposition~\ref{43ineqprop} only
allows us to conclude that $A_1=A_0^0\subset A_0$.

\bigskip

\noindent
{\em Step 1: The construction of a tower.}
Hence, we shall continue a little differently and drop condition
(\ref{yopl}). We shall do this by constructing a `tower' of polynomial-like
mappings. This tower is a limit of the polynomial-like maps
from Theorem~\ref{pl} and corresponds to an
infinitely blown-up neighbourhood of the origin.
Indeed, $R_n|D_n^0=f^{S_n}$ and there exists a univalent extensions of
$$f^{S_n-1}\colon f(D_n^0)\to D_n$$
and
$$f^{S_{n-1}}\colon D_n^1\to D_n$$
onto the disc $D^*(d_n,d_{n-2})$.
By the real bounds, this last disc is a definite proportion
larger than the disc $D_n$ and so
for each fixed $\ell$ -- up to renormalization --
\beq
R_n\colon (D_n^0\cup D_n^1)\to D_n\mbox {
is in a compact family of maps}.
\label{cofa}
\eeq
(Note that the amount of extendability depends heavily on $\ell$
and so this number $N$ might strongly depend on the choice of $\ell$.)
In particular, by Proposition~\ref{corkoebe},
the set $f(D_n^0)$ contains a disc centered at $c_1$
with diameter $k(\ell)\in (0,1)$ times the 
size of $f(D_n^0)\cap \rz=(w_n^f,v_n^f)$.
So pulling back by $f$ gives that
the topological disc $D_n^0$ also contains at least a
disc with radius $\tilde k(\ell)=[k(\ell)]^{1/\ell}\in (0,1)$
times $|v_n|$.
Hence there exists a number $N(\ell)$ such that
\beq
D_{n}\subset D_{n-N}^0
\label{nedi}
\eeq
for each $n\ge n_0+N$.
Note that
\beq
R_n\mbox{ is an iterate of }R_{n-N}
\label{ric}
\eeq
(restricted to its domain).

By (\ref{cofa}), we can take a sequence $n(j)\to \infty$,
such that -- up to scaling --
$$R_{n(j)}\colon (D_{n(j)}^0\cup D_{n(j)}^1)\to D_{n(j)}$$
converges to a map
$$F_0\colon V_0\to W_0$$
In order to fix matters, let $W_0$ be the unit disc
and let $\Lambda_j$ be the scaling map which maps
the disc $D_{n(j)}$ onto $W_0$. This means that
$$\lim_{j\to \infty} \Lambda_j(D_{n(j)}^0\cup D_{n(j)}^1)\to V_0,$$
where $V_0$ consists of two disjoint topological discs $V_0^0,V_0^1$
which are compactly contained in $W_0$
and
$$\Lambda_j\circ R_{n(j)}\circ \Lambda_j^{-1}$$
converges to $F_0$.
Next take a subsequence of this sequence so that -- up to scaling
by the map $\Lambda_j$ --
$$R_{n(j)-N}\colon (D_{n(j)-N}^0\cup D_{n(j)-N}^1)\to D_{n(j)-N}$$
converges to a map
$$F_1\colon V_1\to W_1.$$
Since, by the real bounds,
there exist $1<\kappa_0<\kappa_1<\infty$ (which do depend on $N$
and so on $\ell$ but not on $j$)
such that
$$1< \kappa_0\le \frac{\mbox{ radius of }D_{n(j)}}
{\mbox{ radius of }D_{n(j)+N}}\le \kappa_1,
$$
the radius of the
disc $W_1$ is in $[\kappa_0,\kappa_1]$.
So taking repeatedly subsequences of subsequences
we get that for each $i\ge 0$,
$$R_{n(j)-i \cdot N}\colon (D_{n(j)-i \cdot N}^0\cup D_{n(j)- i\cdot N}^1)
\to D_{n(j)-i \cdot N}$$
converges -- again up to scaling by $\Lambda_j$ -- as $j\to \infty$ to a map
$$F_i\colon V_i\to W_i.$$
One has
\beq
\kappa_0\le \frac{\mbox{ radius of }W_{i+1}}
{\mbox{ radius of }W_i}\le \kappa_1.
\label{ratra}
\eeq
By Proposition~\ref{pl}, the modulus of $W_i\setminus V_i$
is bounded from above and below. By (\ref{cofa}),
the family of maps
$F_i\colon V_i\to W_i$ are
in a compact set (after rescaling so that the image becomes
a unit disc).
Define the {\em filled Julia set} of $F_i$ as
$$J_i=\{z\in V_i \st F_i^k(z)\in V_i\mbox{ for all }k\ge 0\}$$
and define the {\em post-critical} set of $F_i$ as
$$P_i=\mbox{ closure of the iterates of $c$ under $F_i$}.$$
Note that
$$\cz=\cup_{i\ge 0} W_i= \cup_{i\ge 0}V_i$$
and observe that
$$P=\cup_{i\ge 0}P_i$$
is closed because of (\ref{ric}).
Indeed, this implies that $F_i$ is an iterate
of $F_{i-1}$ and so $P_i\subset P_{i-1}$
but since $F_i$ the first return of $F_{i-1}$ to
$W_i$, one also has that $P_i\cap W_i=P_{i-1}\cap W_i$.
Hence $P\cap W_i$ is equal to $P_i$ which is compact.
Thus $P$ is closed.
Moreover, $J_i$ is the complement of an open dense
set, because it is just a rescaled version of a piece
of the Julia set of the Fibonacci map and since
the Fibonacci set has no periodic attractors,
it is the complement of an open dense set.
Hence by Baire, $J=\cup J_i$ is the complement of a generic set.
(Later on, we shall see that -- in spite of this --
$J$ is dense.)
Define the map
$$F\colon \cz\to \cz$$ on the `tower' $\cup_{i\ge 0}V_i$ as follows:
take $z\in \cz$ and let $i\ge 0$ be {\em minimal} such that $z\in V_i$.
Then define
$$F(z)=F_i(z).$$

\bigskip

\noindent
{\em Step 2: The Poincar\'e metric.}
Consider the Poincar\'e metric $d_P$ on $S=\cz\setminus P$.
Then there exists
a universal constant $C>0$ such that
\beq \frac{1}{C} \cdot \frac{|dz|}{d(z,P)}
\le \rho \le C \cdot \frac{|dz|}{d(z,P)}.
\label{comme}
\eeq
Moreover, the diameter of $W_i\setminus V_i\subset S$
is uniformly bounded.

\noindent
{\em Proof of Step 2:}
Step 2 holds because $\omega(c)$ has bounded geometry,
see Theorem~\ref{qsrig}.
One way to prove Step 2, is to use Theorem 1 of \cite{BP},
see also Theorem 2.3 of \cite{McM}. Let us formulate this
result first. Let $U$ be a hyperbolic region $U$ in $\cz$
(so $\cz\setminus U$ consists of at least three points)
and let $d(z,\partial U)$ be the Euclidean distance
between $z$ to the boundary of $U$.
A round annulus $A=\{z\st r<|z-z_0|<s\}\subset U$ is called
{\em essential} if it is not contractible in $U$;
its modulus is equal to $\log|s/r|$.
The {\em core} curve of $A$ is the circle
$|z-z_0|=\sqrt{rs}$. Next define $mod(z,U)$
as the maximal modulus of a essential round annulus in $U$
and whose core passes through $z$.
Then the Poincar\'e metric $\rho$ on $U$
is comparable to
\beq
\rho'=\frac{|dz|}{d(z,\partial U)(1+mod(z,U))}.
\label{beap}
\eeq
This means that $1/C\le \rho/\rho'\le C$ for some universal
constant $C>0$.

Let us take $U=S=\cz\setminus P$. Since the post-critical
set $P_i$ has bounded geometry, see Theorem~\ref{qsrig},
it follows that there
exists a universal upperbound for the modulus of any
round annulus which is essential with respect to $\cz\setminus P_i$
and which is contained in disc of diameter comparable to $W_i$.
Hence $mod(z,S)$ is bounded from above
and (\ref{comme}) follows from (\ref{beap}).
That the diameter of $W_i\setminus V_i$
(in terms of the Poincar\'e metric on $S$)
is uniformly bounded, follows from
(\ref{comme}) and the bounded geometry.
\bigskip

\noindent
{\em Step 3: $F$ expands the Poincar\'e metric.}
We claim that there exists $\epsilon>0$ and $\kappa>1$ such that
$$F_i\colon (V_i\setminus F_i^{-1}(P))\to (W_i\setminus P)$$
expands the Poincar\'e metric in the sense that if
$x,y\in (V_i\setminus F_i^{-1}(P))$ and
$d_P(x,y)\le \epsilon$ then
$$d_P(F_i(x),F_i(y))\ge d_P(x,y).$$
If, moreover, $F_i(x)\in W_i\setminus V_i$
then
$$d_P(F_i(x),F_i(y))\ge \kappa d_P(x,y)$$
and therefore
$$d_P(F(x),F(y))\ge \kappa d_P(x,y).$$

\noindent
{\em Proof of Step 3:}
Let $S_i=\cz\setminus F_i^{-1}(P)$ and $S=\cz\setminus P$.
Let $d_{P,i}$ and $d_P$ be the Poincar\'e metric on these
sets.
By Schwarz,
$$d_P(x,y)\le d_{P,i}(x,y).$$
Since the inverse of
$$F_i\colon (V_i\setminus F_i^{-1}(P))\to (W_i\setminus P)$$
is a holomorphic covering map
(note that $W_i\cap P=P_i$), we get that $F_i$
is a local isometry in the sense that
$d_P(F_i(x),F_i(y))=d_{P,i}(x,y)$
provided $d_P(x,y)\le \epsilon$ where $\epsilon>0$
is number which is independent of $i$.
Hence
$d_P(F_i(x),F_i(y))\ge d_P(x,y)$
for $d_P(x,y)\le \epsilon$.
This implies the first statement.
To prove the second statement,
note that there exists a constant $C<\infty$ (which is independent of $i$) such that
for any $z\in V_i$ with $F_i(z)\in W_i\setminus V_i$,
\beq
d(z,F^{-1}_i(P))\le C\cdot d(z,P)
\label{cdis}
\eeq
where $d$ is the Euclidean metric on $\cz$.
This holds because there are preimages under
$F_i$ of $P_i=P\cap W_i$ in the annulus $W_i\setminus V_i$
(for example one of the preimages of
the critical point under $F_i$ is in this annulus)
and because of the previous step.
Now (\ref{cdis}) implies that there exists a constant $\kappa>1$
(which is independent of $i$) such that
$$d_P(F_i(z),F_i(w))\ge \kappa \cdot d_P(z,w)$$
when $w$ is sufficiently close to $z$.
This completes the proof of Step 3.
\bigskip

\kies{
\begin{figure}[htp]
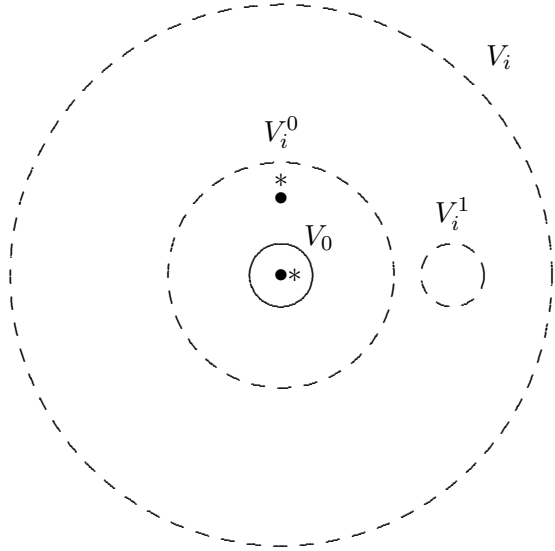
 \hfil
\beginpicture
\dimen0=0.6cm
\setcoordinatesystem units <\dimen0,\dimen0> point at 0 0 
\setplotarea x from -10 to 10, y from -6 to 6
\circulararc 360 degrees from 0.7 0  center at  0 0
\setdashes
\circulararc 360 degrees from 2.5 0  center at  0 0
\circulararc 360 degrees from 4.5 0  center at  3.8 0
\circulararc 360 degrees from 6 0  center at  0 0
\put{\small $\bullet$} at 0 0
\put{\small $*$} at 0.3 0
\put{\small $V_0$} <2mm,2mm> at 0.5 0.5
\put{\small $V_i$} <5mm,5mm> at 4 4
\put{\small $V_i^0$} <0mm,4mm> at 0 2.5
\put{\small $V_i^1$} <0mm,4mm> at 3.8 0.7
\put{\small $*$} at 0 2.1
\put{\small $\bullet$} at 0 1.7
\endpicture

\caption[ ]{\protect{\label{nsf_qcr4}}
{\small The points $c=0,F(c)$ are contained
in the small disc $V_0$; they
are marked with the symbols $\bullet$ and $*$.
The boundaries of the
discs $V_i^0,V_i^1\subset V_i$ are marked with a dashed curve.
Since $F_i$ wraps $V_i^0$ precisely $\ell$ times onto $V_i$
there exists two points $x_1,x_2\in V_i\setminus (V_{i-1}\setminus P)$
which are mapped by $F_i$ to respectively $c$ and $F(c)$.}}
\end{figure}

}

Let $\Delta(x;t)$ be the hyperbolic disc based at $x$
with radius $t$, where we take the Poincar\'e metric on
$S=\cz\setminus P$ from above.
\bigskip

\noindent
{\em Step 4: The set $J$ is uniformly dense in $\cz$.}
We claim that given $\epsilon>0$ there exists $m(\epsilon)$ with
the following properties. For each $x\in W_k$,
$\Delta(x;\epsilon)$ contains $x_0,x_1$ and
$k_0,k_1\le m(\epsilon)$ with
$$F^{k_0}(x_0)=c\mbox{ while }F^{k_0-1}(x_0)\notin P$$
and
$$F^{k_1}(x_1)=F(c)\mbox{ while }F^{k_1-1}(x_1)\notin P.$$
In particular, $\Delta(x;\epsilon)\cap J\ne \emptyset$.

\noindent
{\em Proof of Step 4:}
Since $c=0\in V_0^0$ one has that $F(c)$ is equal to $F_0(c)$
and in $V_0$. First we note that for each $i\ge 1$,
both $c,F(c)\in V_0\subset V_i^0$ have several
preimages under maps $F_i\colon V_i\to W_i$ and we can
choose these preimages outside $P$. (For example just take preimages which
are `near' the imaginary axis, see Figure~\ref{nsf_qcr4}.)
Take two such preimages $x_{i,0}',x_{i,1}'\in V_i\setminus
V_{i-1}$.
This means that
$$F(x_{i,0}')=F_i(x_{i,0}')=c=0\mbox{ and }F(x_{i,1}')=F_i(x_{i,1}')=F(c),$$
while $x_{i,0}',x_{i,1}'\notin P$.
We can choose these points so that
the Euclidean distance of these points to $P$ is of the
same order as the diameter of $W_i$.

So let $n(1)$ be the smallest integer
so that $F^{n(1)}(x)\in W_{k+1}\setminus V_{k+1}$.
If such an integer $n(1)$ does not
exists then $F^{n(1)}(x)\in J_{k+1}$ and so we are done.
Otherwise $F^{n(1)}$ contains a disc of radius $\kappa \epsilon$.
Continuing in this way we get a sequence of integers
$n(1),\dots,n(m)$
so that $F^{n(i)}(x)\in W_{k+i}\setminus V_{k+i}$.
for $i=1,2,\dots,m$ and $F^{n(m)}(B)$
contains a disc of radius $\kappa^m \epsilon$.
But because of the first part of the proof of this claim,
this implies that for $m$ sufficiently large
one has that $F^{n(m)}(B)$ contains a preimage of $c$ and of
$F(c)$ under the map
$F_{k+m}\colon V_{k+m}\to W_{k+m}$.
This concludes the proof of this claim.
\bigskip

\noindent
{\em Step 5: The polynomial-like map has no invariant linefields
on its Julia set.} In Part 2 of Theorem A we proved that
the Julia set of $f$ carries  no measurable $f$-invariant linefield.
Let us show that this implies that 
$$R_i\colon (D_i^0\cup D_i^1)\to D_i$$
also carries no measurable $R_i$-invariant linefield on
its Julia set $J_i$.
So suppose by contradiction that this induced map has such a linefield
$\mu$.
Notice that $R_i$ is a first return map. We will extend
$\mu$ to an $f$-invariant measurable linefield on the
subset 
$\cup_{k\ge 0}f^{-k}(D_i^0\cup D_i^1)$
by $\mu=(f^k)^*\mu$
on $f^{-k}(D_i^0\cup D_i^1)$.
Of course, we have to show that $\mu$ is well defined.
So assume that $x\in 
\cup_{k\ge 0}f^{-k}(D_i^0\cup D_i^1)$
and 
there exists $k'>k$ with 
$$y=f^k(x)\mbox{ and }  y'=f^{k'}(x)=f^{k'-k}(y)\mbox{ are both in}
(D_i^0\cup D_i^1).$$
Since $R_i$ is a first return map,
this implies that $f^{k'-k}$ is an iterate of $R_i$.
Hence, since $\mu$ is $R_i$-invariant one gets that
$(f^{k'-k})^*(\mu(y'))=\mu(y)$. 
This implies that
$$f^k(\mu(y))=f^{k'}(\mu(y'))$$
and hence $\mu(x)$ is well-defined.
Hence, if we define $\mu$ to be zero outside the backward iterates
of $D_i^0\cup D_i^1$ then we get
by construction a linefield which is $f$-invariant.
It is measurable because the original linefield is measurable
and because we have used a countable process to extend its domain.
This implies that the Julia set of $f$ would carry such a linefield,
a contradiction. Therefore $R_i\colon (D_i^0\cup D_i^1)\to D_i$ also
carries no invariant measurable linefield.
Since each of these maps $R_i$ is uniformly quasiconformally
conjugate, these maps are also quasiconformally conjugate to any limit
of $R_i$. Therefore the Julia set of any limit of
$R_i$ also carries no invariant linefield.
\bigskip

\noindent
{\em Step 6: Constructing an invariant linefield
if renormalization does not hold.}
Now we will show how to construct an invariant linefield on
the tower if renormalization does not hold.
Below, we shall show that such an invariant linefield cannot exist,
obtaining a contradiction.

Take $n(j),\tilde n(j),i(j)\to \infty$ with $\tilde n(j)-n(j)\in 4\nz$.
Let $h$ be a quasiconformal conjugacy $h$ between
$$R_{n(j)-i(j)\cdot N}\colon
(D_{n(j)-i(j)\cdot N}^0\cup D_{n(j)-i(j)\cdot N}^1)
\to D_{n(j)-i(j)\cdot N}$$
and
$$R_{\tilde n(j)-i(j)\cdot N}\colon
(D_{\tilde n(j)-i(j)\cdot N}^0\cup D_{\tilde n(j)-i(j)\cdot N}^1)
\to D_{\tilde n(j)-i(j)\cdot N}$$
from Theorem~\ref{qce}.
If the quasi-conformal distortion of $h$ restricted to $D_{n(j)}$
tends to zero as $j\to \infty$, then taking $\Lambda_i\colon \cz\to \cz$
the scaling map which sends $D_i$ onto the unit disc $\Delta$,
$$\Lambda_{\tilde n(j)} \circ h\circ \Lambda_{n(j)}^{-1} \colon \Delta \to \cz$$
tends to a scaling map.
Since $h$ conjugates $$R_{n(j)}\colon (D_{n(j)}^0\cup D_{n(j)}^1)\to D_{n(j)}
\mbox{ to } 
R_{\tilde n(j)}\colon (h(D_{\tilde n(j)}^0)\cup h(D_{\tilde n(j)}^1))\to
h(D_{\tilde n(j)}),$$
it follows that if a subsequence
$R_{n(j)}$ converges to some map $\hat R$ (after rescaling)
then $R_{\tilde n(j)}$ also tends to the same map $\hat R$ after rescaling.
It follows that for each $i_0\in \{0,1,2,3\}$,
the sequence $\{R_{4i+i_0}\}_{i\ge 0}$ (which
is contained in a compact set),
has precisely one limit point (after rescaling).
Hence we are done in this case.

So assume that the quasiconformal distortion of the
conjugacy $h$ does not go to zero.
Then let $\mu$ be the Beltrami-coefficient
of some convergent subsequence of conjugacies
with quasiconformal distortion bounded away form zero
(taking subsequences of the subsequences from above,
so that the sequences of maps from the tower still converge).
Next let $\nu=\pm \mu/|\mu|$ be the corresponding
{\em linefield} defined on the support of $\mu$.
By assumption the support of $\mu$
has positive Lebesgue measure.
Thus we get a measurable linefield $\nu$ on
a set of positive Lebesgue measure in $\cz$
which is invariant under each of
the maps $F_i\colon V_i\to W_i$.
We shall show that this gives a contradiction.

\bigskip

\noindent
{\em Step 7: Constructing a univalent linefield near $c$ and $F(c)$.}

Let us remind the reader that
$z$ is a {\em density point} of a set $E$ if
$$\lim_{t\to 0}\frac{|E\cap B(z;t)|}{|B(z;t)|} = 1$$
where $B(z;t)$ is a discs with centre $z$
and radius $t$ and $|\,\, \cdot \,\,|$ stands
for the Lebesgue measure of a set.
By the Lebesgue Density Theorem, almost every $z$
in $E$ is a density point. Moreover, if $\nu$ on $\cz$
is a measurable function then
almost every $z\in \cz$
is a point of {\em almost continuity} of $\nu$. This means that
for each $\delta>0$ and almost every $z\in \cz$,
$$\lim_{t\to 0} \frac{|\{y\in B(z;t)\st \,\, |\nu(y)-\nu(z)|<\delta\}|}
{|B(x,t)|} = 1.$$
So take a point $z$ which is both a
density point of the
support of the Beltrami coefficient $\nu$
as well as a point of almost continuity of $\nu$.
Since we have already shown that the Julia set of $R_i$
(and also of any of its limits)
does not carry an invariant linefield, see Step 5,
$\nu$ vanishes on $J$ and so we can choose $z\notin J$.
Without loss of generality we can assume that $z\in V_1$.

So define a sequence of integers $k(i)$ so that
$$F^{k(i)-1}(z)\in V_i\mbox{ and }z_i:=F^{k(i)}(z)\in W_i\setminus V_i.$$
Because $z\notin J$, such a sequence $k(i)$ exists.
Now choose $\epsilon>0$ such that there is
a univalent pullback by $F^{k(i)}$ from $\Delta(z_i,2\epsilon)$
to a neighbourhood of $z$. Let $O_i$ be the pullback of
$\Delta(z_i,\epsilon)$ by this map. By Step 3, the diameter
of $O_i$ is exponentially small in terms of $i$.
Since $z$ is a point of almost continuity of $\nu$,
this means that the proportion of the points $y\in O_i$
for which $|\nu(y)-\nu(z)|\ge \delta$ tends to zero.
Define the constant
linefield $\tilde \nu$ on $O_i$ by $\tilde \nu \equiv \nu(z)$.
Next define the linefield
$$\nu_i=(F^{k(i)})_*(\tilde \nu)$$
on $\Delta(z_i,\epsilon)$.
Observe that $F^{k(i)}$ has uniformly bounded distortion
on $O_i$ (this follows by Koebe
since there exists a univalent extension to $\Delta(z_i,2\epsilon)$).
By the invariance of $\nu$
one has
$$\nu=(F^{k(i)})_*(\nu)$$
for all $y\in \Delta(z_i,\epsilon)$
and, combining all this, it follows that for any $\delta>0$,
\beq
\frac
{|\{y\in \Delta(z_i,\epsilon) \st
|\nu(y)-\nu_i(y)|\ge \delta\}|}{|\Delta(z_i,\epsilon)|}\to 0
\label{uni1}
\eeq
as $i\to \infty$.

By Step 4, there exists $x_{i,0}',x_{i,1}'\in \Delta(z_i,\epsilon)$
such that for some $k_0,k_1\le m(\epsilon)$ one has
$$F^{k_0}(x_0)=c\mbox{ while }F^{k_0-1}(x_0)\notin P$$
and
$$F^{k_1}(x_1)=F(c)\mbox{ while }F^{k_1-1}(x)\notin P$$
(in fact, these points even avoid some neighbourhood of $P$).
In particular, this implies that there exists small
neighbourhoods of $x_{i,0},x_{i,1}$ which are mapped
univalently and with bounded distortion to a disc
centered at $c$ respectively $F(c)$ of radius $\rho>0$.
Let $\hat \nu_{i,0}$ be the linefield on
 $\Delta(c;\rho)$ which is defined as the
pushforward by the
univalent maps $F^{k_0}\circ F^{k(i)}$
of the constant linefield $\tilde \nu$
(respectively
$\hat \nu_{i,1}=(F^{k_1}\circ F^{k(i)})_*(\tilde \nu)$
on $\Delta(F(c);\rho)$).
Since $\mu$ is invariant, and since the maps $F^{k_i}$
have bounded distortion, (\ref{uni1}) implies
that the
\beq
\frac
{|\{y\in \Delta(c,\rho) \st
|\nu(z)-\hat \nu_{i,0}(y)|\ge \delta\}|}{|\Delta(c,\rho)|}\to 0
\eeq
as $i\to \infty$ and
similarly for $\nu_{i,1}$ on $\Delta(F(c),\rho)$.
In other words, the restriction of $\mu$ to
$\delta$ discs centered at $c$ and at $F(c)$
is the limit of a sequence of linefields
which are images of a constant linefield under a univalent mapping.
From this it follows that the restriction of
$\mu$ to these discs is actually itself
the image of a constant linefield under a univalent mapping.
(Such linefields are called {\em univalent}.)
\bigskip

\noindent
{\em Step 8: The final contradiction showing that renormalization
holds after all.}

The previous step implies that there are smooth foliations
on a $\delta$ neighbourhoods of $c$ and of $F(c)$
such that the tangent line of the leaves correspond to the linefield
$\mu$ on these neighbourhoods. However, since $\mu$ is invariant
under $F$, the image of the foliation near $c$
must coincide with the foliation near $F(c)$.
This is impossible, because $F$ has a critical point
at $c$. Thus we can conclude that the assumption
we made in Step 6 that renormalization fails, leads
to a contradiction.
\qed

\sect{The random walk argument}
\label{ransec}
In this section we shall state an abstract result
about the evolution of typical points
under a (nearly) Markov map with a kind of random walk structure.
Let $(X,\F,m)$ be some space with probability measure $m$
and $\sigma$-algebra $\F$ and $\A=\{A_k\colon k=0,1,2,\ldots\}$ 
a partition of
$X$ into $\F$-measurable sets. 
$F\colon X\to X$ is a $\F$-measurable transformation, 
$\A_n=\bigvee_{k=0}^{n-1}F^{-k}\A$.
Also assume that there exists $k_0\in \nz$ such that
$$F(A_r)\subseteq\cup_{j=0}^{\infty}A_{r-k_0+j}
\mbox{ for all }r\ge k_0.$$
Observe that $\A$ is a Markov
partition for $F$ if and only if $F^kA$ is an element of 
$\A$ for each $A\in\A_{k+1}$ and each $k\ge 0$.
\par
Define $\ph:X\to\{0,1,2,\ldots\}$ by
\[
\ph(x)=n\mbox{ if }x\in A_n
\]
and
\[
\Delta\ph:=\ph\circ F-\ph\ .
\]
\begin{theo}\label{random}
Assume there are $n_0\in\nz$ and $M>0$ such that for
any $A\in\A_{k+1}$ and any $k\ge 0$ with $\ph_{|F^kA}\ge n_0$
the following inequalities hold:
\beqa\label{drift1}
\int_A(\Delta\ph-1)\circ F^k\, dm &\ge& 0\quad\mbox{ and}
\\ 
\label{secmoments}
\int_A(\Delta\ph)^2\circ F^k\, dm &\le& M\cdot m(A).
\eeqa
Then there exists a set $D\in\F$ with $m(D)>0$ such that
\[
\liminfn\frac{\ph\circ F^j}{j}(x) \ge 1
\quad\mbox{ for each }x\in D,
\]
and
such that for every $x\in D$ the trajectory $x,Fx,F^2x,\ldots$
visits
each set $A_k\in\A$ only finitely often.
\end{theo}
\pr The proof of this theorem is based on a martingale argument
and is due to Gerhard Keller. This proposition and also its proof
can be found as Proposition 4.1 in \cite{BKNS}.
\qed

\sect{A nested sequence of discs}

\begin{figure}[htp]\hfil
\kiesps{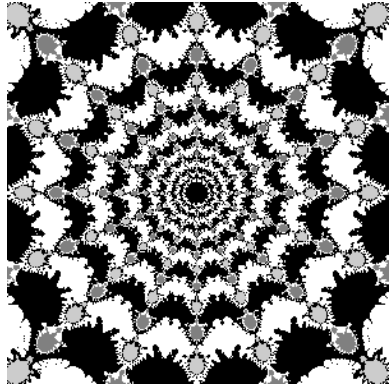} 
\caption{\protect{\small Some levels of the Yoccoz puzzle of the Fibonacci map
$z\mapsto z^\ell+c_1$ when  $\ell=16$. One clearly sees
the long spikes; these spikes are already predicted by
our star-like tips of the set $S_\kappa$. Numerically, the
annuli $D_n\setminus D_{n+1}$ look similar to the
annular regions from the Yoccoz puzzle. This picture was made
by Scott Sutherland in order to simulate a wandering walk
he was studying jointly with Misha Lyubich. The purpose of this
numerical
study was to decide whether the Julia set might have positive Lebesgue measure.
The outcome of this 
turned out to be quite inconclusive. For a discussions
on the reasons for this, see the final section of this paper.}}
\end{figure}

\label{nested}
In this section we shall define an nested sequence
of discs and give some geometric
estimates on these. In the next section we
shall show that these discs can be used to
define an induced mapping with Markov properties.
In fact, these discs are very similar to the
discs constructed in the polynomial-like mapping
from Figure~\ref{nsf_qcr2}. The problem with those discs
is that the inner disc containing $c$ intersects the
real line in the points $v_n,\hat v_n$
and not again points from the sequence
$u_i,\hat u_i$. Therefore we shall take a pull-back of $f^{S_n}$
of a larger disc intersecting the real line in $u_{n-2},\hat u_{n-2}$.
That we can define inductively
a sequence of topological discs 
follows from the real estimate
near $c_1$ which was based on renormalization, see
Proposition~\ref{boundifreno}
and the Lemma of Schwarz.  The resulting partition in annuli is related to
certain annuli from the Yoccoz partition, but
we are not sure whether this Yoccoz partition can be used
directly in our proof. The problem is that we also
need very good estimates on the shape of these annuli.
Although we do not estimate the modulus of these annuli
from below, we do need estimates for their area and we also
need that certain discs are not too `flat'. These estimates
are again based on renormalization by analyzing a
sequence of maps with an almost neutral point. We do not
know whether it is possible to get similar estimates
for the corresponding annuli in the Yoccoz puzzle. 
Therefore we prefer to use our `cruder' partition in annuli.

First we show that one can define a nested
sequence of balls $D_n$ such that
$f^{S_{n+1}}$ maps $D_{n+2}$ as a $\ell$-covering onto
$D_n$. To state the properties of this covering more precisely,
we remind the reader that
$$y_n=f^{S_n}(d_{n+2}),\quad y_n^f=f(y_n)$$
and note that
$|u_n-c|< |y_n-c|< |u_{n-1}-c|$.
Moreover, let
$$a_{n-1}\in \{y_{n-1},\hat y_{n-1}\}
\cap (u_n,\hat u_{n-2}).
$$
Hence $a_{n-1}$ is on the same side of $c$ as $d_n$, $u_n$
and $y_n$. Similarly define
$$b_{n+1}^f\in (c_1,w_n^f)\text{ so that }
a_{n-1}=f^{S_n-1}(b_{n+1}^f).
$$
Then $|u_n-c|< |a_{n-1}-c|=|f^{S_n-1}(b_{n+1}^f)-c| < |u_{n-2}-c|$.
Moreover, let $r_n^f\in (c_1,t_n^f)$ be so that
$$f^{S_{n-1}-1}(r_n^f)=x_{n-1}$$
and therefore so that
$$f^{S_n-1}(r_n^f)=\hat u_{n-2}.$$
(Note that $r_n^f$ is not the image of a point
$r_n\in \rz$; this notation is just to emphasize that
$r_n^f$ is close to $c_1$. Also one should  not
confuse $r_n^f$ with the previously defined point
$w_n^f\in (c_1,t_n^f)$ so that $f^{S_n-1}(w_n^f)=\hat u_{n-1}$.)

\kies{

\begin{figure}[htp]
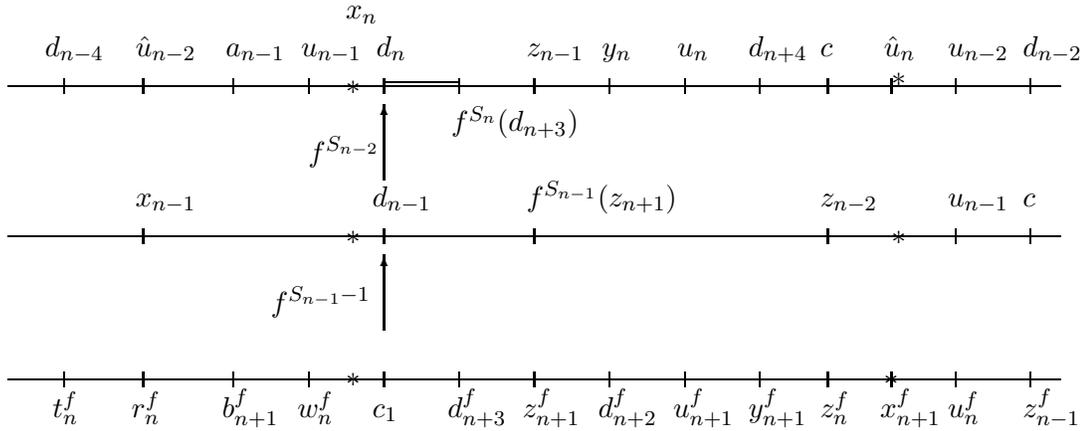

\vskip 1.2cm
\hbox to \hsize{\hss\unitlength=5mm
\beginpic(20,8)(-20,0) \let\ts\textstyle
\put(4,8){\line(-1,0){28}}
\put(3.2,7.8){\line(0,1){0.4}} \put(3,8.8){\small $d_{n-2}$}
\put(1.2,7.8){\line(0,1){0.4}} \put(1,8.8){\small $u_{n-2}$}
\put(-0.5,8){\small $*$}
\put(-0.5,7.8){\line(0,1){0.4}} \put(-0.7,8.8){\small $\hat u_n$}
\put(-2.2,7.8){\line(0,1){0.4}} \put(-2.4,8.8){\small $c$}
\put(-4,7.8){\line(0,1){0.4}} \put(-4.3,8.8){\small $d_{n+4}$}
\put(-6,7.8){\line(0,1){0.4}} \put(-6.2,8.8){\small $u_n$}
\put(-8,7.8){\line(0,1){0.4}} \put(-8.2,8.8){\small $y_n$}
\put(-10,7.8){\line(0,1){0.4}} \put(-10.2,8.8){\small $z_{n-1}$}
\put(-12,7.8){\line(0,1){0.4}}\put(-12.2,6.8){\small $f^{S_n}(d_{n+3})$}
\put(-14,8.1){\line(1,0){2}}
\put(-14,7.8){\line(0,1){0.4}} \put(-14.2,8.8){\small $d_n$}
\put(-15,7.8){\small $*$}\put(-15,9.8){\small $x_n$}
\put(-16,7.8){\line(0,1){0.4}} \put(-16.2,8.8){\small $u_{n-1}$}
\put(-18,7.8){\line(0,1){0.4}} \put(-18.2,8.8){\small $a_{n-1}$}
\put(-20.4,7.8){\line(0,1){0.4}} \put(-20.6,8.8){\small $\hat u_{n-2}$}
\put(-22.5,7.8){\line(0,1){0.4}} \put(-23,8.8){\small $d_{n-4}$}

\put(4,4){\line(-1,0){28}}
\put(3.2,3.8){\line(0,1){0.4}} \put(3,4.8){\small $c$}
\put(1.2,3.8){\line(0,1){0.4}} \put(1,4.8){\small $u_{n-1}$}
\put(-0.5,3.8){\small $*$}
\put(-2.2,3.8){\line(0,1){0.4}} \put(-2.4,4.8){\small $z_{n-2}$}
\put(-10,3.8){\line(0,1){0.4}} \put(-10.2,4.8){\small $f^{S_{n-1}}(z_{n+1})$}
\put(-14,3.8){\line(0,1){0.4}} \put(-14.3,4.8){\small $d_{n-1}$}
\put(-20.4,3.8){\line(0,1){0.4}} \put(-20.6,4.8){\small $x_{n-1}$}
\put(-14,5.5){\vector(0,1){2}}
\put(-15,3.8){\small $*$}
\put(-16,6){\small $f^{S_{n-2}}$}

\put(4,0.2){\line(-1,0){28}}
\put(3.2,0){\line(0,1){0.4}} \put(3,-0.8){\small $z_{n-1}^f$}
\put(1.2,0){\line(0,1){0.4}} \put(1,-0.8){\small $u_{n}^f$}
\put(-0.5,0){\line(0,1){0.4}} \put(-0.8,-0.8){\small $x_{n+1}^f$}
\put(-0.7,0){\small $*$}
\put(-2.2,0){\line(0,1){0.4}} \put(-2.4,-0.8){\small $z_n^f$}
\put(-4,0){\line(0,1){0.4}} \put(-4.3,-0.8){\small $y_{n+1}^f$}
\put(-6,0){\line(0,1){0.4}} \put(-6.3,-0.8){\small $u_{n+1}^f$}
\put(-8,0){\line(0,1){0.4}} \put(-8.3,-0.8){\small $d_{n+2}^f$}
\put(-10,0){\line(0,1){0.4}} \put(-10.3,-0.8){\small $z_{n+1}^f$}
\put(-12,0){\line(0,1){0.4}} \put(-12.3,-0.8){\small $d_{n+3}^f$}
\put(-14,0){\line(0,1){0.4}} \put(-14.3,-0.8){\small $c_1$}
\put(-15,0){\small $*$}
\put(-16,0){\line(0,1){0.4}} \put(-16.3,-0.8){\small $w_n^f$}
\put(-18,0){\line(0,1){0.4}} \put(-18.3,-0.8){\small $b_{n+1}^f$}
\put(-20.4,0){\line(0,1){0.4}} \put(-20.7,-0.8){\small $r_n^f$}
\put(-22.5,0){\line(0,1){0.4}} \put(-22.8,-0.8){\small $t_n^f$}
\put(-14,1.5){\vector(0,1){2}}
\put(-17,2){\small $f^{S_{n-1}-1}$}
\endpic\hss}
\vskip 0.4cm
\caption[ ]{\protect{\label{nsf_nes1}}
{\small Points and their images under $f^{S_{n-1}-1}$ and
$f^{S_{n-2}}$. The slit $Y_{n-1}=[d_n,f^{S_{n}}(d_{n+3})]$ which
will play an important role in the next section, is marked explictly.
We should emphasize that
$D_{n-1}\cap \rz=[\hat u_{n-2},u_{n-2}]$ lives on the top line,
$D_{n-1}^1\cap \rz = [x_{n-1},u_{n-1}]$ on the middle line
and $D_{n+1}\cap \rz=[\hat u_n,u_n]$ on the $f$-preimage of the bottom
line.}}
\end{figure}

}

As before given a bounded real interval $I$, let $D_*(I)$ be the
Euclidean disc
which is symmetric w.r.t. the real line and which intersects
the real axis in $I$. Moreover, if $\kappa\in [0,\pi/2)$
and $y>z>0$ then define $S_\kappa(z,y)^0$ as follows.
Let $l_{\pm \kappa}$ respectively $m_{\pm \kappa}$ be the infinite line through $z>0$
resp. $y$ cutting the real line with angle $\pm \kappa$.
Then define $S_\kappa(z,y)$ to be the closure of the
two components of
$$\{z\in \cz\st |\arg(z)|<2\pi/\ell\}
\,\, \setminus \left( l_{\pm \kappa} \cup m_{\pm} \right)$$
which contain points from $(0,y)$.
Next let $S_\kappa(z,y)^i$ be equal to $S_{\kappa}(z,y)^0$ rotated over
$2\pi i/\ell$ degrees
and let $S_{\kappa}(z,y)=\cup S_{\kappa}(z,y)^i$.

\kies{
\begin{figure}[htp]
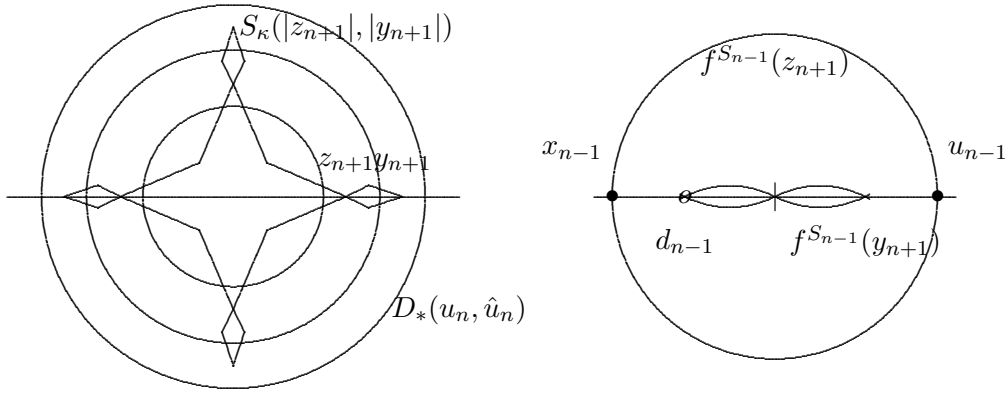
 \hfil \beginpicture
\dimen0=1.5cm \dimen1=0.6cm
\setcoordinatesystem units <\dimen0,\dimen0> point at 0 0
\setplotarea x from -2 to 2, y from -1.8 to 1.8
\plot -2 0 2 0 /
\setlinear
\plot 1 0 0.3 0.3 0 1 -0.3 0.3 -1 0 -0.3 -0.3 0 -1 0.3 -0.3 1 0 /

\plot 1.5 0   1.2 0.1    1 0    1.2 -0.1  1.5 0 /
\plot -1.5 0  -1.2 0.1   -1 0  -1.2 -0.1  -1.5 0 /
\plot 0 1.5   0.1 1.2     0 1   -0.1 1.2   0 1.5 /
\plot 0 -1.5   0.1 -1.2     0 -1   -0.1 -1.2   0 -1.5 /

\setsolid
\circulararc 360 degrees from 1.3 0  center at  0 0
\circulararc 360 degrees from 1.7 0  center at  0 0
\circulararc 360 degrees from 0.8 0  center at  0 0
\put {\small $z_{n+1}$} at 1 0.3
\put {\small $S_{\kappa}(|z_{n+1}|,|y_{n+1}|)$} at 1 1.5
\put {\small $D_*(u_{n},\hat u_{n})$} at 2 -1
\put {\small $y_{n+1}$} at 1.5 0.3

\setcoordinatesystem units <\dimen1,\dimen1> point at -12 0
\setplotarea x from -3 to 3, y from -3 to 3
\plot -4 0 4 0 /
\circulararc 56 degrees from 0 0  center at  -1 -2
\circulararc -56 degrees from 0 0  center at  -1 2
\circulararc -56 degrees from 0 0  center at  1 -2
\circulararc 56 degrees from 0 0  center at  1 2
\circulararc 360 degrees from 3.6 0  center at  0 0
\multiput {\small $\bullet$} at -3.6 0 3.6 0 /
\put {\small $o$} at -2 0
\multiput {\small $|$} at 0 0 /
\put {\small $d_{n-1}$} at -2 -1
\put {\small $x_{n-1}$} at -4.5 1
\put {\small $u_{n-1}$} at 4.5 1
\put {\small $f^{S_{n-1}}(z_{n+1})$} at 0 3
\put {\small $f^{S_{n-1}}(y_{n+1})$} at 2 -1
\endpicture
\caption[ ]{\protect{\label{nsf_nes2}}
{\small The topological ball $D_{n+1}$ is between
$S_{\kappa}(|z_{n+1}|,|y_{n+1}|)$ and the Euclidean
disc $D_*(u_{n},\hat u_{n})$ on the left.
(In fact, the lower bound for $D_{n+1}$ is a pretty
good bound for this set: this set is really squeezed near $z_{n+1}$.
In particular, we should emphasize that we do not have a uniform
lower bound for the moduli of the annuli $A_n$.)
The smaller balls $D_*(u_{n+1},\hat u_{n+1})$
and $D_*(u_{n+2},\hat u_{n+2})$ are also drawn.
If $\ell=4$ then this `star' has $4$ tips (consisting
of `diamonds'). The star does not necessarily
contain the next Euclidean disc $D_*(u_{n+1},\hat u_{n+1})$ completely.
The topological ball $D_{n-1}^1=f^{S_{n-1}}(D_{n+1})$ contains the
interval $[x_{n-1},u_{n-1}]$ and at least the union of
$D((d_{n-1},f^{S_{n-1}}(z_{n+1}));\beta)$ and
$D(f^{S_{n-1}}(z_{n+1}),f^{S_{n-1}}(y_{n+1});\beta)$.
It is inside the disc $D_*(x_{n-1},u_{n-1})$.}}
\end{figure}

}

\begin{theo} \label{complb}
There exist a constant $K<\infty$ and $\ell_0\ge 4$
such that for each $\ell\ge \ell_0$ one has the following properties.
There exists a nested sequence of
open topological balls $D_n$, $D_n^1$ for $n=k_0-2,k_0-1,\dots$
(for some large $k_0$) so that
$D_{k_0-2}$ and $D_{k_0-1}$ are open Euclidean discs,
centered at the critical point
and with boundary through $u_{k_0-3}$ respectively $u_{k_0-2}$
and so that
\begin{enumerate}
\item $D_n\subset D_*(u_{n-1},\hat u_{n-1})$
and $D_n\cap \rz=(u_{n-1},\hat u_{n-1})$;
\item the closure of $D_n$ is inside $D_{n-1}$
and $D_n$ is invariant under a rotation of angle $2\pi/\ell$;
\item $D_n^1\subset D_*(u_n,x_n)$ and $D_n^1\cap \rz=(u_n,x_n)$
(and so this set is in the annulus $A_{n-1}=D_n\setminus D_{n+1}$);
\item $f^{S_{n+1}-1}$ maps $(D_{n+2})^f$ diffeomorphically onto $D_n$;
\item $f^{S_n-1}$ maps $(D_{n+2})^f$ diffeomorphically onto $D_n^1$;
\item $f^{S_{n-1}}$ maps $D_n^1$ diffeomorphically onto $D_n$.
\end{enumerate}
Hence

\bigskip\bigskip

\hbox to \hsize{\hss\unitlength=5mm
\beginpic(20,2)(0,0) \let\ts\textstyle
\put(-5,1){$f^{S_{n}}\colon\,\,$}
\put(-1,1){$D_{n+1}$}
\put(10,1){$D_{n-1}^1$}
\put(19,1){$D_{n-1}$.}
\put(2.5,1){\vector(1,0){6}}\put(5,1.5){$f^{S_{n-1}}$}
\put(11.5,1){\vector(1,0){6}}\put(14,1.5){$f^{S_{n-2}}$}
\endpic\hss}
\vskip 2mm
Moreover, we have the following estimates for the
shape of the topological balls $D_n$ and $D_n^1$, $n\ge k_0$:
there exist universal constants $\kappa>0$ and $\beta>0$
such that
$$S_{\kappa}(|z_n|,|y_n|)\subset
D_n \subset D_*(u_{n-1},\hat u_{n-1})$$
for all $n\ge k_0-2$
and such that
$$\left(D((d_{n-1},f^{S_{n-1}}(z_{n+1}));\beta)\,\, \cup \,\,
D((f^{S_{n-1}}(z_{n+1}),f^{S_{n-1}}(y_{n+1}));\beta)\right)
\subset D_{n-1}^1 \subset D_*(x_{n-1},u_{n-1}).
$$
\end{theo}
\bigskip

Here we remind the
reader that given a real interval $J$ and $\alpha\in (0,\pi)$
we defined a neighbourhood $D(J;\alpha)$ of $J$ in section 4.
This set is a hyperbolic neighbourhood of $J$ in $\cz_J$.
We should also point out that the Properties 1 and 2
stated in this theorem imply that for $n\ge n_0$ the annulus
$A_{n-1}=D_n\setminus D_{n+1}$ intersects
each of the rays  $\rz^+\ni t\mapsto t e^{2\pi i  /\ell}\in \cz$
($i=0,1,\dots, \ell-1$)
in precisely one segments and therefore these rays
divide $A_{n-1}$ into precisely $\ell$ components.

\begin{koro}
For $n$ and $\ell$ sufficiently large,
the Euclidean area of $D_n\setminus D_{n+1}$ is comparable to
the Euclidean area of $D_*(u_{n-1},\hat u_{n-1})\setminus D_*(u_n,\hat u_n)$.
In particular, there exists $\tau>0$ such that
$|D_n\setminus D_{n+1}|> \frac{\tau}{\ell}|D_n|$
for $\ell$ and $n$ large enough.
Similarly, the area of $D_{n-1}^{1,f}$ is at $\tau$ times
the area of $A_{n-2}^f$. The area of $f^{-1}(D_{n-1}^{1,f})$
(consisting of $D_{n-1}$ and all its $\ell$ rotated
versions) is also at least $\tau$ times
the area of $A_{n-2}$.
\end{koro}

\noindent
{\em Proof of the Corollary:}
From the real bounds, see Theorem~\ref{realbounds},
$$\frac{|u_n|-|y_n|}{|y_{n-1}|-|u_n|}$$
is universally bounded from below and above.
It follows that the part of $S_{\kappa}(|z_n|,|y_n|)$
which is outside $D_*(u_n,\hat u_n)$ (i.e., the tips)
has Euclidean area which is of the same order as the
size of the area of
$D_*(u_{n-1},\hat u_{n-1})\setminus D_*(u_n,\hat u_n)$.
The last statement follows since
there are universal constants $C_i$ such that
$$\frac{C_0}{\ell}\le \frac{|u_n-u_{n-1}|}{|u_n-c|}\le \frac{C_1}{\ell}.$$

From the real bounds, the interval $[d_{n-1},f^{S_{n-1}}(z_n)]$
takes up a definite proportion of -- for example --
the interval $[\hat u_{n+1},u_{n-1}]$. Hence, by the
last bound of the previous theorem, the area of the topological
disc $D_{n-1}^1$ is at least a definite proportion of
the area of
$A_{n-2}\cap \{z\in \cz\st \arg(z)<1/\ell\}$. From this the last statement
follows.
\qed
\bigskip

For later references we emphasize that this Corollary implies that
\beq
|D_i\setminus D_{i+i}|=C_i\cdot
\frac{e^{-i/\ell}}{\ell} \cdot |D_n|,
\label{latref}
\eeq
for $i=\ell,\ell+1,\dots$ and where $C_i$ is universally
bounded from below and above for  $n\ge n_0(\ell)$.

\bigskip

\noindent
{\em Proof of Theorem 8.1:}
Take $k_0$ so large that
$$|r_n^f-c_1|<|u_n^f-c_1|$$
for all $n\ge k_0-2$.
By Proposition~\ref{boundifreno}
and Theorem~\ref{reno} such a $k_0$ exists.
(Below we shall increase $k_0$ even further for the last part
of the theorem.)
For $i=k_0-2,k_0-1$ the interval $(u_i,x_i)$ is mapped diffeomorphically
onto $$(u_{i-1},\hat u_{i-1})=D_i\cap \rz$$
by $f^{S_{i-1}}$.
Since $f^{S_{i-1}}$ is a real polynomial, there exists by
Proposition~\ref{schwarz} a set $D_i^1\subset D_*(u_i,x_i)$
with $D_i^1\cap \rz=(u_i,x_i)$
which is mapped by $f^{S_{i-1}}$ diffeomorphically onto $D_i$.
Hence Properties 1-6 are satisfied for $i=k_2-2,k_0-1$.

\kies{
\begin{figure}[htp]
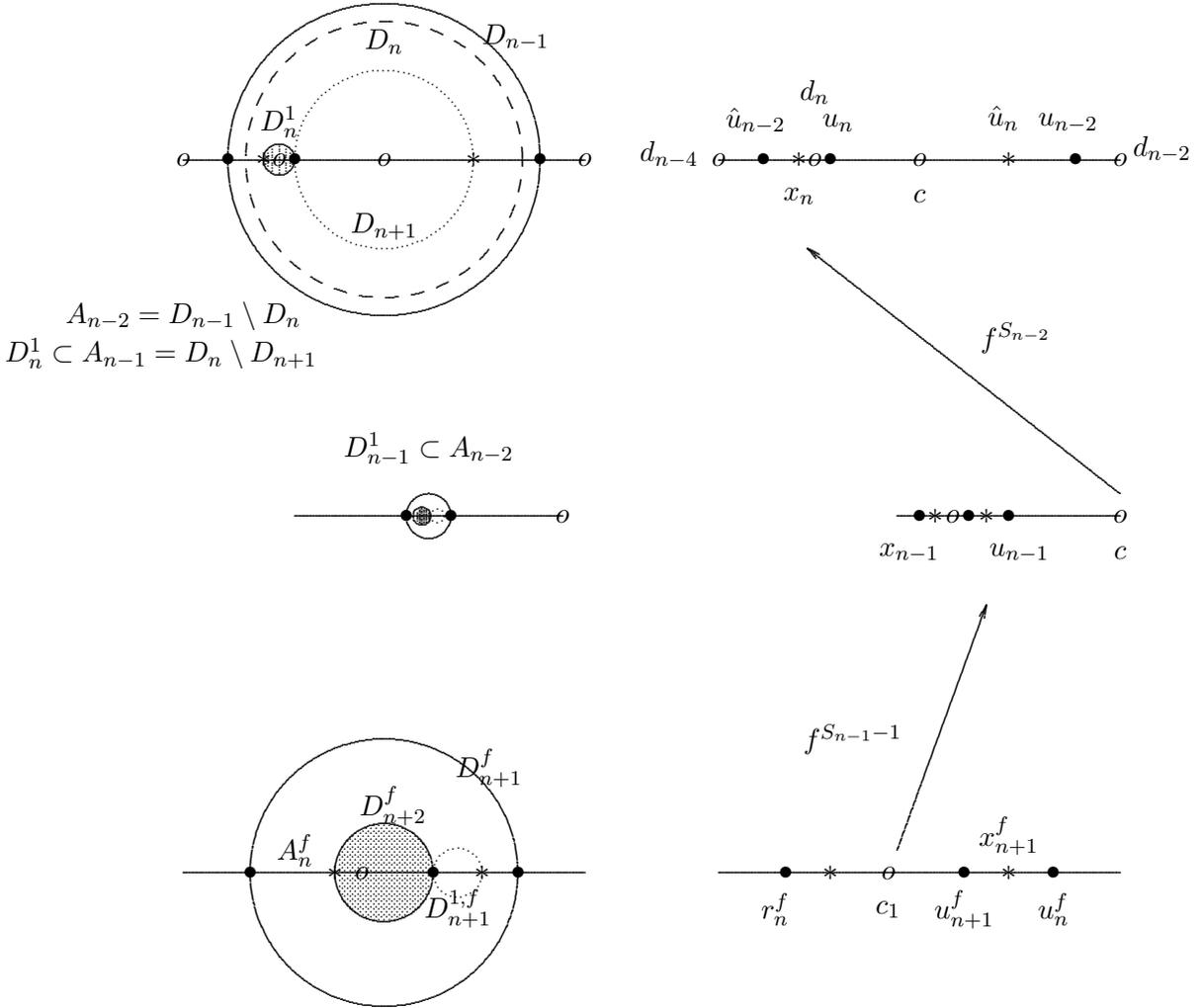
 \hfil
\beginpicture
\dimen0=0.3cm
\setcoordinatesystem units <\dimen0,\dimen0> point at 0 0
\setplotarea x from -9 to 9, y from -7 to 7
\setlinear
\plot -9 0 9 0 /
\put {\small $D_{n-1}$} at 5.9 5.5
\setsolid
\circulararc 360 degrees from 7 0  center at  0 0
\setdashes
\circulararc 360 degrees from 6.2 0  center at  0 0
\setdots <2pt>
\circulararc 360 degrees from 4 0  center at  0 0
\setsolid
\circulararc 360 degrees from -4 0  center at  -4.7 0
\setlinear
\setdots <1pt>
\plot -4.1 -0.2 -4.1 0.2 /
\plot -4.2 -0.3 -4.2 0.3 /
\plot -4.4 -0.5 -4.4 0.5 /
\plot -4.7 -0.7 -4.7 0.7 /
\plot -5   -0.5 -5 0.5   /
\plot -5.2 -0.3 -5.2 0.3 /
\plot -5.3 -0.2 -5.3 0.2 /
\multiput {\small $\bullet$} at -7 0 -4 0 7 0 /
\multiput {\small $o$} at -9 0 0 0 -4.7 0  9 0 /
\multiput {$*$} at -5.4 0 4 0 /
\put {\small $D_n$} at 0 5.3
\put {\small $D_{n+1}$} at 0 -3
\put {\small $D_n^1$} at -4.7 1.7
\put {\small $A_{n-2}=D_{n-1}\setminus D_n$} at -9 -7
\put {\small $D_n^1\subset A_{n-1}=D_{n}\setminus D_{n+1}$} at -10 -8.7
 
\setcoordinatesystem units <\dimen0,\dimen0> point at -24 0
\setplotarea x from -9 to 9, y from -7 to 7
\setsolid
\setlinear
\plot -9 0 9 0 /
\multiput {\small $\bullet$} at -7 0 -4 0 7 0 /
\multiput {\small $o$} at -9 0 0 0 -4.7 0  9 0 /
\multiput {$*$} at -5.4 0 4 0 /
\put {\small $d_{n-2}$} <16pt,4pt> at 9 0
\put {\small $u_{n-2}$}     <-3pt,14pt> at 7 0
\put {\small $\hat u_n$}     <-2pt,16pt> at 4 0
\put {\small $c$} <0pt,-14pt> at 0 0
\put {\small $x_n$} <0pt,-14pt> at -5.4 0
\put {\small $d_n$} <0pt,26pt> at -4.7 0
\put {\small $u_n$}      <3pt,14pt> at -4 0
\put {\small $\hat u_{n-2}$} <-3pt,15pt> at -7 0
\put {\small $d_{n-4}$} <-19pt,2pt> at -9 0
\setcoordinatesystem units <\dimen0,\dimen0> point at 0 16
\setplotarea x from -9 to 9, y from -7 to 7
\setlinear
\plot -4 0 8 0 /
\put {\small $D_{n-1}^1\subset A_{n-2}$} at 2 3
\setsolid
\circulararc 360 degrees from 3 0  center at  2 0
\setdots <2pt>
\circulararc 360 degrees from 2.7 0  center at  2.4 0
\setsolid
\circulararc 360 degrees from 2.1 0  center at  1.7 0
\multiput {\small $\bullet$} at 1 0 3 0  /
\multiput {\small $o$} at  8 0 /
\setdots <0.5pt>
\setlinear
\plot 1.5 -0.2 1.5 0.2 /
\plot 1.6 -0.3 1.6 0.3 /
\plot 1.7 -0.4 1.7 0.4 /
\plot 1.8 -0.3 1.8 0.3 /
\plot 1.9 -0.3 1.9 0.3 /
\plot 2.0 -0.2 2.0 0.2 /
\setcoordinatesystem units <\dimen0,\dimen0> point at -24 16
\setplotarea x from -9 to 9, y from -7 to 7
\setsolid \setlinear
\plot -1 0 9 0 /
\arrow <5pt> [0.2,0.4] from 9 1 to -5 12
\put {\small $f^{S_{n-2}}$} at 4.3 8
\multiput {\small $\bullet$} at 0 0 2.2 0 4 0   /
\multiput {\small $o$} at  1.5 0 9 0  /
\multiput {$*$} at 0.7 0 3 0 /
\put {\small $c$} <0pt,-14pt> at 9 0
\put {\small $u_{n-1}$}  <4pt,-14pt> at 4 0
\put {\small $x_{n-1}$} <-4pt,-14pt> at 0 0
\setcoordinatesystem units <\dimen0,\dimen0> point at 0 32
\setplotarea x from -9 to 9, y from -7 to 7
\setlinear
\plot -9 0 9 0 /
\put {\small $D_{n+1}^f$} at 4.7 4.7
\setsolid
\circulararc 360 degrees from 6 0  center at  0 0
\setsolid
\circulararc 360 degrees from 2.2 0  center at  0 0
\setdots <2pt>
\circulararc 360 degrees from 2.2 0  center at  3.3 0
\multiput {\small $\bullet$} at -6 0 2.2 0 6 0  /
\multiput {\small $o$} at -1 0 /
\multiput {$*$} at -2.2 0  4.4 0 /
\put {\small $A_n^f$} at -4 1
\put {\small $D_{n+2}^f$} at 0.4 3
\put {\small $D_{n+1}^{1,f}$} at 3.3 -1.6
\setquadratic
\setshadegrid span <1pt>
\vshade -2.2 0 0
<,z,,> -2 -0.91 0.91 -1.75 -1.33 1.33
<z,z,,> -1.5 -1.6 1.6 -1 -1.95 1.95
<z,z,,>   0 -2.2 2.2 1 -1.95 1.95
<z,z,,> 1.5 -1.6 1.6 1.75 -1.33 1.33
<z,,,> 2  -0.91 0.91  2.2  0  0 /
\setcoordinatesystem units <\dimen0,\dimen0> point at -24 32
\setplotarea x from -9 to 9, y from -7 to 7
\setsolid \setlinear
\plot -9 0 9 0 /
\arrow <5pt> [0.2,0.4] from -1 1 to 3 12
\put {\small $f^{S_{n-1}-1}$} at -3 6
\multiput {\small $\bullet$} at -6 0 2 0 6 0 /
\multiput {\small $o$} at -1.4 0 /
\multiput {$*$} at -4 0 4 0 /
\put {\small $r_n^f$} <-4pt,-14pt> at -6 0
\put {\small $c_1$} <0pt,-14pt> at -1.4 0
\put {\small $u_{n+1}^f$} <0pt,-14pt> at 2 0
\put {\small $x_{n+1}^f$} <0pt,14pt> at 4 0
\put {\small $u_n^f$}  <0pt,-14pt> at 6 0
\setlinear
\endpicture
\caption[ ]{\protect{\label{nsf_nes3}}
{\small The construction of the new discs $D_{n+1}$.
In the left bottom picture $D_{n+1}^f$ is sketched (in reality
the shape is not round but a smooth version of the star sets
$S_t$ from above). The `dotted' disc
in the annulus $D_{n+1}^f\setminus D_{n+2}^f$ is
denoted $D_{n+1}^f$. The shaded region is mapped into the small shaded
region inside $D_{n-1}^1$ by $f^{S_{n-1}-1}$ and this one is
mapped by $f^{S_{n-2}}$ to $D_n^1$.
On the right the real pullback is drawn -- compare this also
with Figure~\ref{nsf_nes1}, and on the left
the corresponding complex pullback.
If we pullback the bottom topological disc by $f$ then we obtain
the disc $D_{n+1}$ from the top part of the figure. That this
disc is inside $D_n$ follows from Proposition~\ref{boundifreno}.
Note that $D_{n+1}^f\cap \rz=[u_n^f,r_n^f]$, $D_{n+1}^{1,f}\cap \rz=
[x_{n+1}^f,u_{n+1}^f]$. We should emphasize that because of the
real bounds, the real parts of the topological balls
$D_{n+2}^f$, $D_{n+1}^{1,f}$ and $D_{n+1}^f$ have more or less equal length.
Since these balls are mapped to respectively $D_n^1$, $D_{n+1}$ and
$D_{n-1}$ (and the first one of these is small compared to the other
two) this shows that the conformal map $f^{S_n-1}\colon D_{n+1}^f\to D_{n-1}$
is very different from a Moebius transformation. We shall use
this fact in a crucial way at the end of the paper to get
estimates which are much better than those which would follow from Koebe.
Moreover, $D_{n-1}^1$ is in fact far from a real ball: it is squeezed
at $z_{n-2}$ as in Figure~\ref{nsf_nes2}.}}
\end{figure}

}

So assume that $D_i$ and $D_i^1$
satisfying Properties 1-6 already are defined by induction
for $i=k_2-2,k_0-1,\dots,n$.
Since all iterates of the critical point of $f$ are in the
real line, one has that
\beq
\label{lowest}
c_{S_k}\in D_k\text{ and }c_i\notin D_k\text{ for }0<i<S_k.
\eeq
Now we define $D_{n+1}$ as follows. Let $D_{n-1}^1$
be the topological ball which is already defined and
which is mapped diffeomorphically by $f^{S_{n-2}}$ onto $D_{n-1}$.
One has $D_{n-1}^1\cap \rz=(u_{n-1},x_{n-1})$ and $D_{n-1}^1\subset
D_*(u_{n-1},x_{n-1})$.
Because of the results in Section~\ref{seccomb},
$f^{S_{n-1}-1}$ maps $(u_n^f,r_n^f)$ diffeomorphically
onto $(u_{n-1},x_{n-1})$. It follows by Proposition~\ref{schwarz}
that there exists a set $D_{n+1}^f\subset D_*(u_n^f,r_n^f)$
with $D_{n+1}^f\cap \rz=(u_n^f,r_n^f)$ which is mapped
diffeomorphically onto $D_{n-1}^1\subset D_*(u_{n-1},x_{n-1})$
by $f^{S_{n-1}-1}$.

Since $f^{S_{n-1}-1}$ maps $D_n^f\ni c_1$ diffeomorphically onto
$D_{n-2}$ and the set $D_{n+1}^f$ into $D_{n-1}^1$
and because the closure of $D_{n-1}^1$ is contained in
the closure of $D_{n-1}\subset D_{n-2}$,
it follows that
$D_{n+1}^f$ is contained in the interior of $D_n^f$.
Hence
$D_{n+1}$ is contained in the interior of $D_n$.
Since
$D_{n+1}^f\subset D_*(u_n^f,r_n^f)\ni c_1$
and since we have by Proposition~\ref{boundifreno},
$$|r_n^f-c_1|<|u_n^f-c_1|$$
(as $n\ge k_0$)
we even get
$D_{n+1}\subset D_*(u_n,\hat u_n)$.
Similarly, $f^{S_n}$ maps $(u_{n+1},x_{n+1})$ diffeomorphically
to $(u_n,\hat u_n)$. Hence by Proposition~\ref{schwarz}
there exists a set $D_{n+1}^1\subset D_*(u_{n+1},x_{n+1})$
with $D_{n+1}^1\cap \rz=(u_{n+1},x_{n+1})$
which is mapped diffeomorphically onto $D_{n+1}\subset D_*(u_n,\hat u_n)$.
This proves Properties 1-6.

\bigskip

Before proving the last statement
we state a proposition. This proposition will imply the proof of
Theorem~\ref{complb} and the remainder of this section
will be dedicated to the proof of the proposition.

It is convenient to choose $z_{n+1}$ to be on the same side
of $c$ as $a_{n+1}$, $d_{n+2}$, $u_{n+2}$ and $y_{n+2}$.
Then $[z_{n+1},a_{n+1}]$ contains $d_{n+2}$
and this interval does not contain $c$.

In the next proposition we shall show that there exists 
a region $V_{n,\ell}$ which looks like one of  the tips
from $S_\kappa(|z_{n+1}|,|y_{n+1}|)$, see Figure~\ref{nsf_nes2},
which is mapped by $f^{S_n}$ into a similar tip associated
to the next set $S_\kappa(|z_{n-1}|,|y_{n-1}|)$.
\bigskip 
\begin{prop}
For each sufficiently large $\ell$
there exists $n_0(\ell)$ and for each
$n\ge n_0(\ell)$
a closed topological disk $V_{n,\ell}$ containing
$[z_{n+1},a_{n+1}]\ni d_{n+2}$ such that
\begin{itemize}
\item $f^{S_n}$ maps $V_{n,\ell}$ into the interior of $V_{n-2,\ell}$
for each $n\ge n_0(\ell)+1$;
\item there exists a universal number $\alpha\in (\pi/2,\pi)$ 
such that 
$V_{n,\ell}\supset D([z_{n+1},a_{n+1}];\alpha)$ for each $n\ge n_0(\ell)-1$;
\item $V_{n,\ell}\subset D_*(u_{n},\hat u_{n})$ for $n\ge n_0(\ell)-1$.
\end{itemize}
\end{prop}

\bigskip

Before proving this proposition let us show that
it allows us to complete the proof of Theorem 8.1.
That is, we shall show that $D_{n+1}$ contains
$V_{n,\ell}$ and therefore one of the tips of
$S_\kappa(|z_{n+1}|,|y_{n+1}|)$.
An additional argument will then
show that $D_{n+1}$ also contains the main piece
of $S_\kappa(|z_{n+1}|,|y_{n+1}|)$ and also the
other $\ell-1$ tips.

\bigskip
  
\noindent
{\em Conclusion of the proof of Theorem 8.1:}
Let us first prove by induction that
$D_{n+1}$ contains $V_{n,\ell}$ for $n$ sufficiently large.
To do this, choose $n_0$ to be
equal to the integer $k_0$ from the previous theorem.
For $n=n_0-2,n_0-1$ the sets $D_n$ are the Euclidean discs
$D_*(u_{n-1},\hat u_{n-1})$
which by the last property of the proposition
therefore contain $V_{n+1,\ell}$. Hence the inductive statement holds
for $n_0-2,n_0-1$ by the last property stated in the proposition.
Now assume that $D_{n-1}\supset V_{n-2,\ell}$ for some $n\ge n_0$.
Since by definition $f^{S_n}$ maps $D_{n+1}$ onto $D_{n-1}$,
the second assertion of the proposition implies that
$D_{n+1}\supset V_{n,\ell}$, proving the inductive step.
Thus we obtain
by induction that $D_{n+1}\supset V_{n,\ell}$ for all $n\ge n_0-2$.
In other words, $D_{n+1}$ at least contains one of the
small tips of $S_\kappa(|z_{n+1}|,|a_{n+1}|)$.

Since $f^{S_n}$ maps $D_{n+1}$ onto $D_{n-1}\supset V_{n-2,\ell}$
we can get a better lower bound for $D_{n+1}$.
In fact, we want to show that $D_{n+1}$ also contains
the `big' starshaped piece of $S_\kappa(|z_{n+1}|,|a_{n+1}|)$.
Indeed, consider the neighbourhood $V'=D([a_{n-1},z_{n-1}];\alpha)\ni d_n$
of $[a_{n-2},z_{n-1}]$. 
Since $f^{S_n-1}$ is diffeomorphism from some interval
neighbourhood of $c_1$ to $[d_{n-2},d_{n-4}]$
the corresponding inverse branch of $f^{S_n-1}$
on $V'$ has uniformly bounded distortion. It follows from
Lemma~\ref{corkoebe} that there exists a universal number $\alpha'\in (\pi/2,\pi)$
such that the corresponding component of $f^{-(S_n-1)}(V')$
contains $D([c_1,z_{n+1}^f];\alpha')$
(in fact, it contains $V''=D([b_{n+1}^f,z_{n+1}];\alpha')$.)
Hence $D_{n+1}$ contains $f^{-1}(V'')$.
Note that the inverse of $V''$ under $f$
contains the `kite' component of
\beq
\{z\in \cz\st |\arg(z)|<2\pi/\ell\}
\,\, \setminus \left( l_{\pm \kappa}  \right)
\label{kite}
\eeq
containing $(0,z_{k+1})$.
Here $l_{\pm \kappa}$ is the infinite line through $|z_{k+1}|$
with angle $\pm \kappa$ and where $\kappa>0$ is some universal number.
Moreover, as we have proved above,
$D_{n+1}$ contains $V_{n,\ell}$
and since $D_{n+1}$ is invariant under rotation under
$2\pi/\ell$ degrees, $D_{n+1}$ also contains
the rotated versions of these sets and also of the kites from
(\ref{kite}).
Combining this, shows that $D_{n+1}$ contains
a set of the form
$S_{\kappa}(|z_{n+1}|,|y_{n+1}|)$.

Hence $D_{n+1}^f$ contains at least a set of the form
$$D((c_1,z_{n+1}^f);\beta')\cup D((z_{n+1}^f,y_{n+1}^f);\beta')$$
where $\beta>0$ is some universal number.
Since $D_{n-1}^1$ is inside $D_*(x_{n-1},u_{n-1})$
and the map $f^{S_{n-1}-1}$ has uniformly bounded
distortion on $D_{n+1}^f$ (one has uniform Koebe space around
$D_*(x_{n-1},u_{n-1})$), it follows that 
$D_{n-1}^1$ contains
$$\left(D((d_{n-1},f^{S_{n-1}}(z_{n+1}));\beta)\,\, \cup \,\,
D((f^{S_{n-1}}(z_{n+1}),f^{S_{n-1}}(y_{n+1}));\beta)\right).$$
\qed
\bigskip

\kies{
\begin{figure}[htp]
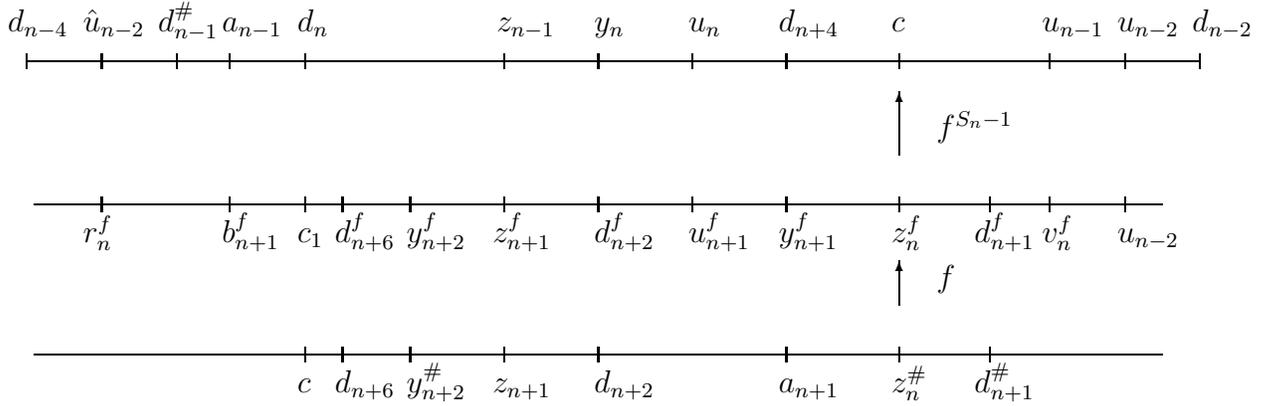

\vskip 0.2cm
\hbox to \hsize{\hss\unitlength=5mm
\beginpic(20,10)(0,-5) \let\ts\textstyle
\put(-5,0.2){\line(1,0){30}}
\put(-3.2,0){\line(0,1){0.4}} \put(-3.7,-0.8){$r_n^f$}
\put(0.2,0){\line(0,1){0.4}} \put(0,-0.8){$b_{n+1}^f$}
\put(2.2,0){\line(0,1){0.4}} \put(2,-0.8){$c_1$} 
\put(3.2,0){\line(0,1){0.4}} \put(3,-0.8){$d_{n+6}^f$}
\put(5,0){\line(0,1){0.4}} \put(4.9,-0.8){$y_{n+2}^f$} 
\put(7.5,0){\line(0,1){0.4}} \put(7.2,-0.8){$z_{n+1}^f$}  
\put(10,0){\line(0,1){0.4}} \put(9.9,-0.8){$d_{n+2}^f$}  
\put(12.5,0){\line(0,1){0.4}} \put(12.4,-0.8){$u_{n+1}^f$}
\put(15,0){\line(0,1){0.4}} \put(14.8,-0.8){$y_{n+1}^f$}
\put(18,0){\line(0,1){0.4}} \put(17.8,-0.8){$z_{n}^f$}
\put(20.4,0){\line(0,1){0.4}} \put(20,-0.8){$d_{n+1}^f$}
\put(22,0){\line(0,1){0.4}} \put(21.8,-0.8){$v_n^f$}
\put(24,0){\line(0,1){0.4}} \put(23.8,-0.8){$u_{n-2}$}

\put(-5.2,4){\line(1,0){31.2}}
\put(-5.2,3.8){\line(0,1){0.4}} \put(-5.7,4.8){$d_{n-4}$}
\put(-3.2,3.8){\line(0,1){0.4}} \put(-3.7,4.8){$\hat u_{n-2}$}
\put(-1.2,3.8){\line(0,1){0.4}} \put(-1.7,4.8){$d_{n-1}^\#$}
\put(0.2,3.8){\line(0,1){0.4}} \put(0,4.8){$a_{n-1}$}
\put(2.2,3.8){\line(0,1){0.4}} \put(2,4.8){$d_n$}
\put(7.5,3.8){\line(0,1){0.4}} \put(7.3,4.8){$z_{n-1}$}
\put(10,3.8){\line(0,1){0.4}} \put(9.9,4.8){$y_n$}
\put(12.5,3.8){\line(0,1){0.4}} \put(12.4,4.8){$u_n$}
\put(15,3.8){\line(0,1){0.4}} \put(14.8,4.8){$d_{n+4}$}
\put(18,3.8){\line(0,1){0.4}} \put(17.8,4.8){$c$}
\put(22,3.8){\line(0,1){0.4}} \put(21.8,4.8){$u_{n-1}$}
\put(24,3.8){\line(0,1){0.4}} \put(23.8,4.8){$u_{n-2}$}
\put(26,3.8){\line(0,1){0.4}} \put(25.8,4.8){$d_{n-2}$}
\put(18,1.5){\vector(0,1){1.7}}
\put(19,2){$f^{S_n-1}$ }

\put(-5,-3.8){\line(1,0){30}}
\put(2.2,-4){\line(0,1){0.4}} \put(2,-4.8){$c$} 
\put(3.2,-4){\line(0,1){0.4}} \put(3,-4.8){$d_{n+6}$}
\put(5,-4){\line(0,1){0.4}} \put(4.9,-4.8){$y_{n+2}^\#$} 
\put(7.5,-4){\line(0,1){0.4}} \put(7.2,-4.8){$z_{n+1}$}  
\put(10,-4){\line(0,1){0.4}} \put(9.9,-4.8){$d_{n+2}$}  
\put(15,-4){\line(0,1){0.4}} \put(14.8,-4.8){$a_{n+1}$}
\put(18,-4){\line(0,1){0.4}} \put(17.8,-4.8){$z_{n}^\#$}
\put(20.4,-4){\line(0,1){0.4}} \put(20,-4.8){$d_{n+1}^\#$}
\put(18,-2.5){\vector(0,1){1.2}}
\put(19,-2){$f$ }
\endpic\hss}
\vskip 4mm
\caption[ ]{\protect{\label{nsf_nes4}}
{\small Points and their images under $f^{S_{n}}$.
Here $p^\#$ is either $p$ or $\hat p$.
Note that $D_{n+1}\cap \rz=(\hat u_n,u_n)$,
$D_{n+2}\cap \rz=(\hat u_{n+1},u_{n+1})$.
The region $R_{n+1}$ is contained in $D_{n+1}$ by the proof below
the statement of Proposition 8.1. Note that $R_{n+1}\cap \rz\supset
[\hat z_n,z_n]$. In the complex plane, $R_{n+1}$ will contain
a diamont-shaped neighbourhood of $(z_{n+1},a_{n+1})$ and also
a neighbourhood of the half-open interval $[c,z_{n+1})$.
Of course, the scales are quite different for $\ell$ large;
for example $|\hat u_n-d_{n+6}|/|u_n-c|$ of order $1/\ell$.}}
\end{figure}
}

\medskip
In the remainder of this section we shall prove Proposition 8.1.
\bigskip\medskip

\noindent
{\em Proof of Proposition 8.1:}
First note that
$$Df^{S_{n+1}}(z_{n+1})=Df^{S_{n-1}}(z_{n-1}) Df^{S_n}(z_{n+1})$$
and that $f^{S_{2k}}$ and
$f^{S_{2k+1}}$ both converge (up to scaling and orientation)
in the $C^1$ topology by the renormalization results from Section
\ref{sectqcr}, see Theorem~\ref{reno}.
It follows that the left and middle terms of the last inequality
have the same limit (in fact, up to a minus sign as the maps
$f^{S_{n+1}}$ and $f^{S_{n-1}}$ have opposite orientations
but the points $z_{n+1}$ and $z_{n-1}$ are on opposite sides
of $c$)
and therefore that
$$Df^{S_n}(z_{n+1})\to -1.
$$
Notice that $[a_{n+1},z_{n+1}]\ni d_{n+2}$ and 
$[a_{n-1},z_{n-1}]\ni d_n$ lie on opposite sides
of $c=0$. Define $h_n$ to be the (orientation reversing)
affine map with 
$$h_n(z_{n-1})=z_{n+1}\text{ and with }h_n(a_{n-1})=a_{n+1}.$$
By the renormalization result, $h_n$
converges to a scalar contraction map
(with contraction factor of the order $-(1-C/\ell)$).

Before continuing with the proof of the proposition
we shall investigate some properties of
$\psi_n=h_n\circ f^{S_n}\colon [c,a_{n+1}]\to \rz$.
Note that $\psi_n$ also depends on $\ell$.
\medskip

\begin{lemma} \label{psin} The map
$\psi_n=h_n\circ f^{S_n}$ is an orientation reversing
diffeomorphism from $[c,a_{n+1}]$ into $[c,a_{n+1}]$
with fixed point $z_{n+1}$.
Moreover, for each $\ell$ there exists $n_0(\ell)$
such that for $n\ge n_0(\ell)$,
\begin{itemize}
\item
There exists a sequence
$C_n>0$ converging to a positive constant as $n\to \infty$
such that
$$D\psi_n(z_{n+1})=-1+C_n/\ell\,\, ;$$
\item
the basin of $z_{n+1}$ under the map $\psi_n$
contains at least the interval $[c,a_{n+1}]$;
\item
all points in $[c,a_{n+1}]$ are mapped by a uniformly
bounded number of iterates of $\psi_n$
in a neighbourhood of $z_{n+1}$
of size $C/\ell$. More precisely, given
$\epsilon>0$ there exists $m(\epsilon)$
(not depending on $n$ and $\ell$)
such that  $|\psi_n^{m(\epsilon)}(p)-z_{n+1}|\le
\epsilon |z_{n+1}-a_{n+1}|$
for each $p\in [c,a_{n+1}]$.
\end{itemize}
\end{lemma}

\medskip
Note that the choice of
$|z_{n+1}-a_{n+1}|$ in the last part of this lemma is more or
less arbitrary: we could also take $|z_{n+1}-d_{n+2}|$ (or something else)
because of the real bounds.

\bigskip
\noindent
{\em Proof of Lemma 8.1:}
The first assertion follows since  $f^{S_n}$ maps
$[a_{n+1},c]$ onto $[d_n,d_{n+4}]$ and because this
last interval is `well inside'  $[a_{n-1},c]$
in view of the real bounds from Theorem~\ref{realbounds}.
Therefore $\psi_n$ maps $[c,a_{n+1}]$ inside $[c,a_{n+1}]$.
In fact, the image of this interval has length of order $1/\ell$.
Since $Df^{S_n}(z_{n+1})\to -1$,
the orientation reversing diffeomorphism $\psi_n$
has a contracting fixed point in $z_{n+1}$.
Since $|a_{n-1}-c|/|a_{n+1}-c|$ is of order $1+C/\ell$
because of the real bounds, this implies that
$D\psi_n(z_{n+1})=-1+C_n/\ell$ where $C_n>0$ is
uniformly bounded and bounded away from zero
for all large $n$ and $\ell$. Because of the renormalization
result $C_n$ converges to a positive constant as $n$ tends
to infinity.

Let us show that $\psi_n$ attracts all points
in $[c,a_{n+1}]$. Since $\psi_n$ is orientation reversing,
$\psi_n$ would otherwise have a periodic two orbit
in this interval.
But the Schwarzian derivative of $f\colon \rz\to \rz$ is negative.
Hence the Schwarzian derivative of
$\psi_n$ is also negative.
Hence
$$[c,a_{n-1}]\ni z\mapsto |D\psi_n^2(z)|$$
has no positive local minima. It follows that
$\psi_n^2$ cannot have three fixed points
and that
each point in $[c,a_{n+1}]$ is in the basin of $z_{n+1}$.

To prove the last assertion of the lemma, let us assume by contradiction
that such a bound $m(\epsilon)$ did not exist.
Then there exists a sequence of integers
$n,\ell\to \infty$ and a sequence of
points $p_n\in [c,a_{n+1}]$ with
$|p_n-z_{n+1}|\ge \epsilon |z_{n+1}-a_{n+1}|$
(depending on $n$ and $\ell$)
which are almost saddle-nodes:
$$|D\psi^2_n(p_n)|\to 1\text{ and }\frac{|\psi_n^2(J_n)|}{|J_n|}\to 1$$
where $J_n=[z_{n+1},p_n]$.
(The last statement is illustrated in Figure~\ref{nsf_nes5} below.)
To show this is impossible we consider the following cross-ratio:
\beq
B(\psi^2_n;J_n):=\frac{\left[\frac{\psi_n^2(J_n)|}{|J_n|}\right]^2}
{D\psi^2_n(z_{n+1})D\psi_n^2(p_n)}.
\label{cra}
\eeq
Since $\psi_n(z_{n+1})=z_{n+1}$
and $|D\psi_n(z_{n+1})|=1-C_n/\ell$, the previous properties
of $p_n$ imply that the previous cross-ratio gets
arbitrarily close to one for some sequence of integers $n$ and $\ell$
tending to infinity.
However, since $\psi_n=h_n\circ f^{S_n}$,
$$B(\psi_n^2;J_n)=B(\psi_n;\psi_n(J_n))\cdot
B(h_n;f^{S_n}(J_n))\cdot B(f^{S_n-1};f(J_n))\cdot
B(f;J_n).$$
Since $h_n$ is affine (and therefore preserves cross-ratios),
and since $Sf<0$ (and $f$ therefore expands cross-ratios),
the previous cross-ratio is bounded from below by
the
cross-ratio
$$B(f;J_n)=\frac{\left[\frac{|f(J_n)|}{|J_n|}\right]^2}{Df(z_{n+1})Df(p_n)}$$
where $J_n=[z_{n+1},p_n]\subset \rz_{\pm}$. Of course,
this cross-ratio is at least one because $Sf<0$, but in fact more holds:
write $p_n=t_nz_{n+1}$ and $t_n=1+\kappa_n/\ell$.
Then the last cross-ratio
is equal to
$$\frac{\left[\frac{|t_n^\ell-1|}{|t_n-1|}\right]^2}{\ell^2t_n^{\ell-1}}.$$
From the choice of $p_n$
(i.e., from $|p_n-z_{n+1}|\ge \epsilon |z_{n+1}-a_{n+1}|$)
and from the real bounds
one has $\kappa_n\ge C'\epsilon$ for some universal constant
$C'>0$. Hence the limit of the previous expression
as $\ell\to \infty$ is equal to
$\frac{(e^{\kappa_n}-1)^2}{{\kappa_n}^2e^{\kappa_n}}$
which is at least $1+C_0\kappa_n^2$
where $C_0>0$ is a uniform number.
It follows that
$$\frac{\left[\frac{|f(J_n)|}{|J_n|}\right]^2}{Df(z_{n+1})Df(p_n)}
\ge 1+C_0\kappa_n^2\ge
1+ C_0'\epsilon^2$$
where $\kappa_n$ is as above. But as
we had shown before (\ref{cra}) tends to one
if the last property of the lemma does not hold,
contradicting the last inequality.
\qed
\bigskip

\kies{
\begin{figure}[htb]
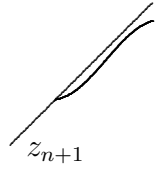
 \hfil \beginpicture
\dimen0=3cm
\setcoordinatesystem units <\dimen0,\dimen0>
\setplotarea x from -0.2 to 0.4, y from -0.2 to 0.4
\setlinear \plot
-0.2 -0.2 0.43 0.43 /
\plot 
0.000  0.000
0.002  0.000
0.004  0.001
0.006  0.001
0.008  0.001
0.010  0.002
0.012  0.002
0.014  0.003
0.016  0.003
0.018  0.004
0.020  0.004
0.022  0.005
0.024  0.006
0.026  0.006
0.028  0.007
0.030  0.008
0.032  0.008
0.034  0.009
0.036  0.010
0.038  0.011
0.040  0.012
0.042  0.012
0.044  0.013
0.046  0.014
0.048  0.015
0.050  0.016
0.052  0.017
0.054  0.018
0.056  0.019
0.058  0.020
0.060  0.021
0.062  0.023
0.064  0.024
0.066  0.025
0.068  0.026
0.070  0.027
0.072  0.029
0.074  0.030
0.076  0.031
0.078  0.032
0.080  0.034
0.082  0.035
0.084  0.037
0.086  0.038
0.088  0.039
0.090  0.041
0.092  0.042
0.094  0.044
0.096  0.045
0.098  0.047
0.100  0.049
0.102  0.050
0.104  0.052
0.106  0.053
0.108  0.055
0.110  0.057
0.112  0.059
0.114  0.060
0.116  0.062
0.118  0.064
0.120  0.066
0.122  0.067
0.124  0.069
0.126  0.071
0.128  0.073
0.130  0.075
0.132  0.077
0.134  0.078
0.136  0.080
0.138  0.082
0.140  0.084
0.142  0.086
0.144  0.088
0.146  0.090
0.148  0.092
0.150  0.094
0.152  0.096
0.154  0.098
0.156  0.100
0.158  0.103
0.160  0.105
0.162  0.107
0.164  0.109
0.166  0.111
0.168  0.113
0.170  0.115
0.172  0.117
0.174  0.120
0.176  0.122
0.178  0.124
0.180  0.126
0.182  0.128
0.184  0.131
0.186  0.133
0.188  0.135
0.190  0.137
0.192  0.140
0.194  0.142
0.196  0.144
0.198  0.146
0.200  0.149
0.202  0.151
0.204  0.153
0.206  0.155
0.208  0.158
0.210  0.160
0.212  0.162
0.214  0.165
0.216  0.167
0.218  0.169
0.220  0.171
0.222  0.174
0.224  0.176
0.226  0.178
0.228  0.181
0.230  0.183
0.232  0.185
0.234  0.188
0.236  0.190
0.238  0.192
0.240  0.194
0.242  0.197
0.244  0.199
0.246  0.201
0.248  0.204
0.250  0.206
0.252  0.208
0.254  0.210
0.256  0.213
0.258  0.215
0.260  0.217
0.262  0.219
0.264  0.222
0.266  0.224
0.268  0.226
0.270  0.228
0.272  0.230
0.274  0.233
0.276  0.235
0.278  0.237
0.280  0.239
0.282  0.241
0.284  0.243
0.286  0.245
0.288  0.248
0.290  0.250
0.292  0.252
0.294  0.254
0.296  0.256
0.298  0.258
0.300  0.260
0.302  0.262
0.304  0.264
0.306  0.266
0.308  0.268
0.310  0.270
0.312  0.272
0.314  0.274
0.316  0.276
0.318  0.278
0.320  0.279
0.322  0.281
0.324  0.283
0.326  0.285
0.328  0.287
0.330  0.289
0.332  0.290
0.334  0.292
0.336  0.294
0.338  0.296
0.340  0.297
0.342  0.299
0.344  0.301
0.346  0.302
0.348  0.304
0.350  0.305
0.352  0.307
0.354  0.308
0.356  0.310
0.358  0.312
0.360  0.313
0.362  0.314
0.364  0.316
0.366  0.317
0.368  0.319
0.370  0.320
0.372  0.321
0.374  0.323
0.376  0.324
0.378  0.325
0.380  0.326
0.382  0.328
0.384  0.329
0.386  0.330
0.388  0.331
0.390  0.332
0.392  0.333
0.394  0.334
0.396  0.335
0.398  0.336
0.400  0.337
0.402  0.338
0.404  0.339
0.406  0.340
0.408  0.341
0.410  0.342
0.412  0.343
0.414  0.344
0.416  0.344
0.418  0.345
0.420  0.346
0.422  0.347
0.424  0.347
0.426  0.348
0.428  0.348
0.430  0.349
0.432  0.350
0.434  0.350
/
\put {$z_{n+1}$} <0mm,-7mm> at 0 0
\endpicture
\caption[ ]{\protect{\label{nsf_nes5}}{\small
An almost saddle-node point implies that the cross-ratio expansion
on some interval is close to one.}}	
\end{figure}
}

Now we will continue with the proof of 
the proposition. For this we will analyze the map $\psi_n$
considered as a map on a complex neighbourhood of
$z_{n+1}$.

\begin{lemma}\label{psin2}
For each sufficiently large $\ell$
there exists $n_0(\ell)$ and
a set $\hat V_{n,\ell}$ for $n\ge n_0(\ell)$.
This set intersects the real line in a closed interval containing
$z_{n+1}$ and
\begin{itemize}
\item $\psi_n$ maps $\hat V_{n,\ell}$ into the interior of $\hat V_{n,\ell}$
for each $n\ge n_0(\ell)$;
\item
there exists $\alpha\in (\pi/2,\pi)$, $\epsilon>0$
and two real intervals $J_1,J_2$ of length $\ge \epsilon |z_{n+1}-a_{n+1}|$
with unique common point $z_{n+1}$ such that
$$\hat V_{n,\ell}\supset D(J_1,\alpha) \bigcup D(J_2,\alpha)$$
provided $n\ge n_0(\ell)$;
\item $\hat V_{n,\ell}$ is contained in
$D_*(\hat V_{n,\ell}\cap \rz)$.
\end{itemize}
\end{lemma}
\noindent
{\em Proof of Lemma 8.2:}
Consider the affine map
$sc_n\colon \cz\to \cz$
which sends $[z_{n+1},a_{n+1}]$ to $[0,1]$.
Note that $|D sc_n|$ is of order $\ell/|z_n-c|$
because of the real bounds
and, in fact,
$D sc_n(c)\to -\infty$ as $n$ or $\ell$ tends to infinity.
Now define
$$\Psi_n=sc_n\circ \psi_n\circ sc_n^{-1}.$$
Then
$\Psi_n(0)=0$, $D\Psi_n(0)=(-1+C/\ell)$ and
by renormalization $\Psi_n$ converges as $n\to \infty$
to some function $\Psi$.
(Of course, $\Psi$ might depend on $\ell$.)
Also $S\Psi_n<-\delta<0$. Indeed,
\beq
S(g_1\circ g_2)=Sg_1(g_2)(Dg_2)^2+Sg_2.
\label{sco}
\eeq
Moreover, an explicit calculation gives that there exists $C>0$
with
$$Sf(z)\le -C \frac{\ell^2}{|z_n-c|^2}$$
and, therefore,
$$Sf^{S_n}(z)\le -C \frac{\ell^2}{|z_n-c|^2}$$
for each $z\in [u_{n+4},u_{n-4}]$ (or in fact,
any similar interval).
Since $h_n$ is an affine map close to the identity,
also
$$S\psi_n \le -C \frac{\ell^2}{|z_n-c|^2}.$$
By the real bounds, the length of $[z_{n+1},a_{n+1}]$
and therefore $|D\, sc_n^{-1}|$ is of
the order of $|z_n-c|/\ell$. It follows by (\ref{sco})
and the definition of $\Psi_n$ that there exists $C>0$ such that
$$S\Psi_n(z)\le -C$$
for all $z$ in, say, $[-1,1]$.

By the renormalization result,
\beq
\Psi_n(z)=(-1+\alpha_n'/\ell)z+\beta_n' z^2 + \gamma_n' z^3 + O(z^4)
\label{bcoe}
\eeq
where the coefficients $\alpha_n',\beta_n',\gamma_n'$ do depend
on $\ell$ but for each fixed $\ell$ converge to constants
as $n\to \infty$.
In fact, $\alpha_n'$ converges to a positive constant.

\noindent
{\em Claim:} the coefficients as well as the remainder term is bounded uniformly
at the origin: $|O_n(z^4)|\le C|z^4|$ for some universal $C$
for all $z$ in, say, $sc_n([d_{n+6},a_{n+1}])\ni 0$
(notice that the distance of each of
the endpoints of the interval $sc_n([d_{n+6},a_{n+1}])\ni 0$
to the origin is by the real bounds of order one).

To prove this claim, note that
the coefficients and the remainder term can be estimated uniformly in
$n$ and $\ell$ because of the Taylor theorem.
Indeed, this theorem shows that the remainder term in (\ref{bcoe})
is bounded by the fourth derivative of $\Psi_n$
in the $sc_n([d_{n+6},a_{n+1}])$.
Now $f^{S_n}$ has a univalent extension
from a neighbourhood $V$ of $[d_{n+6},a_{n+1}]$
onto a (bounded) disc $\hat W$ containing the disc
$$W=D_*(f^{S_n}(d_{n+6}),f^{S_n}(a_{n+1})]=D_*(f^{S_n}(d_{n+6}),d_{n+4})$$
and with the same centre. Due to the real bounds, see Figures~\ref{nsf_nes3}
and \ref{nsf_nes4}, we can make
choose $\hat W$ so that its radius is a definite factor larger
than the radius of $W$. Hence by the Koebe Lemma,
the distance of the boundary of $W$ is of the same order as the size of,
say, the interval $[d_{n+6},a_{n+1}]$.
Hence, by applying the scaling map $sc_n$,
all these discs and intervals get a size of unit order
and the distance of the boundary
of $\tilde V=sc_n(V)$ to $sc_n([d_{n+6},a_{n+1}])$
is also of the order one.
Moreover, $\tilde V$ is mapped by $\Psi_n$
into a uniformly bounded disc. By the Cauchy integral
formula,
$$\Psi_n^{(k)}(z)=\frac{k!}{2\pi i}\oint
\frac{\Psi_n(t)}{(t-z)^{k+1}} \, dt.$$
Because the distance of $\tilde V$ to a point in $sc_n([d_{n+6},a_{n+1}])$
is of order one, we can take a circle around $z$ inside $\tilde V$
with a radius which is uniformly bounded from below. Since $\Psi_n$
is uniformly bounded on $\tilde V$, it follows that
$\Psi_n^{(i)}$,  $i=1,2,3,4$ are bounded on $\tilde V\supset
sc_n([d_{n+6},a_{n+1}])$. This completes the proof of the
claim.

Because of (\ref{bcoe}), by an explicit calculation,
\beqas
\Psi_n^2(z)
&=&
(-1+\alpha_n'/\ell)^2z+(\beta_n'\alpha_n'/\ell) z^2-\gamma_n'z^3+O(z^4)
\\
&=&(1-\frac{\alpha_n}{\ell})z+\frac{\beta_n}{\ell}z^2
-\gamma_n z^3 + O(z^4)
\eeqas
where $\alpha_n>0$ and $\gamma_n$ converge to positive constants
and $\beta_n$ also converges to a constant for each fixed $\ell$ as
$n\to \infty$.
(That $\gamma_n$ converges to a {\it positive} constant as $n\to \infty$
provided $\ell$ is large, follows
from the fact that $S\Psi_n\le -\delta<0$ and from the fact
that $\frac{\alpha_n}{\ell}\to 0$ and
$\frac{\beta_n}{\ell}\to 0$ as $n,\ell\to \infty$.)

This family of maps seem similar as $\ell\to \infty$ 
to the well-known map
$$z\mapsto z-z^3$$
which has a neutral fixed point at $0$, having a basin containing two
petals attached to $0$.
Now we will study our family of maps $\Psi_n^2$
in the same way as is commonly done
for this `limit' map and show that
 the basin of the fixed point is also not too small.

For this we will introduce new coordinates
$w=1/z^2$ and send the origin in a two-fold way to infinity.
In these new $w$ coordinates our map becomes 
\beq
\Theta(w)   = 
w+ \frac{\tilde \alpha_n}{\ell}w
-\frac{\tilde \beta_n}{\ell}\sqrt{w}
+ \tilde \gamma_n
+ O(\sqrt{w^{-1/2}}).
\label{tr1}
\eeq
Here $\tilde \alpha_n,\tilde \beta_n,\tilde \gamma_n$
are constants with $\tilde \gamma_n$ and $\tilde \alpha_n$ converging to
positive limits (which are uniformly bounded and bounded away from
zero for all $\ell$ large).
It follows from (\ref{tr1})
that if $\gamma_1$ is sufficiently large,
then there exists $\delta_1$ such that
\beq
|w|\ge \gamma_1\ell\text{ implies }
|\Theta(w)| \ge |w|+\delta_1
\label{tr2}
\eeq
for $\ell$ sufficiently large (here $\gamma_1$ and
$\delta_1>0$ are universal constants).
This holds because for such $w$, the second
term on the right hand side of (\ref{tr1}) dominates
the last three terms since
$\liminf \tilde \alpha_n>0$ (provided $\gamma_1$ and $\ell,n$ are large).
On the other hand, (\ref{tr1})
also implies that if $|w|\le \gamma_1\ell$ and $Re(w)=\gamma_2$
then there exists a universal constant $\delta_2>0$ such that 
\beq
Re(\Theta(w))\ge Re(w)+\delta_2
\label{tr3}
\eeq
provided $\gamma_2$ and $\ell$ are sufficiently large.
Indeed, for $|w|\le \gamma_1\ell$ and $Re(w)=\gamma_2>0$,
for the second, third and last term on the right hand side of (\ref{tr1})
one has
\beq
Re\left(\frac{\tilde \alpha_n}{\ell}w\right) >0\text{ , }
Re\left(\frac{\tilde \beta_n}{\ell}\sqrt{w}\right)\to 0 \text{ and }
Re\left(O(w^{-1/2})\right)\le \text{Const} |\gamma_2|^{-1/2}
\label{tr4}
\eeq
for $n$ and $\ell$ large. Because the third term in (\ref{tr1})
is uniformly bounded away from zero, one gets (\ref{tr3}).
Combining (\ref{tr2}) and (\ref{tr3}) one gets that
the basin of $\infty$ contains points outside a big circle
with radius $\ell$ (and centered at the origin)
and also the points to the right of the vertical line
$l$ given by $Re(w)=\gamma_1$.
Hence the basin of the origin under the map $\Psi_n$ is not too small.
In fact, it is the union of a disc of radius of the
order $1/\sqrt{\ell}$ (a strongly related and crucial estimate
will reappear in Section~\ref{asytex}!!) and a rotated figure eight region
(which does not depend on $\ell$), see Figure~\ref{nsf_nes6}.
The figure eight is the inverse of the vertical line under the transformation
$z\mapsto 1/z^2$. Let us call this set $\hat V$.

Now $sc_n^{-1}(\hat V)$ satisfies by construction the
first and second property announced in the lemma.
The remaining task is show that
$\hat V\subset D_*(\hat V\cap \rz)$.
But this can be seen as follows. The inverse of the
region to the right of the
line $l$ (so the region $Re(w)\ge \gamma_1$) under the map
$z\mapsto 1/z$ is a symmetric disc through $0$ and $1/\gamma_1$.
Similarly, the inverse  of the region outside the
circle with radius $\gamma_1\ell$ (centered at
$0$) under the map $z\mapsto 1/z$
is a disc centred in $0$ and with radius $1/\ell$.
Hence the union of the inverses under $z\mapsto 1/z^2$ of these regions
is contained in $D_*(\hat V\cap \rz)$.
\qed

\kies{
\begin{figure}[htp]
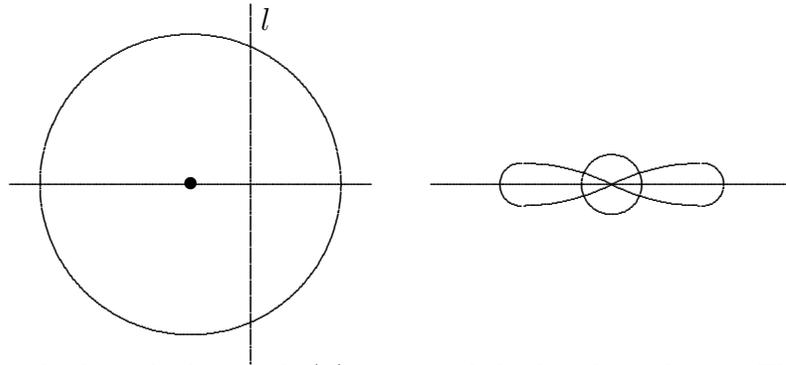
 \hfil
\beginpicture
\dimen0=0.4cm
\setcoordinatesystem units <\dimen0,\dimen0>
\setplotarea x from -6 to 6, y from -5 to 5
\setlinear
\plot -6 0 6 0 /
\put {$\bullet$} at 0 0
\put {$l$} at 2.5 5.5
\setsolid
\circulararc 360 degrees from 5 0  center at  0 0
\plot 2 -6 2 6 /
\setcoordinatesystem units <\dimen0,\dimen0> point at -14 0
\setplotarea x from -6 to 6, y from -5 to 5
\setlinear
\plot -6 0 6 0 /
\setsolid
\circulararc 360 degrees from 1 0  center at  0 0
\circulararc 180 degrees from 3 -0.72  center at  3 0
\circulararc 180 degrees from -3 0.72  center at  -3 0
\circulararc -26 degrees from 0 0 center at 3 -6
\circulararc 26 degrees from 0 0 center at 3 6
\circulararc 26 degrees from 0 0 center at -3 -6
\circulararc -26 degrees from 0 0 center at -3 6
\endpicture
\caption[ ]{\protect{\label{nsf_nes6}}
{\small On the left, the big circle $|z|=\gamma_1 \ell$
and the line $l$ are drawn.
The the right, the images under the map $w=1/z^2$ of these regions
are drawn schematically.}}
\end{figure}

}

\bigskip

\noindent
{\em Proof of Proposition 8.1:}
The previous lemma has still one shortcoming for our purpose:
it is still not guaranteed that $\hat V_{n,\ell}$ contains
the interval $[z_{n+1},a_{n+1}]$.
However, because of the last statement of Lemmma 8.1,
there exists a universal integer $m$ and a constant $\gamma_2>0$
in the proof of the
previous lemma so that
$\hat V_{n,\ell}\cap \rz=\psi_n^m[z_{n+1},a_{n+1}]$.
Now the inverse of $\psi_n^m$ has
bounded distortion on $D_*([z_{n+1},a_{n+1}])$
(it has a univalent extension on a definite neighbourhood of
this disc because of the real bounds).

\kies{
\begin{figure}[htp]
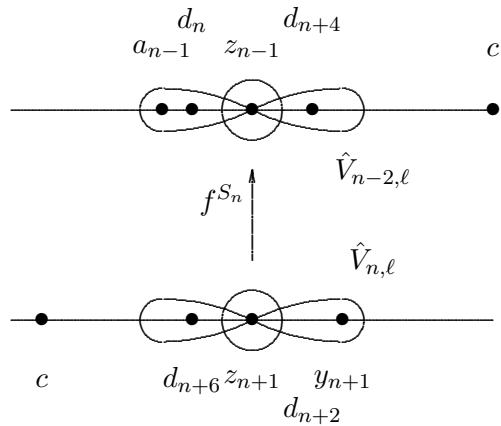
 \hfil
\beginpicture
\dimen0=0.4cm
\setcoordinatesystem units <\dimen0,\dimen0>
\setplotarea x from -8 to 10, y from -5 to 5
\setlinear
\plot -8 0 8 0 /
\setsolid
\circulararc 360 degrees from 1 0  center at  0 0
\circulararc 180 degrees from 3 -0.72  center at  3 0
\circulararc 180 degrees from -3 0.72  center at  -3 0
\circulararc -26 degrees from 0 0 center at 3 -6
\circulararc 26 degrees from 0 0 center at 3 6
\circulararc 26 degrees from 0 0 center at -3 -6
\circulararc -26 degrees from 0 0 center at -3 6
\multiput {$\bullet$} at -3 0 -2 0  0 0  2 0 8 0 /
\put {\small $a_{n-1}$} at -3 2
\put {\small $d_n$} at -2 3
\put {\small $z_{n-1}$} at 0 2
\put {\small $d_{n+4}$} at 2 3
\put {\small $c$} at 8 2
\put {\small $\hat V_{n-2,\ell}$} at 4 -2
\setcoordinatesystem units <\dimen0,\dimen0> point at 0 7
\setplotarea x from -10 to 6, y from -5 to 5
\setlinear
\plot -8 0 8 0 /
\setsolid
\circulararc 360 degrees from 1 0  center at  0 0
\circulararc 180 degrees from 3 -0.72  center at  3 0
\circulararc 180 degrees from -3 0.72  center at  -3 0
\circulararc -26 degrees from 0 0 center at 3 -6
\circulararc 26 degrees from 0 0 center at 3 6
\circulararc 26 degrees from 0 0 center at -3 -6
\circulararc -26 degrees from 0 0 center at -3 6
\multiput {$\bullet$} at -7 0 -2 0  0 0  3 0 /
\put {\small $c$} at -7 -2
\put {\small $d_{n+6}$} at -2 -2
\put {\small $y_{n+1}$} at 3 -2
\put {\small $d_{n+2}$} at 2 -3
\put {\small $z_{n+1}$} at 0 -2
\put {\small $\hat V_{n,\ell}$} at 4 2
\arrow <5pt> [0.2,0.4] from 0 2 to 0 5
\put {\small $f^{S_n}$} at -1 4
\endpicture
\caption[ ]{{\small The region $V_{n,\ell}$ mapped into the
region $V_{n-2,\ell}$ for large $n$.}}
\end{figure}

}

Now $\hat V_{n,\ell}$ fits inside $D_*(\hat V_{n,\ell}\cap \rz)$.
It follows by Proposition~\ref{corkoebe} that
$V_{n,\ell}\, :=\psi_n^{-m}(\hat V_{n,\ell})$
satisfies
\begin{itemize}
\item $\psi_n$ maps $V_{n,\ell}$ into the interior of $V_{n,\ell}$
for large $n$;
\item
there exists $\alpha\in (\pi/2,\pi)$, such that
$$\hat V_{n,\ell}\supset D([z_{n+1},a_{n+1}];\alpha)$$
provided $n\ge n_0(\ell)$ (in fact, it also contains the other
part of the figure eight, but this part is not needed here.)
\end{itemize}
Let us interpret this information on $\psi_n=h_n\circ f^{S_n}$
for the maps $f^{S_n}$.
Because of the last property we get
for each $\ell$ sufficiently large an integer
$n_0(\ell)$ such that 
$f^{S_n}$ maps $V_{n,\ell}$ into the interior of $V_{n-2,\ell}$
for $n\ge n_0(\ell)$.
(Here we use the renormalization result that $f^{S_n}$ converges
up to scaling.)
Thus we have proved Proposition 8.1.
\qed

\sect{An induced mapping with Markov properties}
\label{markov}
In this section we shall use the previous
discs to define an induced map
which has nice Markov properties.
Let $D_n$ be the discs from the previous section
and define for $n\ge k_0$,
\[
A_n=D_{n+1}\setminus D_{n+2}
\]
and let
$A_n'$ to be the annulus $A_n$ minus the disc
$D_{n+1}^1\subset D_{n+1}\setminus D_{n+2}$.
Then
\[
f^{S_{n}-1}\text{ maps }D_{n+1}^f\text{ diffeomorphically onto }D_{n-1}
\]
and
\[
f^{S_{n}-1}\text{ maps }D_{n+2}^f\text{ diffeomorphically onto }D_{n}^1
\]
So
\[
f^{S_n}\text{ maps }A_n\text{ as an }\ell\text{-fold covering
onto }D_{n-1}\setminus D_{n}^1.
\]
This last set is equal to
$$\left(\cup_{i\ge n}A_i \right)
\,\, \bigcup \,\, A_{n-1}'\,\, \bigcup \,\, A_{n-2}.$$
Moreover, $f^{S_n}$ maps $A_n'$ again onto
$\left(\cup_{i\ge n}A_i \right)
\,\, \bigcup \,\, A_{n-1}'\,\, \bigcup \,\, A_{n-2}$
but now the map only
$(\ell-1)$-covers $\cup_{i\ge n}A_i$
(because the missing disc would also have been mapped
diffeomorphically onto $\cup_{i\ge n}A_i$)
while it is still $\ell$-covers the remaining part of the
target.

To formalize all this we define
$X_n$ to be the disjoint union of
$A_n$ and $A_n'$ and $X=\cup X_n$.
Define
\[
F\,\,\colon \,\,\cup X_n \to
\cup X_n
\]
by
\[
F|A_n=f^{S_n}\ \text{ and }F|A_n'=f^{S_n}.
\]
for $n\ge k_0$ and $F|(A_{k_0-1}\cup A_{k_0-2})=id$.
Moreover, let $\A$ be the partition
of $X$ into sets $A_n$ and $A_n'$.
Then $F\colon X\to X$ is Markov map with respect
to this partition. It sends each
element of the partition $\A$ as
a covering map onto a union of elements.
Now we will iterate $F$.
So define the partition
\[
\A_{n+1}=\A\vee F^{-1}(\A) \vee \dots \vee F^{-n}(\A).
\]
Then $\A_1=\A$. If $A$ is an element from
$\A_2$ then $B=F(A)\in \A_1$
and $F\colon A\to B$ is a {\it covering map}:
\begin{itemize}
\item $F\colon A\to B$ is a local homeomorphism;
\item there exists $k$ such that for each $y\in B$,
$$\# (F|A)^{-1}(y)=k.$$
\end{itemize}
Indeed, assume $A\subset A_r$ or $A\subset A_r'$.
Then
$$k=\left\{
\begin{array}{ll}
{\ell}   &{\text{ if }B=A_{r-2},}\\
{0}      &{\text{ if }B=A_{r-2}',}\\
{\ell}  &{\text{ if }B=A_{r-1}',}\\
{0}     &{\text{ if }B=A_{r-1},}\\
{1}     &{\text{ if }B=A_{r+j},}\\
{0}     &{\text{ if }B=A_{r+j}',}
\end{array}
\right.
$$
where $j\ge 0$.
Of course, there are $\ell$ components $A\in \A_2$ inside
$(F|A_r)^{-1}(A_{r+j})$
whereas there are only $\ell-1$ components
$A\in \A_2$ in $(F|A_r')^{-1}(A_{r+j})$.
However, this will play no role in the future.
In fact,
$$F\colon A_r\to \left(\cup_{i\ge r}A_i \right)
\,\, \bigcup \,\, A_{r-1}'\,\, \bigcup \,\, A_{r-2}$$
is also a covering map, while
$$F\colon A_r'\to \left(\cup_{i\ge r}A_i \right)
\,\, \bigcup \,\, A_{r-1}'\,\, \bigcup \,\, A_{r-2}$$
is not. This is because in the latter case, points from
$\left(\cup_{i\ge r}A_i \right)$ are covered precisely $\ell-1$ times
(due to the missing disc in $A_r'$)
whereas each point from $A_{r-1}'\,\, \bigcup \,\, A_{r-2}$ is
covered $\ell$ times.

It follows from this that for each component $A$ from an element of
$\A_{k+1}$, \,\, $F^{k}$ is also a covering map
from $A$ to one element of $\A$. If $F^k(A)=A_r$ or $F^k(A)=A_r'$ then
$F^{k+1}$ maps $A$ onto
$$\left(\cup_{i\ge r}A_i \right)
\,\, \bigcup \,\, A_{r-1}'\,\, \bigcup \,\, A_{r-2}.$$
In fact, this map need not be a covering map
when $F^k(A)=A_r'$ because -- as we pointed out above --
the map $F|A_r'$ is not.

\kies{
\begin{figure}[htp]
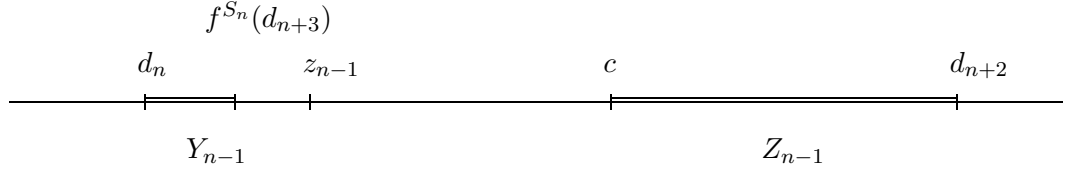

\vskip 0.7cm
\hbox to \hsize{\hss\unitlength=5mm
\beginpic(20,4)(-20,7) \let\ts\textstyle
\put(4,8){\line(-1,0){28}}
\put(1.2,7.8){\line(0,1){0.4}} \put(1,8.8){\small $d_{n+2}$}
\put(-4,6.5){\small $Z_{n-1}$} \put(-8,8.1){\line(1,0){9.2}}
\put(-8,7.8){\line(0,1){0.4}} \put(-8.2,8.8){\small $c$}
\put(-16,7.8){\line(0,1){0.4}} \put(-16.2,8.8){\small $z_{n-1}$}
\put(-18,7.8){\line(0,1){0.4}} \put(-18.8,10){\small $f^{S_{n}}(d_{n+3})$}
\put(-19.3,6.5){\small $Y_{n-1}$} \put(-20.4,8.1){\line(1,0){2.4}}
\put(-20.4,7.8){\line(0,1){0.4}} \put(-20.6,8.8){\small $d_{n}$}
\endpic\hss}
\vskip 0.4cm
\caption[ ]{{\small The intervals $Y_{n-1}$  and $Z_{n-1}$.
Note that $Y_{n+1}\subset Z_{n-1}$, see Figure~\ref{nsf_nes1}.
Whether the intervals $Y_n$ and $Y_{n+1}$ lie on the same side of $c$
is determined by the parity of $n$.}}
\end{figure}

}
 
We will have to analyze the distortion of iterates of $F^{k+1}|A$.
For this we want to apply the Koebe Lemma
and so we would like to see how much one can extend
$F$. Therefore we
associate to $A_n$ and $A_n'$ the following slit-regions.
Define
$$Z_n=[d_{n+3},c]\text{ and }
Y_n=[d_{n+1},f^{S_{n+4}}(d_{n+1})]=[d_{n+1},f^{n+1}(d_{n+4}],$$
(so $Y_n\subset D_{n+1}^1\subset A_n$ and $Z_n\subset \cup_{i\ge n+1}A_i$).
Note that
$$f^{S_n}(Z_n)=[d_{n},f^{S_{n+3}}(d_n)]=Y_{n-1}
\text{ and }
f^{S_n}(Y_n)=[d_{n+2},y_{n+2}]\subset [d_{n+2},c]=Z_{n-1}.$$
Moreover, let
$$
Slit_n=\cz_{(d_n,\hat d_n)} \setminus Z_n
$$
and
$$
Slit'_n=
\cz_{(d_n,\hat d_n)} \setminus
\left( Y_n \cup Z_n \right)
$$
and
$$Slit^*_n=
\cz_{(d_n,\hat d_n)} \setminus Y_{n+1}
$$
Note that
$$A_n\subset Slit_n\text{ , }A_n'\subset Slit'_n$$
and
$$\cup_{i\ge n}A_i\cup A_{n-1}'\cup A_{n-2}\subset Slit^*_{n-2}.$$
Let us study the extensions of $F|A_n$ and $F|A_n'$.
\bigskip

\kies{
\begin{figure}[htp]
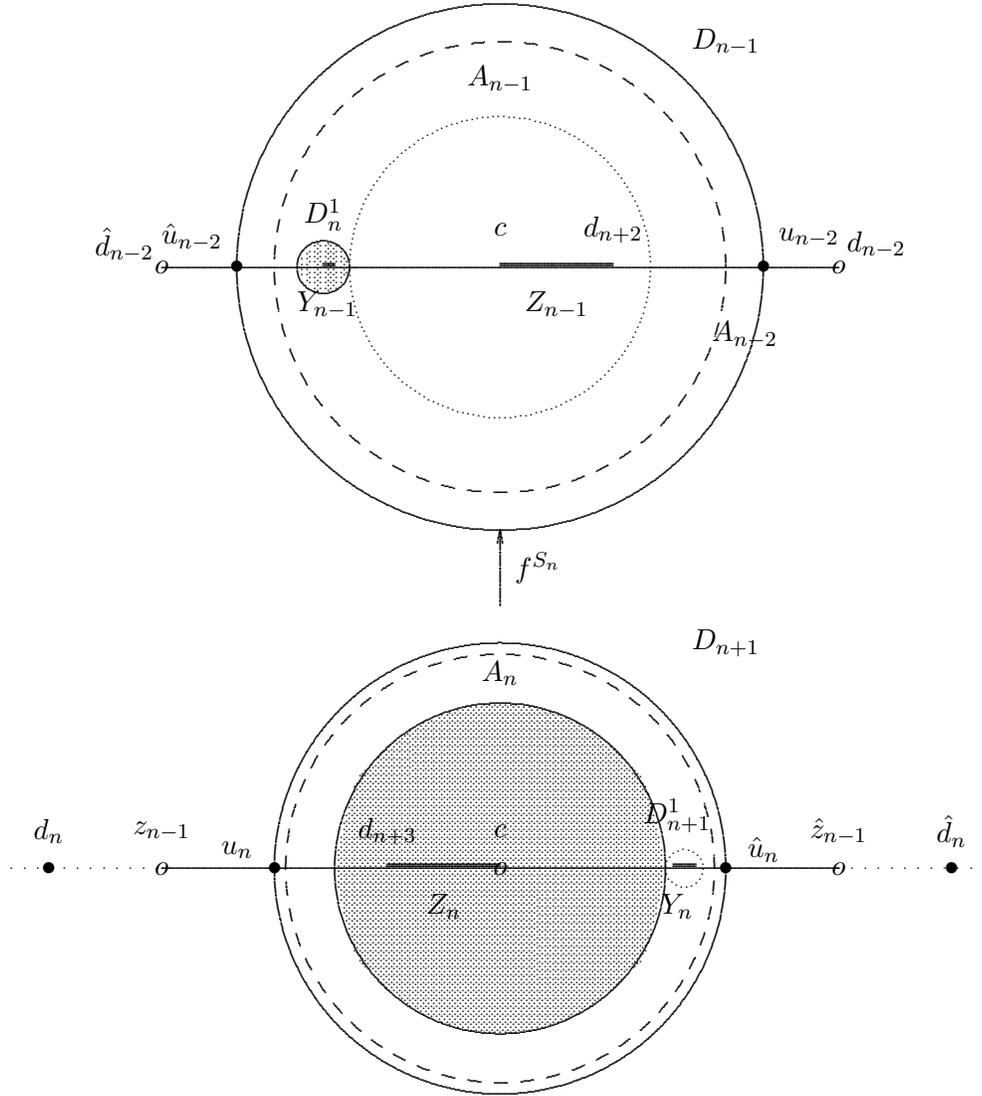
 \hfil
\beginpicture
\dimen0=0.5cm
\setcoordinatesystem units <\dimen0,\dimen0>
\setplotarea x from -9 to 9, y from -7 to 7
\setlinear
\plot -9 0 9 0 /
\setsolid
\circulararc 360 degrees from 7 0  center at  0 0
\setdashes
\circulararc 360 degrees from 6 0  center at  0 0
\setdots <2pt>
\circulararc 360 degrees from 4 0  center at  0 0
\setsolid
\circulararc 360 degrees from -4 0  center at  -4.7 0
\multiput {\small $\bullet$} at -7 0 7 0 /
\multiput {\small $o$} at -9 0 9 0 /
\put {\small $c$} at 0 1
\put {\small $d_{n+2}$} at 3 1
\put {\small $Z_{n-1}$} at 1.5 -1
\put {\small $u_{n-2}$} at 8.2 0.8
\put {\small $\hat u_{n-2}$} at -8.2 0.8
\put {\small $A_{n-2}$} at 6.5 -1.8
\put {\small $A_{n-1}$} at 0 5
\put {\small $D_n^1$} at -4.7 1.4
\put {\small $D_{n-1}$} at 6 6
\put {\small $d_{n-2}$} at 10 0.6
\put {\small $\hat d_{n-2}$} at -10 0.6
\setshadegrid span <1pt>
\setlinear \vshade  -5.4 0 0 -5.2 -0.4 0.4 -5 -0.6 0.6 -4.7 -0.7 0.7 -4.3 -0.6 0.6 -4 0 0 /
\setsolid
\setlinear
\plot 3 0.05 0 0.05 /
\plot -4.4 0.05 -4.7 0.05 /
\put {\small $Y_{n-1}$} at -4.6 -1
\plot 3 0.1 0 0.1 /
\plot -4.4 0.1 -4.7 0.1 /
\arrow <5pt> [0.2,0.4] from 0 -9 to 0 -7
\put {\small $f^{S_n}$} at 1 -8
\setcoordinatesystem units <\dimen0,\dimen0> point at 0 16
\setplotarea x from -9 to 9, y from -7 to 7
\setlinear
\plot -9 0 9 0 /
\setdots
\plot -13 0 -9 0 /
\plot 9 0 13 0 /
\setsolid
\circulararc 360 degrees from 6 0  center at  0 0
\setsolid
\circulararc 360 degrees from 4.4 0  center at  0 0
\setdots <2pt>
\circulararc 360 degrees from 4.4 0  center at  4.9 0
\setdashes
\circulararc 360 degrees from 5.7 0  center at  0 0
\setsolid
\setlinear
\plot -3 0.05 0 0.05 /
\plot 4.6 0.05 5.2 0.05 /
\put {\small $D_{n+1}$} at 6 6
\put {\small $Y_n$} at 4.7 -1
\put {\small $Z_n$} at -1.5 -1
\put {\small $d_{n+3}$} at -3 1
\plot -3 0.1 0 0.1 /
\plot 4.6 0.1 5.2 0.1 /
\multiput {\small $\bullet$} at -6 0 6 0  -12 0 12 0 /
\multiput {\small $o$} at 0 0 -9 0 9 0 /
\put {\small $z_{n-1}$} at -9 1
\put {\small $\hat z_{n-1}$} at 9 1
\put {\small $\hat d_n$} at 12 1
\put {\small $d_n$} at -12 1
\put {\small $u_{n}$} at -7 0.5
\put {\small $\hat u_{n}$} at 7 0.5
\put {\small $A_n$} at 0 5.2
\put {\small $D_{n+1}^1$} at 4.7 1.4
\put {\small $c$} at 0 1
\setquadratic
\vshade -4.4 0 0
<,z,,> -4 -1.82 1.82 -3.5 -2.66 2.66
<z,z,,> -3 -3.2 3.2 -2 -3.9 3.9
<z,z,,>   0 -4.4 4.4 2 -3.9 3.9
<z,z,,> 3 -3.2 3.2 3.5 -2.66 2.66
<z,,,> 4  -1.82 1.82  4.4  0  0 /
\endpicture
\caption[ ]{{\small The image under $f^{S_n}=F$ of the annulus
$A_n$ on the bottom is an `asymetric' annulus with a `small' disc
$D^1_n$ in the annulus $A_{n-1}=D_n\setminus D_{n+1}$
(bounded by a dotted and dashed curve)
removed from the disc $D_{n-1}=\cup_{i\ge n-2}A_i$.
So $F(A_n)=\cup_{i\ge n}A_i\cup A_{n-1}'\cup A_{n-2}$.
The slits $Z_{n-1}=[d_{n+2},c]$ and $Y_{n-1}=[d_n,f^{S_{n+3}}(d_n)]
\subset D^1_n\subset A_{n-1}$ are drawn with thicker lines.}}
\end{figure}

}

\begin{theo}\label{extend}
If $F(A_n)\supset A_m$
then there exists a region $E_n$ with
\begin{itemize}
\item $E_n\subset Slit_n$;
\item $F(E_n)=Slit_m$ and
\item $(F|A_n)^{-1}(A_m)\supset E_n$.
\end{itemize}
The same statement holds if we replace $A_n$ by $A_n'$
or $A_m$ by $A_m'$ provided we then replace
$Slit_n$ by $Slit_n'$ respectively $Slit_m$ by $Slit_m'$.

Moreover, if $H\in \A_{k+1}$ is contained in
$A_n$ then there exists $r$ such that $F^k(H)=A_r$
or $F^k(H)=A_r'$ and there exists a region $E$
such that
\begin{itemize}
\item $H \subset E\subset Slit_n$
\item $F^{k+1}$ is a covering from $H$ onto
$$\cup_{i\ge r}A_i\cup A_{r-1}'\cup A_{r-2}\subset Slit^*_{r-2}$$
and this map covers $E$ onto $Slit^*_{r-2}$.
\item $F^{k+1}(E)=Slit^*_{r-2}$.
\end{itemize}
\end{theo}
\pr
$f^{S_n-1}$ maps $[z_{n-1}^f,t_n^f]$ diffeomorphically
to $[d_{n-2},d_{n-4}]$.
Therefore, $f^{S_n-1}$ maps
some region in
$$\cz_{[z_{n-1}^f, t_n^f]}
$$
onto
$$\cz_{[d_{n-2},\hat d_{n-4}]}.$$
Since $Y_n^f,Z_n^f\subset [z_{n-1}^f,t_n^f]$
and since
$$f^{S_n}(Y_n)\subset Z_{n-1}\text{ and } f^{S_n}(Z_n)= Y_{n-1}$$
it follows that
$f^{S_n-1}$ maps
some region in
$$\cz_{[z_{n-1}^f,t_n^f]}\setminus Y_n^f $$
onto
$$\cz_{[d_{n-2},d_{n-4}]}\setminus Z_{n-1}.$$
Similarly, $f^{S_n-1}$ maps some region in
$$\cz_{[z_{n-1}^f,t_n^f]}\setminus (Y_n^f \cup Z_n^f)$$
onto
$$\cz_{[d_{n-2},d_{n-4}]}
\setminus \left(Z_{n-1}\cup Y_{n-1}\right).$$
It follows that $f^{S_n}$ maps some region in
$$\cz_{[z_{n-1},\hat z_{n-1}]}
\setminus Y_n $$
as an $\ell$-cover onto $Slit_{n-1}$.
Because $f^{-1}(f(Z_n))$ consists of
$\ell$ lines through the critical point,
this region has a slit in all of these $\ell$ lines.
Similarly, $f^{S_n}$ maps some region in
$$\cz_{[z_{n-1},\hat z_{n-1}]}
\setminus (Y_n \cup Z_n)$$
onto
$$Slit_{n-1}'=\cz_{[d_{n-2},\hat d_{n-2}]}
\setminus \left(Z_{n-1}\cup Y_{n-1}\right).$$
Since $[d_n,\hat d_n]\supset [z_{n-1},\hat z_{n-1}]
\supset [u_n,\hat u_n]$ the first part of the theorem follows.

The last part of the theorem
holds because if $F^k(A)=A_r$ or $F^k(A)=A_r'$
then $F^{k+1}$ covers
$$\cup_{i\ge r}A_i\cup A_{r-1}'\cup A_{r-2}\subset Slit^*_{r-2}.$$
Therefore the last assertion follows immediately
by induction from the first part of the theorem.
\qed

\begin{theo}\label{slitkoebe}
There exists a universal constant $C>0$ such that if we define
$Y_n$ and $Z_n$ as above
(so $Y_n\subset D_{n+1}^1 \subset A_n$),
then
$$\text{dist}\left((A_n\setminus D_{n+1}^1),Y_n\right)
\ge C|A_n|$$
and
$$\text{dist}\left(A_n,Z_n\right)
\ge C|A_n|$$
provided $\ell$ is sufficiently large and $n\ge n_0(\ell)$
is sufficiently large.
\end{theo}
\pr As we have seen in Theorem~\ref{complb},
$f^{S_n}$ maps $D_{n+1}^1$ as a univalent map onto
$D_n$. Moreover, $f^{S_n}(Y_n)=[d_{n+2},y_{n+2}]$.
Now because of the last part of Theorem~\ref{complb}
$D_n$ contains a neighbourhood of $[d_{n+2},y_{n+2}]$
with thickness is comparable to the size of this interval:
there exists $C>0$ such that each point in a neighbourhood
of $Y_n$ of the form
$$N_{C|d_{n+2}-y_{n+2}|}:=
\{x \st d\left(x,[d_{n+2},y_{n+2}]\right)\le C\cdot |d_{n+2}-y_{n+2}|\}$$
is contained in $D_n$.
In fact, this proves the last assertion of the
theorem (but with $n$ replaced by $n-1$).
Moreover, the Koebe Lemma
this gives that $f^{S_n}$ has uniformly bounded distortion
on the subset of $D_{n+1}^1$ (containing $Y_n$)
which is mapped diffeomorphically onto
the slightly thinner set
$$N_{(C/2)|d_{n+2}-y_{n+2}|}.$$
Since $f^{S_n}(Y_n)=[d_{n+2},y_{n+2}]$,
it follows that $Y_n$ also has a reasonably thick neighbourhood
inside $D_n^1$.
\qed
\bigskip

Let us now combine the results of this section in the
following way.
Define $A_n^i$ be
the part of $A_n$ which is between
the rays $l_i$ and $l_{i+1}$
where $l_i$ is given by
$\rz^+\ni t\mapsto t e^{2\pi i  /\ell}\in \cz$.
As explained below Theorem~\ref{complb},
$A_n^i$ consists of one component and
$f$ maps $A_n^i$ diffeomorphically onto
$f(A_n)$ (appart from the image of $l_i$).
\bigskip
\begin{theo}\label{maco}
Take $A\in \A_{k+1}$. Then there exists $r$ such that
$F^k$ maps $A$ diffeomorphically to $\tilde A_r$
where $\tilde A_r$ is either $A_r$ or to $A_r'$.
Moreover, let $A^i=F^{-k}(A_r^i)$. Then
$$F^k\colon A^i\to A_r^i$$
has uniformly bounded distortion.
\end{theo}
\pr
The diameter of $A_r^i$ is comparable to the
distance of $A_r^i$ to the nearest critical value of $F^k$.
Therefore the result follows from Koebe.
\qed

\sect{An asymptotic expression for the real induced map}
\label{asytex}

In this section we shall give a very good estimate
for the diffeomorphism
$f^{S_n-1}|[r_n^f,u_n^f]\to [\hat u_{n-2},u_{n-2}]=
\cup_{i\ge -2}(A_{n+i}\cap \rz)$. Here $$A_{n+i}\cap \rz=
[u_{n+i},\hat u_{n+i}]\setminus [u_{n+i+1},\hat u_{n+i+1}].$$
\bigskip

\begin{theo}\label{miracle}
There exists a constant $C$ and $\ell_0$
such that for each $\ell\ge \ell_0$ and
each $n\ge n_0(\ell)$,
the following estimates for the derivative of
the diffeomorphisms
$f^{S_n-1}|[r_n^f,u_n^f]\to [\hat u_{n-2},u_{n-2}]$
and
$f^{S_{n-2}}|[x_{n-1},u_{n-1}]\to [\hat u_{n-2},u_{n-2}]$
hold. Let $i\ge 2$ and take $x\in [r_n^f,u_n^f]$
and $y\in [x_{n-1},u_{n-1}]$ so that
$f^{S_n-1}(x)\in A_{n+i}$ and $f^{S_{n-2}}(y)\in A_{n+i}$.
Then
\beq
|Df^{S_n-1}(x)|\le C \frac{i^{3/2}}{\ell}\frac{|u_{n-2}-\hat
u_{n-2}|}{|r_n^f-u_n^f|}
\label{mi1}
\eeq
respectively
\beq
|Df^{S_{n-2}}(y)|\le C \frac{i^{3/2}}{\ell}\frac{|u_{n-2}-\hat
u_{n-2}|}{|x_{n-1}-u_{n-1}|}.
\label{mi2}
\eeq
\end{theo}
\medskip

The remainder of this section shall be occupied with the
proof of this theorem.
\bigskip

First we should point out that since $\cup_{i\ge n}
A_i$ has a definite
amount of Koebe space around $\cup_{i\ge n+\ell} A_i$, it follows that
the maps from the previous theorem
have uniformly bounded distortion on the
piece that maps into $A_{n+\ell}$. Hence we could replace
the $i^{3/2}/\ell$ term in (\ref{mi1}) and (\ref{mi2})
by $\min(i^{3/2}/\ell,\sqrt{\ell})$.

Secondly, we should emphasize that these estimates are far better than
would obtain from a Koebe estimate: using Koebe we would
have to replace $i^{3/2}$ by $i^2$ (for $i\le \ell$).
As we shall show
below these `good' estimates hold because -- by the
real bounds -- the maps $f^{S_{n-2}}$ are far away from
Moebius transformations. We shall come back to this
issue below. In fact, even though we shall prove that
(\ref{mi2}) holds for $i=\ell$ by a simple
cross-ratio argument, to show that (\ref{mi2}) holds
for all $i\ge 2$ we shall use renormalization in a crucial way.

\kies{

\begin{figure}[htp]
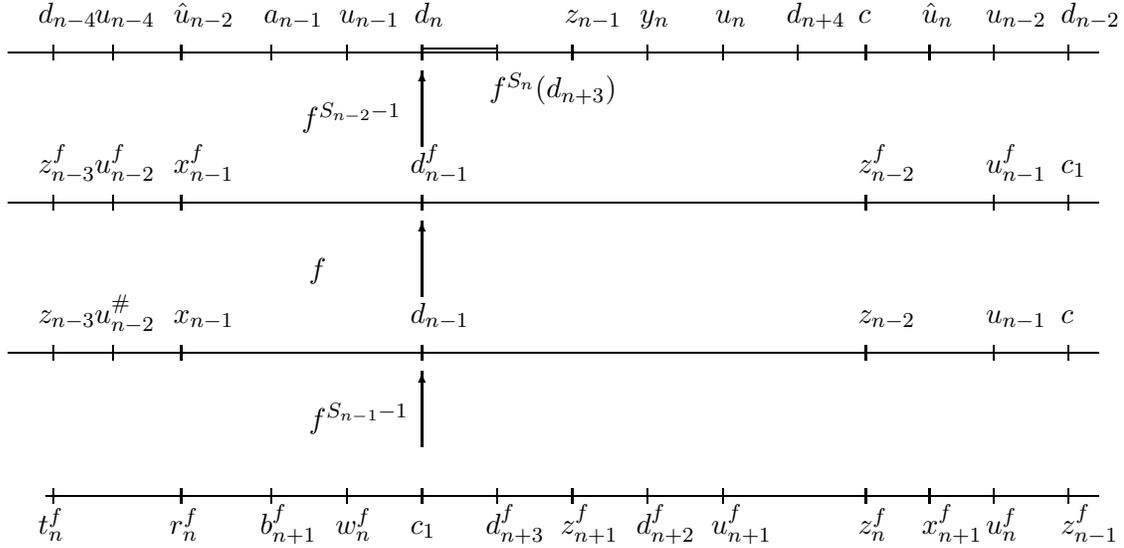

\vskip 0.7cm
\hbox to \hsize{\hss\unitlength=5mm
\beginpic(20,12)(-20,0) \let\ts\textstyle
\put(4,12){\line(-1,0){29}}
\put(3.2,11.8){\line(0,1){0.4}} \put(3,12.8){\small $d_{n-2}$}
\put(1.2,11.8){\line(0,1){0.4}} \put(1,12.8){\small $u_{n-2}$}
\put(-0.5,11.8){\line(0,1){0.4}} \put(-0.7,12.8){\small $\hat u_n$}
\put(-2.2,11.8){\line(0,1){0.4}} \put(-2.4,12.8){\small $c$}
\put(-4,11.8){\line(0,1){0.4}} \put(-4.3,12.8){\small $d_{n+4}$}
\put(-6,11.8){\line(0,1){0.4}} \put(-6.2,12.8){\small $u_n$}
\put(-8,11.8){\line(0,1){0.4}} \put(-8.2,12.8){\small $y_n$}
\put(-10,11.8){\line(0,1){0.4}} \put(-10.2,12.8){\small $z_{n-1}$}
\put(-12,11.8){\line(0,1){0.4}}\put(-12.2,10.8){\small $f^{S_n}(d_{n+3})$}
\put(-14,12.1){\line(1,0){2}}
\put(-14,11.8){\line(0,1){0.4}} \put(-14.2,12.8){\small $d_n$}
\put(-16,11.8){\line(0,1){0.4}} \put(-16.2,12.8){\small $u_{n-1}$}
\put(-18,11.8){\line(0,1){0.4}} \put(-18.2,12.8){\small $a_{n-1}$}
\put(-20.4,11.8){\line(0,1){0.4}} \put(-20.6,12.8){\small $\hat u_{n-2}$}
\put(-22.2,11.8){\line(0,1){0.4}} \put(-22.7,12.8){\small $u_{n-4}$}
\put(-23.8,11.8){\line(0,1){0.4}} \put(-24.2,12.8){\small $d_{n-4}$}

\put(4,8){\line(-1,0){29}}
\put(3.2,7.8){\line(0,1){0.4}} \put(3,8.8){\small $c_1$}
\put(1.2,7.8){\line(0,1){0.4}} \put(1,8.8){\small $u_{n-1}^f$}
\put(-2.2,7.8){\line(0,1){0.4}} \put(-2.4,8.8){\small $z_{n-2}^f$}
\put(-14,7.8){\line(0,1){0.4}} \put(-14.3,8.8){\small $d_{n-1}^f$}
\put(-20.4,7.8){\line(0,1){0.4}} \put(-20.6,8.8){\small $x_{n-1}^f$}
\put(-22.2,7.8){\line(0,1){0.4}} \put(-22.7,8.8){\small $u_{n-2}^f$}
\put(-23.8,7.8){\line(0,1){0.4}} \put(-24.2,8.8){\small $z_{n-3}^f$}
\put(-14,9.5){\vector(0,1){2}}
\put(-17.2,10){\small $f^{S_{n-2}-1}$}

\put(4,4){\line(-1,0){29}}
\put(3.2,3.8){\line(0,1){0.4}} \put(3,4.8){\small $c$}
\put(1.2,3.8){\line(0,1){0.4}} \put(1,4.8){\small $u_{n-1}$}
\put(-2.2,3.8){\line(0,1){0.4}} \put(-2.4,4.8){\small $z_{n-2}$}
\put(-14,3.8){\line(0,1){0.4}} \put(-14.3,4.8){\small $d_{n-1}$}
\put(-20.4,3.8){\line(0,1){0.4}} \put(-20.6,4.8){\small $x_{n-1}$}
\put(-22.2,3.8){\line(0,1){0.4}} \put(-22.7,4.8){\small $u_{n-2}^\#$}
\put(-23.8,3.8){\line(0,1){0.4}} \put(-24.2,4.8){\small $z_{n-3}$}
\put(-14,5.5){\vector(0,1){2}}
\put(-17,6){\small $f$}

\put(4,0.2){\line(-1,0){28}}
\put(3.2,0){\line(0,1){0.4}} \put(3,-0.8){\small $z_{n-1}^f$}
\put(1.2,0){\line(0,1){0.4}} \put(1,-0.8){\small $u_{n}^f$}
\put(-0.5,0){\line(0,1){0.4}} \put(-0.7,-0.8){\small $x_{n+1}^f$}
\put(-2.2,0){\line(0,1){0.4}} \put(-2.4,-0.8){\small $z_n^f$}
\put(-6,0){\line(0,1){0.4}} \put(-6.3,-0.8){\small $u_{n+1}^f$}
\put(-8,0){\line(0,1){0.4}} \put(-8.3,-0.8){\small $d_{n+2}^f$}
\put(-10,0){\line(0,1){0.4}} \put(-10.3,-0.8){\small $z_{n+1}^f$}
\put(-12,0){\line(0,1){0.4}} \put(-12.3,-0.8){\small $d_{n+3}^f$}
\put(-14,0){\line(0,1){0.4}} \put(-14.3,-0.8){\small $c_1$}
\put(-16,0){\line(0,1){0.4}} \put(-16.3,-0.8){\small $w_n^f$}
\put(-18,0){\line(0,1){0.4}} \put(-18.3,-0.8){\small $b_{n+1}^f$}
\put(-20.4,0){\line(0,1){0.4}} \put(-20.7,-0.8){\small $r_n^f$}
\put(-23.8,0){\line(0,1){0.4}} \put(-24.2,-0.8){\small $t_n^f$}
\put(-14,1.5){\vector(0,1){2}}
\put(-17,2){\small $f^{S_{n-1}-1}$}
\endpic\hss}
\vskip 0.4cm
\caption[ ]{{\small Points and their images under $f^{S_{n-1}-1}$ and
$f^{S_n-2}$. The interval near $c_1$ shown in the second lowest line
should be thought of as ordered in the opposite direction
(since $c_1$ is the mimumum of $f$).
Here $u_{n-2}^\#$ stands for either $u_{n-2}$ or $\hat u_{n-2}$
depending on whether $n$ is even or odd.}}
\end{figure}

}

Thirdly, we should note that the map
$f^{S_{n-1}-1}$ restricted to $[r_n^f,u_n^f]$ has uniformly
bounded distortion because the image of this interval
is equal to $[x_{n-1},u_{n-1}]$ and because this
map extends diffeomorphically to an interval with
image $[d_{n-5},c]$. So by the
real Koebe Principle, see Proposition 3.1,
and the real bounds, see Theorem 3.1,
it follows that the distortion of $f^{S_{n-1}-1}|[r_n^f,u_n^f]$
is uniformly bounded (for all $\ell \ge 4$ and all
$n$ sufficiently large). So (\ref{mi1}) follows from (\ref{mi2}).

Hence we shall be interested in the distortion of the map
$f^{S_{n-2}}\colon [x_{n-1},u_{n-1}]\to [\hat u_{n-2},u_{n-2}]$.
This map only extends diffeomorphically to an interval
whose image extends on each side a fraction $1/\ell$.
So from the real Koebe Principle we can deduce that
the distortion of the map is bounded by $\ell^2$.
In this section we shall improve this bound.

Let $T_{n-1}=[x_{n-1},u_{n-1}]$, then $f^{S_{n-2}}(T_{n-1})=[\hat u_{n-2},u_{n-2}]$.
We first show that
\beq
|Df^{S_{n-2}}(u_{n-1})|=\frac{C_n}{\ell}\frac{|f^{S_{n-2}}(T_{n-1})|}{|T_{n-1}|},
\label{asy1}
\eeq
\beq
|Df^{S_{n-2}}(x_{n-1})|=\frac{C_n'}{\ell}
\frac{|f^{S_{n-2}}(T_{n-1})|}{|T_{n-1}|}
\label{asy2}
\eeq
where $C_n,C_n'>0$ are uniformly bounded and bounded away from zero.
This follows easily from the real Koebe Principle and the real bounds:
all the intervals connecting the points
$u_{n-2}^\#$, $x_{n-1}$, $d_{n-1}$, $z_{n-2}$ and $u_{n-1}$ are of the
same order (and $c$ is `far away').
But the interval $[u_{n-2}^\#,x_{n-1}]$ is mapped diffeomorphically
to $[u_{n-4},\hat u_{n-2}]$ whose size is order $\frac{1}{\ell}$ times
$|f^{S_{n-2}}(T)|$.
Because there is Koebe space around $[u_{n-4},\hat u_{n-2}]$,
formula (\ref{asy1}) follows. Similarly, (\ref{asy2}) follows by considering
the interval $[f^{S_{n-1}-1}(x_{n+1}^f),u_{n-1}]$ (because
the size of this interval is also of the same order as $|T|$
and its image under $f^{S_{n-2}}$ has a size of the order
$(1/\ell)|f^{S_{n-2}}(T_{n-1})|$). Since one has again
Koebe space around this interval, (\ref{asy2}) follows.

Since the Koebe space around the interval
$f^{S_{n-2}}(T_{n-1})$ is only of order $1/\ell$, we get from the
real Koebe Principle that
$$|Df^{S_{n-2}}(x)|\le C\ell^2 |Df^{S_{n-2}}(u_{n-1})|
\frac{|f^{S_{n-2}}(T_{n-1})|}{|T_{n-1}|}
\le C\ell \frac{|f^{S_{n-2}}(T_{n-1})|}{|T_{n-1}|}.$$
In this section we shall show that this estimate is far from optimal.
First we shall give this estimate in a special case
using cross-ratios only; this proof does not require renormalization.
Since it is not strictly needed in this paper, the reader
can skip the next subsection.

\subsection{A partial result using cross-ratios}

In this subsection we shall show

\begin{prop}
\label{rootl}
$$|Df^{S_{n-2}}(z_{n-2})| \le C\sqrt{\ell}\frac{|f^{S_{n-2}}(T_{n-1})|}{|T_{n-1}|}.$$
\end{prop}

In fact, in Proposition~\ref{comdv}, we shall also obtain
estimates for $|Df^{S_{n-2}}(\gamma)|$ when $\gamma$ is some
arbitrary point in $T$. Those estimates imply
Proposition~\ref{rootl} but unlike the proof
of Proposition~\ref{rootl} are based on renormalization results.
In fact, we feel that the proof of Proposition~\ref{rootl}
should have much wider applications: philosophically
speaking the proof shows that if a map is not too close to
a Moebius transformation then one has improved Koebe estimates!!
For this reason we have included Proposition~\ref{rootl}.

In the proof of Proposition~\ref{rootl} we use the following
cross-ratio operator: if $J=[\alpha,\beta]$
is an interval and $g\colon J\to \rz$
a diffeomorphism, define
$$A(g,J)=\frac{[\frac{|g(J)|}{|J|}]^2}{|Dg(\alpha)||Dg(\beta)|}.$$

\noindent
{\em Proof of Propositon~\ref{rootl}:}
Consider $T'=[x_{n-1},z_{n-2}]$ and $T''=[z_{n-2},u_{n-1}]$.
These intervals have one point in common and their union is
equal to $T=T_{n-1}$.
Because of the next lemma and because $A(g_1\circ g_2,J)=A(g_1,g_2(J))\cdot
A(g_2,J)$ one has
$$A(f^{S_{n-2}},T')A(f^{S_{n-2}},T'')\ge \sqrt{A(f^{S_{n-2}},T)}.$$
By (\ref{asy1}) and (\ref{asy2}),
$A(f^{S_{n-2}},T)$ is of order $\ell^2$. Hence
\beq
A(f^{S_{n-2}},T')A(f^{S_{n-2}},T'')\ge C \ell.
\label{asy3}
\eeq
Now $|f^{S_{n-2}}(T_i)|/|T_i|$, $i=1,2$ are both of the
same order as $|f^{S_{n-2}}(T)|/|T|$ because by the real bounds
$T_i,T$ and also $|f^{S_{n-2}}(T_i)|, |f^{S_{n-2}}(T)|$
are of the same order. Using this and (\ref{asy1}), (\ref{asy2}),
$$A(f^{S_{n-2}},T')A(f^{S_{n-2}},T'')$$
is of the same order as
$$\left[\ell\frac{\frac{|f^{S_{n-2}}(T)|}{|T|}}{|Df^{S_{n-2}}(z_{n-2})|}
\right]^2.$$
Using (\ref{asy3}) this gives that
$$
|Df^{S_{n-2}}(z_{n-2})|\le
\sqrt{\ell}\frac{|f^{S_{n-2}}(T)|}{|T|}.
$$
\qed

In the proof of the previous proposition we used the following
lemma:

\begin{lemma}
Assume that $f(z)=z^\ell+c_1$. Then there exists $\ell_0$
such that for each $\ell\ge \ell_0$ and each two intervals
$T',T''$ in one component of
$\rz\setminus \{0\}$ with a unique common endpoint,
\beq
A(f,T')A(f,T'')\ge \sqrt{A(f,T'\cup T'')}
\label{asy4}
\eeq
\end{lemma}

\noindent
{\em Proof of the lemma:} Write $T'=[\alpha,x]$,
$T''=[x,\beta]$ and assume for simplicity
that $\alpha<x<\beta$. First we {\bf claim} that the left hand side of
(\ref{asy4}) is minimal if $x=\sqrt{\alpha \beta}$. To prove this,
first notice that this expression is invariant if we multiply
all the points by the same factor. Hence we may assume that $\alpha=1$
and write $x=e^{\tau b}$ and $\beta=e^b$ where $\tau\in (0,1)$
and $b>0$. With this notation the left hand side of (\ref{asy4})
becomes
$$\frac{(\frac{e^{\ell b \tau}-1}{e^{b \tau}-1})^2
(\frac{e^{\ell b }-e^{\ell b \tau}}{e^{b \tau}-e^{b\tau}})^2}
{\ell \cdot \ell^2 e^{2(\ell-1)b\tau}\cdot \ell \cdot e^{(\ell-1)b}}
=
\frac{(\frac{e^{\ell b \tau}-1}{e^{b \tau}-1})^2
(\frac{e^{\ell b (1-\tau)}-1}{e^{b (1-\tau)}-1)^2}}
{\ell^4 e^{(\ell-1)b}}.$$
Now the denominator of this term,
$$(\frac{e^{\ell b \tau}-1}{e^{b \tau}-1})^2
(\frac{e^{\ell b (1-\tau)}-1}{e^{b (1-\tau)}-1})^2
$$
is equal to
$$
G(\tau):=(1+e^{b\tau}+\dots+e^{(\ell-1)b\tau})
(1+e^{b(1-\tau)}+\dots+e^{(\ell-1)b(1-\tau)})
=
\mbox{ const }+ P(e^{b\tau})+P(e^{b(1-\tau)})$$
where $P$ is a polynomial with positive constants.
Hence
$$G'(\tau)=P'(e^{b\tau})\cdot e^{b\tau} - P'(e^{b(1-\tau)})\cdot
e^{b(1-\tau)}.$$
Since $P$ has positive coefficients,
$\tau\mapsto P'(e^{b\tau})\cdot e^{b\tau}:=G_1(\tau)$
is increasing,
$G_1(\tau)=G_1(1-\tau)$ holds only if $\tau=1/2$.
Since $G'(\tau)=0$ is equivalent to
$G_1(\tau)=G_1(1-\tau)$ the first claim holds.

The previous claim implies that it suffices to consider
the situation that $\alpha=1$, $x=e^a$ and $\beta=e^{2a}$.
So we need to prove that
$$A(f,T')A(f,T'')=
\frac{(\frac{e^{\ell a}-1}{e^{a}-1})^2
(\frac{e^{\ell 2a }-e^{\ell a}}{e^{2a}-e^{a}})^2}
{\ell \cdot \ell^2 e^{2(\ell-1)a}\cdot \ell \cdot e^{(\ell-1)2a}}
=
\frac{(\frac{e^{\ell a}-1}{e^{a}-1})^4}{\ell^4 e^{2(\ell-1)a}}$$
is greater or equal than
$$\sqrt{A(f,T)}=
\frac{\frac{e^{2\ell a}-1}{e^{2a}-1}}{\ell e^{(\ell-1)a}}.$$
Hence, writing
$$g_1(t)=(e^t-1)^3\mbox{ and }g_2(t)=(e^t+1)e^t t^3/2,
$$
and using $e^{2t}-1=(a^t-1)(e^t+1)$ the required estimate
(\ref{asy4})
is equivalent to
\beq
g_1(\ell a) g_2(a) \ge g_2(\ell a) g_1(a).
\label{asy4a}
\eeq
To show that there exists $\ell_0$ such that
this inequality holds for each $a>0$ and each $\ell\ge \ell_0$
we proceed as follows. First notice that
\beq
g_1(t)-g_2(t)=(e^t-1)^3-(e^t+1)e^t t^3/2=
(e^{3t}-3e^{2t}+3e^t-1)-(e^{2t}+e^t)t^3/2.
\label{asy5}
\eeq
The coeffient corresponding to the $t^n$-th term in
the Taylor expansion of (\ref{asy5}) is equal to
$$\frac{1}{n!}[(3^n-3 \, 2^n +3) - \frac{n(n-1)(n-2)}{2}(2^{n-3}+1)].$$
So these coefficients are zero for $n=2,3,4,5,6$ and strictly positive
for $n\ge 7$. In fact, the coefficient corresponding to $t^7$ is equal
to $1/240$ and so we get
$$g_1(t)-g_2(t)\ge \frac{1}{240}t^7$$
for $t\ge 0$.
Hence
\beq
\frac{g_1(t)}{g_2(t)}-1\ge \frac{\frac{1}{240}t^7}{\frac{(e^t+1)e^t}{2}t^3}
\ge \frac{1}{120}\frac{t^4}{e^{2t}+e^t}.
\label{asy6}
\eeq
Moreover, there exists $t_0$ such that
\beq
\frac{g_1(t)}{g_2(t)}-1=
\frac{(e^t-1)^3-t^3e^{2t}/2-t^3e^t/2}{\frac{(e^t+1)e^t}{2}t^3}
\ge e^{(1/2)t}t^4
\label{asy7}
\eeq
for all $t\ge t_0$. Combining (\ref{asy6}) and (\ref{asy7})
we get that
\beq
\frac{g_1(t)}{g_2(t)}-1 \ge \mbox{ const } \cdot e^{(1/2)t}t^4
\label{asy8}
\eeq
for all $t\ge 0$ where $\mbox{const}$ is a positive constant.
Similarly one gets that
\beq
\frac{g_1(t)}{g_2(t)}-1 \le \mbox{ const' } \cdot  e^{2t}t^4
\label{asy9}
\eeq
for all $t\ge 0$ where $\mbox{ const' }$ is a finite constant.
Applying (\ref{asy8}) and (\ref{asy9}) it follows that
$$\frac{g_1(\ell a)}{g_2(\ell a)}
\frac{g_2(a)}{g_1(a)}
\ge
\frac
{1+\text{ const }\cdot (\ell a)^4 e^{(1/2)\ell a}}
{1+\text{ const' } \cdot a^4 e^{2 a}}
\ge (1+\text{ const'' } a^4)\ge 1$$
for each $a\ge 0$
provided $\ell$ is sufficiently large.
This concludes the proof of (\ref{asy4a}) and the proof of the
lemma.
\qed

We would also like to remark that (\ref{asy4}) and (\ref{asy4a})
also hold for the quadratic case $\ell=2$.
Indeed, in this case (\ref{asy4a})
is equivalent to showing that for all $a\ge 0$,
$$e^{4a}+6e^{2a}+1\ge 4 e^{3a}+4 e^a.$$
So it suffices to show that
the coefficients of the power series of the left hand side dominates those
of the right hand side. This means that we have to show that
$$4^n+6\, 2^n\ge 4\, 3^n+4$$
for each $n\ge 1$. This is readily checked.
Presumably, this lemma holds for all $\ell\ge 2$.

Moreover, the same ideas also work for $C^2$ maps
because in this case there exists a universal constant such that
$$
A(f,T')A(f,T'')\ge \sqrt{A(f,T'\cup T'')}(1-C|T'\cup T''|).
$$
So if we have that $\sum_i |f^i(T'\cup T'')|$ is bounded
then one can proceed as before. In particular, such a bound
holds for the Fibonacci map (see Section 2 of \cite{BKNS})
Proposition~\ref{rootl} also holds for $C^2$ Fibonacci maps with
a critical point of order $\ell$.
\bigskip

\subsection{An asymptotic expression for $f^{S_{n-2}}$ on
$[x_{n-1},u_{n-1}]$}

Now we will give a more precise
version of the last proposition. In fact, we will
obtain an asymptotic expression for large $\ell$
of the limit of the sequence of diffeomorphism
$$f^{S_{n-2}}\colon T_{n-1}\to [u_{n-2},\hat u_{n-2}],$$
$n\in 2\N$
where as before $T_{n-1}=[x_{n-1},u_{n-1}]\ni z_{n-2}$.
Note that this situation is quite remarkable: as
$\ell$ increases the non-linearity increases and the amount of
Koebe space decreases. Evenso, we are able
to determine the limit function (it is far from linear).
So let
$$H_{n-1}\colon T_{n-1}\to \rz$$ be the orientation
preserving affine map with
$$H_{n-1}(z_{n-2})=0$$
and such that
$$|H_{n-1}(u_{n-1})|=1.$$
Moreover, let $L_n\colon [-1,1]\to [u_{n-2},\hat u_{n-2}]$
be the linear orientation preserving bijection.
\bigskip

\begin{theo}
\label{miracle2}
There exists $K(\ell)>0$ which is universally bounded and bounded away
from zero such that if we write
$$\Gamma(x)=\frac{x \sqrt{K(\ell)\ell+1}}{\sqrt{K(\ell) \ell x^2 +1}}$$
then for $n>\ell^{3/2}$, $n\in 2\nz$
and for each $x\in [x_{n-1},u_{n-1}]$,
$$ L_n\circ  \Gamma\circ H_{n-1}(x)\cdot
(1-o(1/\ell)) \le
f^{S_{n-2}}(x) \le L_n\circ \Gamma\circ H_{n-1}(x)\cdot
(1+o(1/\ell))$$
and
$$
D\left(L_n \circ \Gamma \circ H_{n-1}\right)(x)\cdot (1-o(\ell)) \le
Df^{S_{n-2}}(x) \le
D\left(L_n\circ \Gamma \circ H_{n-1}\right)(x)\cdot (1+o(\ell)).$$
Here $o(t)$ stands for a function which tends to zero
as $t$ tends to zero.
For $n\in 2\nz+1$ there exists a similar constant $K(\ell)$.
\end{theo}

\begin{remark} We should note that
$\Gamma$ can be written as a composition
of the square map $x\mapsto x^2$,
the Moebius map
$$M_\ell\colon t\mapsto \frac{t(K(\ell)\ell +1)}{K\ell t +1}$$
and the root map $x\mapsto \sqrt{x}$.
This Moebius map $M_\ell$ send the interval $[0,1]$ onto itself
but pushing points extremely far to the right when $\ell$ is large.
The good bounds for the distortion of
$f^{S_{n-2}}$ come from the fact that it is close
to $M_\ell$ {\it up to a conjugation with the square map}.
Presumably, using Ecalle cylinders or so, one can
also give good estimates on the domain in the complex plane
for which this asymptotic expression holds.
\end{remark}

The idea of this result is
related to the so-called Ecalle-cylinders
used in \cite{Sh1}.
Before proving this theorem let us show how Theorem~\ref{miracle}
follows from it.

\bigskip

\noindent
{\em Proof of Theorem~\ref{miracle}:}
Note that $|u_{n-2}|=|\hat u_{n-2}|$ and the symmetry
of the map $\Gamma$ in the previous theorem implies that
$1-o(1/\ell)\le |H_{n-1}(x_{n-1})|\le 1+o(1/\ell)$
for large $n$. In other words,
$$H_{n-1}[x_{n-1},u_{n-1}]\to [-1\pm o(1/\ell),1 \pm o(1/\ell)].$$
We need to determine the preimages of the intervals
$A_{n+i}\cap \rz$ under $f^{S_{n-2}}\colon T_{n-1}=[x_{n-1},u_{n-1}]
\to [u_{n-2},\hat u_{n-2}]$ and then determine the size of
the derivative of this map in these intervals.
First note that the intervals $A_{n+j}\cap \rz$
have a size of the order $1/\ell$ of the size of $[u_{n-2},\hat u_{n-2}]$
when $-2\le j\le \ell$. So $H_{n-1}^{-1}(A_{n+j})$
is of the
form
$$\pm [(1-C(j+3)/\ell),(1-C (j+2)/\ell)].$$
for $-2\le j\le \ell$.
The inverse of $\Gamma$ is
$$\Gamma^{-1}(y)=\frac{y}{\sqrt{K \ell(1-y^2)+1}}.$$
This implies that the preimage
of $\pm (1-C(j+2)/\ell)$ is of the form
$\frac{C_0}{\sqrt{j+3}}$, $-2\le j\le \ell$.
So if $y$ is in the annulus $A_{n+j}$
then its preimage
$\tilde x$ under $H_{n-1}\circ \Gamma$
is in an interval of the form $[\frac{C_0}{\sqrt{j+4}},
\frac{C_0}{\sqrt{j+3}}]$.
Since
$$D\Gamma(\tilde x)=\frac{\sqrt{K\ell+1}}{(K\ell \tilde x^2+1)^{3/2}}
\le C\frac{\sqrt{\ell}}{(K\ell\frac{1}{j}+1)^{3/2}}
\le C\frac{\sqrt{\ell}}{(\ell+j)^{3/2}}j^{3/2}
\le C\frac{j^{3/2}}{\ell}.$$
This gives that if $x$ is the preimage of
$y$ under $f^{S_{n-2}}\colon [x_{n-1},u_{n-1}]\to [u_{n-2},\hat u_{n-2}]$
then $x$ is also in scaled down interval of this form
and the previous estimate gives
$$|Df^{S_{n-2}}(x)|\le C(1+o(1/\ell))\frac{j^{3/2}}{\ell}
\frac{|u_{n-2}-\hat u_{n-2}|}{|x_{n-1}-u_{n-1}|}.$$
This is the required estimate for $j<\ell$.
If $y\in A_{n+j}$ with $j\ge \ell$, then for its
preimage $x$ one has $|x|\le \frac{C}{\sqrt{\ell}}$ and therefore
we get also
$$D \Gamma(x)\le C\sqrt{\ell}.$$
\qed

Now we will start with the proof of  Theorem~\ref{miracle2}.
Notice that $f^{S_{n-2}}$ can be written
as
$$f^{S_{i_0}}\circ f^{S_{i_0+1}}\circ f^{S_{i_0+3}}\circ \dots \circ
f^{S_{n-7}}\circ f^{S_{n-5}}\circ f^{S_{n-3}}$$
where $i_0<n$ and $n-i_0$ is even.

\kies{

\begin{figure}[htp]
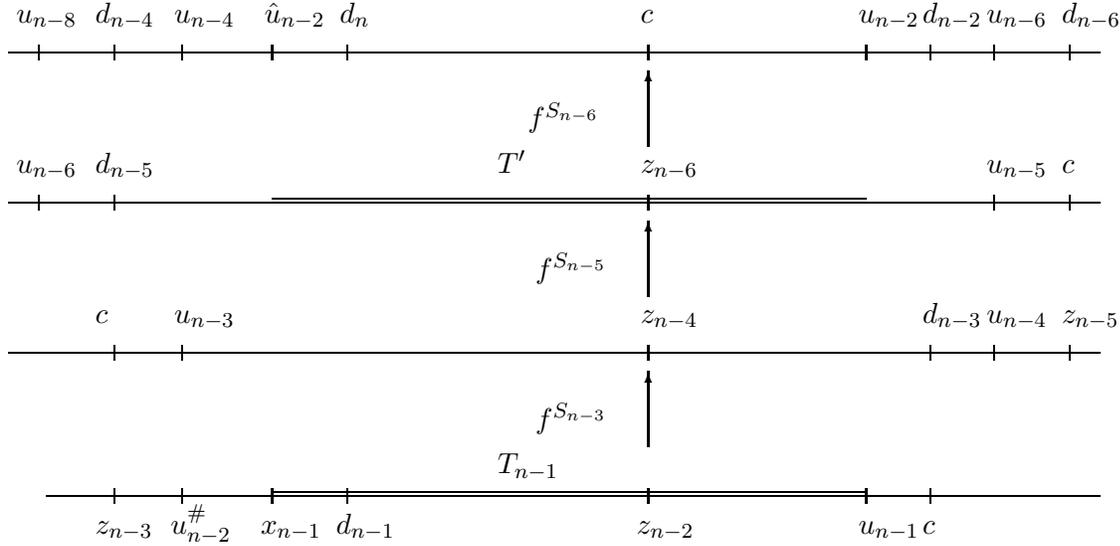

\vskip 0.7cm
\hbox to \hsize{\hss\unitlength=5mm
\beginpic(20,12)(-20,0) \let\ts\textstyle
\put(4,12){\line(-1,0){29}}
\put(3.2,11.8){\line(0,1){0.4}} \put(3,12.8){\small $d_{n-6}$}
\put(1.2,11.8){\line(0,1){0.4}} \put(1,12.8){\small $u_{n-6}$}
\put(-0.5,11.8){\line(0,1){0.4}} \put(-0.7,12.8){\small $d_{n-2}$}
\put(-2.2,11.8){\line(0,1){0.4}} \put(-2.4,12.8){\small $u_{n-2}$}
\put(-8,11.8){\line(0,1){0.4}} \put(-8.2,12.8){\small $c$}
\put(-16,11.8){\line(0,1){0.4}} \put(-16.2,12.8){\small $d_n$}
\put(-18,11.8){\line(0,1){0.4}} \put(-18.2,12.8){\small $\hat u_{n-2}$}
\put(-20.4,11.8){\line(0,1){0.4}} \put(-20.6,12.8){\small $u_{n-4}$}
\put(-22.2,11.8){\line(0,1){0.4}} \put(-22.7,12.8){\small $d_{n-4}$}
\put(-24.2,11.8){\line(0,1){0.4}} \put(-24.8,12.8){\small $u_{n-8}$}

\put(4,8){\line(-1,0){29}}
\put(3.2,7.8){\line(0,1){0.4}} \put(3,8.8){\small $c$}
\put(1.2,7.8){\line(0,1){0.4}} \put(1,8.8){\small $u_{n-5}$}
\put(-8,7.8){\line(0,1){0.4}} \put(-8.2,8.8){\small $z_{n-6}$}
\put(-22.2,7.8){\line(0,1){0.4}} \put(-22.7,8.8){\small $d_{n-5}$}
\put(-24.2,7.8){\line(0,1){0.4}} \put(-24.8,8.8){\small $u_{n-6}$}
\put(-18,8.1){\line(1,0){15.8}}
\put(-12,8.8){\small $T'$}
\put(-8,9.5){\vector(0,1){2}}
\put(-11.2,10){\small $f^{S_{n-6}}$}

\put(4,4){\line(-1,0){29}}
\put(3.2,3.8){\line(0,1){0.4}} \put(3,4.8){\small $z_{n-5}$}
\put(1.2,3.8){\line(0,1){0.4}} \put(1,4.8){\small $u_{n-4}$}
\put(-0.5,3.8){\line(0,1){0.4}} \put(-0.7,4.8){\small $d_{n-3}$}
\put(-8,3.8){\line(0,1){0.4}} \put(-8.2,4.8){\small $z_{n-4}$}
\put(-20.4,3.8){\line(0,1){0.4}} \put(-20.6,4.8){\small $u_{n-3}$}
\put(-22.2,3.8){\line(0,1){0.4}} \put(-22.7,4.8){\small $c$}
\put(-8,5.5){\vector(0,1){2}}
\put(-11,6){\small $f^{S_{n-5}}$}

\put(4,0.2){\line(-1,0){28}}
\put(-18,0.3){\line(1,0){15.8}}
\put(-12,0.8){\small $T_{n-1}$}
\put(-0.5,0){\line(0,1){0.4}} \put(-0.7,-0.8){\small $c$}
\put(-2.2,0){\line(0,1){0.4}} \put(-2.4,-0.8){\small $u_{n-1}$}
\put(-8,0){\line(0,1){0.4}} \put(-8.3,-0.8){\small $z_{n-2}$}
\put(-16,0){\line(0,1){0.4}} \put(-16.3,-0.8){\small $d_{n-1}$}
\put(-18,0){\line(0,1){0.4}} \put(-18.3,-0.8){\small $x_{n-1}$}
\put(-20.4,0){\line(0,1){0.4}} \put(-20.7,-0.8){\small $u_{n-2}^{\#}$}
\put(-22.2,0){\line(0,1){0.4}} \put(-22.7,-0.8){\small $z_{n-3}$}
\put(-8,1.5){\vector(0,1){2}}
\put(-11,2){\small $f^{S_{n-3}}$}
\endpic\hss}
\vskip 0.4cm
\caption[ ]{{\small The map $f^{S_{n-2}}|T$ can be factored as shown
(or as a longer composition).}}
\end{figure}

}

Our aim is to give an
asymptotic expression for this composition
by comparing it with the solution of some particular
differential equation which has an almost neutral attracting
singularity. To do this
we preceed similarly as in Lemma~\ref{psin} and Lemma~\ref{psin2}.
Note that $f^{S_i}$ maps $z_{i+1}$ to $z_{i-1}$
and $[u_{i+2},x_{i+2}]\ni z_{i+1}$ into $[u_i,x_i]\ni z_{i-1}$.
So let $H_i\colon [u_i,x_i]\to \rz$
be the orientation preserving affine map
with
$$H_i(z_{i-1})=0$$
so that
$$|H(u_i)|=1.$$
By the real bouds, $|H_i(x_i)|$
is also uniformly bounded and bounded away from zero.
Let
$$\Psi_i=H_i\circ f^{S_i}\circ H_{i+2}^{-1}\colon \rz\to \rz.$$
As in the proof of Lemma~\ref{psin2},
$$\Psi_i(z)=(-1+\alpha_i'/\ell)z+\beta_i' z^2 + \gamma_i' z^3 + O(z^4)$$
where the coefficients $\alpha_i',\beta_i',\gamma_i'$ do depend
on $\ell$ but for each fixed $\ell$ converge to constants
as $i\to \infty$ with either $i\in 2\nz$ or $i\in 2\nz+1$. 
In fact, $\liminf \alpha_i'$
is uniformly bounded away from $0$.
Hence, in the same way as in Lemma~\ref{psin2},
$$\Theta_i:=\Psi_{i-2}\circ \Psi_i(z)
=(1-\frac{\alpha_i}{\ell})z+\frac{\beta_i}{\ell}z^2
-\gamma_i z^3 + O(z^4)$$
where $\alpha_i,\gamma_i$ converge to positive constants (again as
$i\in 2\nz$ or $i\in 2\nz+1$). The expression $O(z^4)$
stands for a function which in norm is dominated by
$C|z|^4$; that this last bound holds is explained below (\ref{bcoe}).
These limits are uniformly bounded and bounded away for all large $\ell$
and $\alpha_i,\beta_i,\gamma_i$ converge to constants which are
uniformly bounded in norm.

Because of the convergence
of the sequence of renormalizations, as $i\in 2\nz$, tends to infinity,
$[H_i(u_i),H_i(x_i)]$ tends to an interval $M$ (having $1$ as an
endpoint and containing $0$ in its interior) and -- up to scaling --
$f^{S_{i-2}}\circ f^{S_i}$ also converges
(provided $i$ runs through either the even or
the odd integers).
In particular, $\Psi_{i-2}\circ \Psi_i$
tends to some fixed map $\Theta\colon M\to M$
of the form
$$\Theta
=(1-\frac{\alpha_{0,\ell}}{\ell})z+\frac{\beta_{0,\ell}}{\ell}z^2
-\gamma_{0,\ell} z^3 + O(z^4)$$
on a neighbourhood of $M\subset \cz$ where
the coefficients $\alpha_{0,\ell},\beta_{0,\ell},
\gamma_{0,\gamma}$. In fact, this limit
depends on whether $i$ runs through
$i\in 2\nz$ or $i\in 2\nz+1$, but to avoid needless
repetition we will not mention this anymore in the
remainder of this section.
Of course $\Theta$ presumably depends on $\ell$.
So there exists $i_0$ (depending on $\ell$) so that
for all $i\ge i_0$ the maps
$$f^{S_{i-2}}\circ f^{S_i}\colon
[x_{i+1},u_{i+1}] \to [x_{i-3},u_{i-3}]$$
are -- up to scaling -- in a given
neighbourhood of the limiting map $\Theta$.

Choose $\hat \iota\in \{i_0,i_0+1,i_0+2,i_0+3\}$
so that $n-\hat \iota\in 4\nz$.
Then
$$f^{S_{n-2}}=f^{S_{\hat \iota-2}}\circ f^{S_{\hat  \iota -1}}\circ f^{S_{\hat 
\iota+1}}
\circ \dots \circ
f^{S_{n-7}}\circ f^{S_{n-5}}\circ f^{S_{n-3}}$$
is a composition of a large number of maps of the form
$f^{S_{i-2}}\circ f^{S_i}$ with $i\ge i_0$ and the map $f^{S_{i_0}}$.
The image of $f^{S_{n-2}}$ of the
interval $T_{n-1}=[x_{n-1},u_{n-1}]$ is equal to $[\hat u_{n-2},u_{n-2}]$.
Let
$$T'=f^{S_{\hat \iota-1}}\circ \dots \circ f^{S_{n-5}}\circ f^{S_{n-3}}
(T_{n-1})$$
and let $H$ be the maximal interval containing $T'$ on which
$f^{S_{\hat \iota-2}}$
is a diffeomorphism.
The endpoints of $f^{S_{\hat \iota-2}}(H)$ 
consist of points of the form $d_j$ with $j<2i_0-2$
and so the diffeomorphic image
contains a $k(n)$-scaled neighbourhood of
$[\hat u_{n-2},u_{n-2}]$ where $k(n)\to \infty$
as $n\to \infty$.
It follows that
$f^{S_{\hat \iota-2}}|T'$
has uniformly bounded distortion
and that this distortion even disappears as $n$ tends to infinity:
$$
1-o(1/n)\le
\frac
{|Df^{S_{\hat \iota-2}}(x)|}
{|Df^{S_{\hat \iota-2}}(y)|}
\le 1+o(1/n)
$$
for each $x,y\in T'$.

So it suffices to describe the limit of the sequence of maps
$$f^{S_{\hat \iota-1}}\circ \dots \circ f^{S_{n-5}}\circ f^{S_{n-3}}$$
on $T_{n-1}$. So remember that
$$f^{S_{i-2}}\circ f^{S_i}=H_{i-2}^{-1}\circ \Psi_{i-2}\circ \Psi_i\circ
H_{i+2}= H_{i-2}^{-1}\circ \Theta_i\circ H_{i+2}.$$
Hence
\beq
f^{S_{\hat \iota-1}}\circ \dots \circ f^{S_{n-5}}\circ f^{S_{n-3}}=
H_{\hat \iota+1}^{-1} \circ \Theta_{\hat \iota+1} \circ \dots \circ
\Theta_{n-7}\circ \Theta_{n-3}
\circ H_{n-1}
\label{lode1}
\eeq
and each of these
maps is up to scaling near $\Psi_i$ because $i\ge i_0$.

Before continuing with our proof, let us describe the
idea. By our identification of
$z_i$ with the origin, the composition $f^{S_{i-2}}\circ f^{S_i}$
is identified with a map $\Theta_i$ having an almost neutral
fixed point. So we shall analyse a long composition
of such maps which are all close to a given map.  
Now an orientation preserving  one-dimensional map can
be essentially imbedded in a flow. So we shall be
able to get good estimates for a high iterate of the
map, by integrating a certain vector field explicitly.
Since the maps $\Theta_i$ are not all identical,
we shall only be able to use this comparison
up to a certain number of iterates.
We shall choose $n=\ell^{3/2}$ iterates, because it will
turn out that for the remaining iterates the relevant
restriction has a distortion which disappears
as $n$ and $\ell$ tend to infinity.
Now comes the miracle:
the limit of the first $n=\ell^{3/2}$ iterates
is a function which only depends on the coefficients
$\alpha_{0,\ell}$ and $\gamma_{0,\ell}$ and not
on the remainder of the Taylor series of $\Theta$.

To explain this and to analyze such a composition of maps
we make a digression to solutions of a particular
differential equation.

\begin{lemma}
Consider the following differential equation
\beq
x'(t)=-(\alpha/\ell)x(t)-\gamma\cdot [x(t)]^3,
\label{dv1}
\eeq
where $\alpha,\gamma>0$ are positive constants
and let $\phi_t(x)$ be its flow.
Then
\beq
\phi_t(x)=\frac{x\exp(-\alpha t /\ell)}
{\sqrt{(\gamma\ell/\alpha)\left[1-\exp(-2 \alpha t/\ell)\right]x^2+1}}.
\label{dv2}
\eeq
In particular, the Taylor expansion of $x\mapsto \phi_t(x)$
at $x=0$ is $x(1-\frac{\alpha t}{\ell})(1-\gamma't x^2+ \dots )$
where $\gamma'=\gamma+ o(1/\ell)$.
Moreover, there exist universal constants $C$ such that
\beq
\sum_{i=0}^\infty |\phi_i(x)| \le C \sqrt{\ell},\,\,
\sum_{i=0}^\infty |\phi_i(x)|^2 \le C \log(\ell)
\mbox{ and }
\sum_{i=0}^\infty |\phi_i(x)|^3 \le C\, ,
\label{dv3}
\eeq
provided $|x|\le 1$ (or $|x|$ is universally bounded).
\end{lemma}

\bigskip

\begin{remark} Essentially the reason for the
miracle mentioned above
is that for given $x,y\in [1/2,1]$
$$\lim_{t\to \infty}
\frac{\phi_t(x)}{\phi_t(y)}=
\frac{x}{\sqrt{\gamma \ell x^2 /\alpha +1}}
\frac{\sqrt{\gamma \ell y^2 /\alpha +1}}{y}
= 1 \pm o(1/\ell).$$
So an error in the initial condition becomes
less and less important as $t,\ell\to \infty$.
We should also remark that
$$\phi_t(x)=\pm \sqrt{M_t(x^2)}$$
where $M_t\colon \R^+\to \R^+$ is a Moebius transformation
which becomes increasingly degenerate as $t\to \infty$.
\end{remark}

\noindent
\pr
The general solution of (\ref{dv1}) is of the form
\beq
\frac{1}{(x(t))^2}=\frac{-\gamma \ell}{\alpha} +
\exp(2 \alpha t /\ell)\cdot c_0
\label{dv6}
\eeq
(where $c_0$ is an integration constant).
Indeed, differentiation of (\ref{dv6}) gives:
$$-\frac{2x'}{x^3}=\frac{2\alpha}{\ell}\exp(2\alpha t/\ell)\cdot c_0=
\frac{2\alpha}{\ell}\left[ \frac{\gamma}{\alpha}\ell
+ \frac{1}{x^2}\right].$$
This last expression is  (\ref{dv1}) rewritten.
It follows that the integration constant
is equal to $c_0=\frac{\gamma \ell}{\alpha}+\frac{1}{[x(0)]^2}$,
which gives the required expression.

Since
$$\frac{1}{\ell(1-\exp(-2\alpha t/\ell))}
\le \max(C/t,C/\ell),$$
the last three inequalities of this lemma
can be derived from (\ref{dv2}):
$$\sum_{i=0}^\infty |\phi_i(x)|
\le
C+ C\sum_{i=1}^\ell \frac{1}{\sqrt{i}} +
C\frac{1}{\sqrt{\ell}} \sum_{i=\ell+1}^\infty \exp(-a i /\ell)
\le C \sqrt{\ell} + C \frac{1}{\sqrt{\ell}} \ell.$$
$$\sum_{i=0}^\infty |\phi_i(x)|^2
\le
C+ C\sum_{i=1}^\ell \frac{1}{i} +
C\frac{1}{\ell} \sum_{i=\ell+1}^\infty \exp(-a i /\ell)
\le C \log(\ell) + C \frac{1}{\ell}\ell .$$
$$\sum_{i=0}^\infty |\phi_i(x)|^3
\le
C+ C\sum_{i=1}^\ell \frac{1}{i^{3/2}} +
C\frac{1}{\ell^{3/2}} \sum_{i=\ell+1}^\infty \exp(-a i /\ell)
\le C .$$
\qed

Let us write
$\theta_i=\Theta_{n-2i-1}$.
Then
\beq
\Theta_{\hat \iota+1}\circ \dots \circ \Theta_{n-7}\circ \Theta_{n-3}=
\theta_m\circ \dots \theta_1
\label{lode2}
\eeq
where $m=(n-\hat \iota-4)/4$.
This brings us to the following
abstract situation, where are before $\phi_t$ is the solution
of the previous differential equation~\ref{dv2}:
\bigskip

\begin{theo}\label{comdv}
Consider a sequence of analytic maps
$\theta_{i,\ell}\colon \rz\to \rz$ for $i=0,1,\dots$ 
such that
$$\theta_{i,\ell}(x)=
\left(1-\frac{\alpha_{i,\ell}}{\ell}\right)\,
x+\frac{\beta_{i,\ell}}{\ell}\, x^2-
\gamma_{i,\ell}\, x^3 + O(|x|^4).$$
Assume that
\begin{itemize}
\item
$\alpha_{0,\ell}$, $\gamma_{0,\ell}$ 
are positive, uniformly bounded away from zero and bounded from above
and that $\beta_{0,\ell}$ uniformly bounded in norm;
\item
$\theta_{i,\ell}$ is a diffeomorphism from $(-1,1)$ into itself with
$|\theta_{i,\ell}(x)|<|x|$
and such that for each given $\epsilon>0$ there exists
$n$ (which does not depend on $i$ and $\ell$)
such that $\theta_{i,\ell}^n(-1,1)\subset (-\epsilon,\epsilon)$;
\item for each
$|\alpha_{i,\ell}-\alpha_{0,\ell}|<C/\ell$,
$|\beta_{i,\ell}-\beta_{0,\ell}|<C/\ell$,
$|\gamma_{i,\ell}-\gamma_{0,\ell}|<C/\ell$
and $\theta_{i,\ell}$ is in a compact set of maps.
\end{itemize}
Then writing $F_m=\theta_m\circ \dots \circ \theta_1$
one has for each $m\ge m(\ell)=\ell^{3/2}$
and each $x\in (-1,1)$,
\beq
1-o(1/\ell)\le \frac{F_m(x)}{\phi_{m}(x)}\le 1+o(1/\ell)
\label{dv7}
\eeq
and
\beq
1-o(1/\ell)\le \frac{DF_m(x)}{D\phi_{m}(x)}\le 1+o(1/\ell),
\label{dv8}
\eeq
where $o(s)$ is some universal function which tends
to zero as $s\to 0$
and where $t \mapsto \phi_t$ is the flow
of the differential equation from the previous lemma, with
$\alpha=\alpha_{0,\ell}$ and $\gamma=\gamma_{0,\ell}$.
\end{theo}

Because of (\ref{lode1}) and (\ref{lode2}),
Theorem~\ref{miracle2} follows from this theorem.
So it remains to prove Theorem~\ref{comdv}.
For this we need the following two lemmas.
\medskip

\begin{lemma}
\label{comdv1}
The assertion of the previous proposition holds
for $m=\ell^{3/2}$ in which case one also has
for $|x|<1$,
\beq
|F_m(x)|<\frac{1}{\ell}.
\label{dv9}
\eeq
\end{lemma}
\pr
In order to be definite choose $x>0$. The case that $x<0$
goes similarly.

\medskip

\noindent
Step 1. First we claim that if
\beq
0<x<\frac{1}{\sqrt{j}}
\label{dv10}
\eeq
and $i\in \nz$ then there exists
$t'<1<t$ with
\beq
|t-1|,|t'-1|\le C\max\left(\frac{1}{\sqrt{j}},\frac{1}{\ell}\right)
\label{dv11}
\eeq
such that
\beq
\phi_t(x)\le \theta_i(x) \le \phi_{t'}(x).
\eeq
Here as before $C$ is a universal constant
(not depending on $j,i,\ell$).
This can be seen as follows: if $O(x)$ is a bounded function
then
\beq
x^4 \cdot O(x) \le |t'-1| \cdot |x|^3
,
\beta |x|^2\le |t'-1|\cdot |x|
\mbox{ and }
\frac{C_0}{\ell}\le |t'-1|
\label{dv12}
\eeq
provided $x$ is as in (\ref{dv10}) and $t'<1$
so that equality holds in (\ref{dv11}).
Therefore, taking $t'<1$ in this way,
we get from (\ref{dv12}) and from the second assumption of this lemma that
$$\theta_{i,\ell}(x)=
\left(1-\frac{\alpha_{i,\ell}}{\ell}\right)\,
x+\frac{\beta_{i,\ell}}{\ell}\, x^2-
\gamma_{i,\ell}\, x^3 + x^4\cdot O(x)$$
is bounded from above by
$$
x(1-\frac{\alpha t'}{\ell})(1-\gamma' t' x^2+ \dots )
=
\frac{x\exp(-\alpha t' /\ell)}
{\sqrt{(\gamma\ell/\alpha)\left[1-\exp(-2 \alpha t'/\ell)\right]x^2+1}}
=\phi_t(x)
$$
where $\gamma'=\gamma+ o(1/\ell)$ and where we take
$\alpha=\alpha_{0,\ell}$, $\gamma=\gamma_{0,\ell}$.
As usual, $O(x)$ is some bounded function of $x$.

\medskip

\noindent
Step 2. By assumption $\theta_i(x)\in (0,x)$
and there exists a universal number $j$ such that
$\theta_j \circ \dots \circ \theta_1(x)$ is
inside a given neighbourhood of $0$.
So it follows from the previous step
that for $i$ sufficiently large:
$$F_i(x)=\theta_i\circ \dots \circ \theta_1(x) \le \phi_{i/2}(x).$$
Using the explicit formula for $\phi_t(x)$
one sees that $F_i(x)\le \phi_{i/2}(x)\le C \frac{1}{\sqrt{i}}$.

\medskip

\noindent
Step 3. Using Step 2,
$F_i(x)\le C \frac{1}{\sqrt{i}}$
and therefore
there exists $t'_i<1<t_i$ as in Step 1 such that
$$\phi_{t_i}\circ F_i(x)\le
\theta_{i+1}\circ F_i(x)=F_{i+1}(x)\le
\phi_{t_{i+1}}\circ F_i(x).$$
By induction we get that
\beq
\phi_{T_i}(x)\le F_i(x)\le \phi_{T'_i}(x)
\label{betw}
\eeq
with $T'_i<i<T_i$ such that $T_i=t_i+\dots + t_1$,
$T'_i=t'_i+\dots + t'_1$  and $|t_i'-t_i|$ as above.
Hence
$$|T_i-T'_i|\le \sum_{j=0}^i|t_i-t_i'|
\le
C\left[ \sum_{j=0}^i \left(\frac{1}{\sqrt{j}} + \frac{1}{\ell}\right)
\right].$$
For $i\le \ell^{3/2}$ this gives
$$|T'_i-T_i|\le
C[\sqrt{i}+i/\ell]\le C[\ell^{3/4} + \ell^{1/2}].$$
Since
$$\phi_t(x)=\frac{x\exp(-\alpha t /\ell)}
{\sqrt{(\gamma\ell/\alpha)\left[1-\exp(-2 \alpha t/\ell)\right]x^2+1}}$$
it follows that if $i\le \ell^{3/2}$
for such $T=T_i, T'=T_i'$,
$$
\frac{\phi_{T'}(x)}{\phi_{T}(x)}
\le \exp(\alpha |T-T'| /\ell)
\frac
{\sqrt{(\gamma\ell/\alpha)\left[1-\exp(-2 \alpha T/\ell)\right]x^2+1}}
{\sqrt{(\gamma\ell/\alpha)\left[1-\exp(-2 \alpha T'/\ell)\right]x^2+1}}.
$$
Because
$$\sqrt{(\gamma\ell/\alpha)\left[1-\exp(-2 \alpha T'/\ell)\right]x^2+1}
\ge C_0\sqrt{\min(|T'|,\ell)}$$
this implies that
$$
\frac{\phi_{T'}(x)}{\phi_{T}(x)}
\le \left( 1+ C\frac{|T-T'|}{\ell}\right)
\left( 1+ C\ell\frac{\left[
\exp(-2\alpha T/\ell)-\exp(-2\alpha T'/\ell)\right]}{\sqrt{\min(|T'|,\ell)}}
\right).$$
For $i\le \ell^{3/2}$, $T'<i<T$ with $|T'-T|\le \sqrt{i}$,
the last factor is at most
$$
\le 1+ C\frac{|T-T'|}{\sqrt(|T'|,\ell)}\le
1+ \max\left(\frac{1}{\sqrt{i}},\frac{1}{\sqrt{\ell}}\right)
.$$
From this it follows that
\beq
\frac{\phi_{T'}(x)}{\phi_{T}(x)}
\le 1+ C\left[\frac{\sqrt{i}}{\ell} + \frac{1}{\sqrt{i}}
+ \frac{1}{\sqrt{\ell}}\right]
\label{dv13}
\eeq
for $i\le \ell^{3/2}$ and the the upper bound in (\ref{dv7}).
A lower bound is shown in the same way.
Therefore, using (\ref{betw}), one gets
\beq
\big| \frac{F_i(x)}{\phi_i(x)}-1
\big|
\le C\left[\frac{\sqrt{i}}{\ell} + \frac{1}{\sqrt{i}}
+ \frac{1}{\sqrt{\ell}}\right]
\label{dv14}
\eeq
which gives the first inequality (\ref{dv7}) claimed in the
theorem.
Note that the right hand side is  at most
$C \ell^{-(1/4)}$ when $i\le \ell^{3/2}$.

\medskip

\noindent
Step 4. Now we shall prove the second inequality (\ref{dv8})
claimed in the theorem. For this note
\beq
\frac{D\theta_{i+1}(F_i(x))}{D\phi_{1}(\phi_i(x))}
\le
\left( 1+ C|D\theta_{i+1}(F_i(x))-D\phi_1(F_i(x))|
+ C|D\phi_1(F_i(x))-D\phi_1(\phi_i(x))|\right).
\label{dies}
\eeq
This expression is at most
$$
\left( 1 + C/\ell^2 + (C /\ell) |F_i(x)|
+ C\left[|F_i(x)|^2-|\phi_{i}(x)|^2 \right] + C |F_i(x)|^3 \right) .$$
By the last part of the lemma,
$$\sum |F_i(x)|\le \sum |\phi_i(x)|\le  \sqrt{\ell}.$$
Moreover, by (\ref{dv14})
$$
\left[|F_i(x)|^2-|\phi_{i}(x)|^2 \right]
\le C\cdot
\left[\frac{|F_i(x)|}{|\phi_{i}(x)}-1| \right]
\cdot |\phi_i(x)|^2$$
whne $i\le \ell^{3/2}$
which gives that
$$
\sum_{0}^{\ell^{3/2}}\left[|F_i(x)|^2-|\phi_{i}(x)|^2 \right]
\le C\cdot
\sum_{0}^{\ell^{3/2}}
\left[\frac{\sqrt{i}}{\ell} + \frac{1}{\sqrt{i}}
+ \frac{1}{\sqrt{\ell}}\right] \,\,
\left( \frac{1}{\sqrt{i}}\right)^2
\le C \frac{\sqrt{\ell^{3/2}}}{\ell}\le C \ell^{-1/4}.$$
Combining this proves (\ref{dv8}).
\qed

\medskip

The remaining iterates have a non-linearity which vanishes
as $\ell\to \infty$:

\begin{lemma}
\label{comdv2}
Let $\theta_{i,\ell}$ be as in the previous lemma.
Assume that $|y|<1/\ell$.
Then for any $m,k\in \nz$
one has, writing  $F_{m,k}=\theta_{m+k}\circ \dots \circ \theta_{m+1}$,
$$1-o(1/\ell)
\le \frac{|D \hat F_{m,k}(y)|}{|D\hat F_{m,k}(0)|}\le 1+o(1/\ell).$$
This means that the distortion of $F_{m,k}$ on $[-1/\ell,1/\ell] $
is small for large $\ell$.
\end{lemma}
\pr
One has
$$|D\theta_i(y)-D\theta_i(0)|\le \frac{\beta}{\ell}|y|+\gamma|y|^2 + O(|y|^3).$$
Since $|y|<1/\ell$, therefore
$$|D\theta_i(y)-D\theta_i(0)|\le \frac{C}{\ell}|y|$$
and
$$|D\theta_i(z)|\le (1-C_0/\ell)$$
for each $z\in (0,y)$.
It follows that
$$
\left|
\frac{|D\hat F_{m,k}(y)|}{|D\hat F_{m,k}(0)|}-1 \right|
\le \frac{C}{\ell} \sum_{i=0}^{k}|F_{m,i}(y)| \le
\frac{C}{\ell} \sum_{i=0}^m \frac{1}{\ell}
(1-C_0/\ell)^i \le \frac{C}{\ell}.$$
\qed

\noindent
{\em Proof of Theorem~\ref{comdv}:}
Take as before $m(\ell)=\ell^{3/2}$ and write
for $m\ge m(\ell)$, $F_m=F_{m,m(\ell)}\circ F_{m(\ell)}$.
Because of Lemma~\ref{comdv1}, one has that the
first map $F_{m(\ell)}$ can be compared very well
with $D\phi_{m(\ell)}$. Since this lemma also asserts
that $F_{m(\ell)}(-1,1)\subset (-1/\ell,1/\ell)$,
the last map $F_{m,m(\ell)}$ can be compared
better and better (as $\ell$ tends to infinity) with a
linear map because of Lemma~\ref{comdv2}.
Combined, this gives the required estimates.
\qed

\sect{The proof of the Main Theorem}
In this section we shall complete the proof of the Main Theorem.
First we should remark that the filled Julia set of
$f$ is nowhere dense. This simply follows from the fact that
the critical point is recurrent, $c\in \omega(c)$, and therefore
$f$ has no periodic attractors or neutral periodic points,
see for example \cite{Blan}, \cite{L0} or \cite{Mil}.
Actually, this also implies that the Julia set of $f$
is connected and that  $\cup f^{-k}(c)$ is dense
in the Julia set. Rather than showing that the Julia set of $f$
has positive Lebesgue measure we shall prove the following
stronger theorem. (In fact, this theorem is equivalent to
the main theorem because of \cite{L0}, \cite{L2}.)

\begin{theo}\label{mainth}
For all $x$ from a set of positive Lebesgue measure
$$\omega(x)\subset \omega(c).$$
\end{theo}
 
\par
Let $A_k,A_k'$, $F$, $\A_n$ be defined as in the previous section.
In this section we shall show that $F$ and $\A$ satisfy the
assumptions of Theorem~\ref{random}.

\begin{theo}\label{drift}
For all sufficiently large $\ell$ holds:
The set $D$ of all points $x$
for which the trajectory $(F^kx)_{k>0}$ visits
$A_n$ at most finitely often, has positive Lebesgue measure.
\end{theo}
Let us first show that this result implies our Theorem~\ref{mainth}.
\bigskip

\noindent
{\em Proof of Theorem~\ref{mainth}}.
Because $f^{-k}(c)$ is dense in the Julia set,
the length of a maximal of interval of monotonicity of $f^{S_n}|\rz$
is at most $\delta(n)$ where $\delta(n)\to 0$ as $n\to \infty$.
Consider a point $x\in X$ for which $(F^kx)_{k>0}$
visits each annulus $A_n$ at most finitely often, and denote by
$t_1<t_2<t_3<\ldots$ the 
sequence of times for which $F^kx=f^{t_k}x$. We have to show that
$\lim_{t\to\infty}\mbox{dist}(f^tx,\omega_f(c))=0$. Along the
subsequence $t_k$ 
this holds as
$\lim_{k\to\infty}f^{t_k}x=\lim_{k\to\infty}F^kx=c\in\omega_f(
c)$.
Consider now $t_k<t<t_{k+1}$ and suppose that $F^kx\in A_n\subset D_{n+1}$.
Write $i:=t-t_k$ and note that $0< i< S_k$.
Now $f^{S_{n}-1}$ is a diffeomorphism from
$f(D_{n+1})$ to $D_{n-1}\subset D_*(u_{n-2},\hat u_{n-2})$
and using the Lemma of Schwarz $f^i(A_n)\subset f^i(D_{n+1})$
is contained inside a disc $D_*(h_i,j_i)$ where $(h_i,j_i)\subset \rz$
is an interval of monotonicity of $f^{S_n-1}$.
So the diameter of this disc is at most $\delta(n)$
and since $c_{S_{n+1}}\in A_n$ it follows that
$$\text{dist}(f^ix,\omega(c))\le \text{diam}(f^i(A_n))\le
\delta(n)\to 0\text{ as }n\to \infty.$$
\qed
\par

\medskip

\noindent
{\em Proof of Theorem~\ref{drift}:\quad}
Let $k_0$ be a large integer as before, for $k\ge k_0$ let
$X_k$ be the disjoint union of $A_k$ and $A_k'$ and
finally let $X=\cup_{k\ge k_0}X_k$.
So $X$ can be considered as the disjoint union of
the disc $\cup A_k$ and the disc with holes $\cup A_k'$.
Define $\X$ to be the partition of $X$ in elements $X_k$
and let $m$ the Lebesgue measure on $X$
(i.e., the Lebesgue measure on $\cup A_k$ and on $\cup A_k'$).
Let us show that we can apply Theorem~\ref{random}
with these choices. So take $n_0\ge k_0+2$,
take $A\in \X_{k+1}$ and assume that $n:=\phi(F^k(A))\ge n_0$.
Then $F^k$ either maps $A$ onto $\tilde A_n$
where $\tilde A_n$ is the annulus $A_n$ or the
annulus $A_n'$ with a hole.

Because of Theorem~\ref{extend}, the map $F^k$ extends in a
univalent way to a map which maps respectively onto
the slit $Slit_n$ or onto the slit $Slit_n'$.
Because of Theorem~\ref{slitkoebe} and because of
the Koebe Lemma this implies that $F^k\colon A\to \tilde A_n$
has uniformly bounded distortion.

Now we will check that for $\ell$ sufficiently large
the first condition
\beq
\int_A(\Delta\ph-1)\circ F^k\, dm \ge 0
\eeq
from the random walk result, Theorem~\ref{random}, is satisfied,
where
\[
\Delta\ph:=\ph\circ F-\ph\
\]
and $\ph(x)=i$ if $x\in A_i$.
In other words, we need to show
that
\beq
\frac{1}{m(A)}\int_A \Delta\ph\circ F^k\, dm \ge 1.
\label{pro1z}
\eeq
Let $A^i$ be the part of $A$ which $F^k$ sends to $\tilde A_n^i$.
Note that this is the part of $A_n$ which is between two infinite rays
$l_i$ and $l_{i+1}$; also note that this piece is connected since
$A_n\cap \rz$ is equal to $[u_n,\hat u_n]\setminus [u_{n+1},\hat u_{n+1}]$
and since $A_n$ is rotational symmetry.
From Theorem~\ref{maco} one has that the diffeomorphism
$F^k\colon A^i\to A_n^i$ has uniformly bounded distortion.
Since $f(A_n^i)=f(A_n)$ this implies that (\ref{pro1z})
follows from
\beq
\frac{1}{m(\tilde A_n)}
\int_{\tilde A_n}\Delta\ph\, dm
\ge \gamma(\ell).
\label{pro1a}
\eeq
where $\gamma$ stands for some function such that
$\gamma(\ell)\to \infty$ as $\ell\to \infty$.
Now let $\hat A_n=A_n\setminus f^{-1}(D_{n+1}^{1,f})$.
So $\hat A_n$ is the annulus $A_n$ with {\it all} the
$\ell$ symmetrics of $D_{n+1}^1$ removed (instead of only
one such disc $D_{n+1}^1$ removed as is the case with $A_n'$).
Using this notation,
\beqas
\int_{A_n}\Delta\ph\, dm &=&
\int_{\hat A_n}\Delta\ph\, dm +
\ell \int_{D_{n+1}^1}\Delta\ph\, dm\\
&\ge& -2 m(A_n) + \ell \int_{D_{n+1}^1}\Delta\ph\, dm
\eeqas
and
\beqas
\int_{A_n'}\Delta\ph \, dm &=&
\int_{\hat A_n}\Delta\ph \, dm +
(\ell-1) \int_{D_{n+1}^1}\Delta\ph \, dm\\
&\ge& -2 m(A_n) + (\ell-1)
\int_{D_{n+1}^1}\Delta\ph \, dm
\eeqas
where we have used $\Delta\ph\ge -2$.
From this it follows that it is enough to prove that (\ref{pro1a})
holds for $\tilde A_n=A_n$ and therefore it suffices to show that
\beq
\frac{1}{m(A_n)}
\int_{A_n}(\ph\circ F-n)\, dm
=\frac{1}{m(A_n)}
\int_{A_n}\Delta\ph\, dm
\ge \gamma(\ell)
\label{pro1b}
\eeq
Note that
$$F|A_n=f^{S_n}=f^{S_{n-2}}\circ f^{S_{n-1}}$$
and that
$f|A_n$ has uniformly bounded distortion
since $A_n$ is between two discs centered at $0$
of radius $|u_n|$ respectively $|u_n|(1-C/\ell)$.
This implies that it is enough to show that
$$
\frac{1}{m(A_n^f)}
\int_{A_n^f}(\ph\circ f^{S_n-1}\,-\, n)\, dm
\ge \gamma(\ell).
$$
(Here $\gamma$ is a function with the same properties as before.)

\kies{

\begin{figure}[htp] \hfil
\beginpicture
\dimen0=0.3cm
\setcoordinatesystem units <\dimen0,\dimen0> point at 0 0
\setplotarea x from -8 to 9, y from -7 to 7
\setlinear
\plot -9 0 9 0 /
\put {\small $D_{n-1}$} at 5.9 5.5
\setsolid
\circulararc 360 degrees from 7 0  center at  0 0
\setdashes
\circulararc 360 degrees from 6.2 0  center at  0 0
\setdots <2pt>
\circulararc 360 degrees from 4 0  center at  0 0
\setsolid
\circulararc 360 degrees from -4 0  center at  -4.7 0
\setlinear
\setdots <1pt>
\plot -4.1 -0.2 -4.1 0.2 /
\plot -4.2 -0.3 -4.2 0.3 /
\plot -4.4 -0.5 -4.4 0.5 /
\plot -4.7 -0.7 -4.7 0.7 /
\plot -5   -0.5 -5 0.5   /
\plot -5.2 -0.3 -5.2 0.3 /
\plot -5.3 -0.2 -5.3 0.2 /
\multiput {\small $\bullet$} at -7 0 -4 0 7 0 /
\multiput {\small $o$} at -9 0 0 0 -4.7 0  9 0 /
\multiput {$*$} at -5.4 0 4 0 /
\put {\small $D_n$} at 0 5.3
\put {\small $D_{n+1}$} at 0 -3
\put {\small $D_n^1$} at -4.7 1.7
\put {\small $A_{n-2}=D_{n-1}\setminus D_n$} at -13 -7
\put {\small $D_n^1\subset A_{n-1}=D_{n}\setminus D_{n+1}$} at -14 -8.7
\put {\small $A_i\subset D_{n+1}\,\,\, , i\ge n$} at -13 -10.4
 
\setcoordinatesystem units <\dimen0,\dimen0> point at 0 15
\setplotarea x from -8 to 9, y from -7 to 7
\setlinear
\plot -9 0 9 0 /
\put {\small $D_{n+1}^f$} at 4.7 4.7
\setsolid
\circulararc 360 degrees from 6 0  center at  0 0
\setsolid
\circulararc 360 degrees from 2.2 0  center at  0 0
\setdots <2pt>
\circulararc 360 degrees from 2.2 0  center at  3.3 0
\multiput {\small $\bullet$} at -6 0 2.2 0 6 0  /
\multiput {\small $o$} at -1 0 /
\multiput {$*$} at -2.2 0  4.4 0 /
\put {\small $A_n^f$} at -4 1
\put {\small $D_{n+2}^f$} at 0.4 3
\put {\small $D_{n+1}^{1,f}$} at 3.3 -1.6
\setquadratic
\setshadegrid span <1pt>
\vshade -2.2 0 0
<,z,,> -2 -0.91 0.91 -1.75 -1.33 1.33
<z,z,,> -1.5 -1.6 1.6 -1 -1.95 1.95
<z,z,,>   0 -2.2 2.2 1 -1.95 1.95
<z,z,,> 1.5 -1.6 1.6 1.75 -1.33 1.33
<z,,,> 2  -0.91 0.91  2.2  0  0 /
\endpicture
\caption[ ]{{\small The annulus $A_n^f$ and its image under the
map $f^{S_n-1}$. Moreover, $B_i^i\subset D_{n+1}^{1,f}$ is the preimage
of the region $A_{n+i}^i\subset A_{n+i}$.}}
\end{figure}

}

Since $(\ph\circ f^{S_n-1}\,-\, n)\ge -2$,
since $(\ph\circ f^{S_n-1}\,-\, n)\ge 0$ on $D_{n+1}^{1,f}$
and since the area of $D_{n+1}^{1,f}\subset A_n^f$ occupies a definite
proportion of the area of $A_n^f$ it suffices to show that
\beqa
&& \frac{1}{m(D_{n+1}^{1,f})}
\sum_{i\ge 0}
\,\, i \cdot m(\{x\in D_{n+1}^{1,f}\st \ph(\circ f^{S_n-1}x)=n+i\}) \\
&& \quad = \,\,
\frac{1}{m(D_{n+1}^{1,f})}
\int_{D_{n+1}^{1,f}}(\ph\circ f^{S_n-1}\,-\, n)\, dm
\ge \gamma(\ell).
\label{gga}
\eeqa
For $i\ge 1$, write as before
$$A_{n+i}^i=A_{n+i}\cap \{z\in \cz \st |\arg(z)|< i/\ell\}.$$
For $i$ small this is a very small piece of the annulus $A_{n+i}$
but for $i\approx \ell$, this piece occupies a definite proportion
of $A_{n+i}$. Let $B_i^+$ be the preimage of this set under the map
$f^{S_n-1}\colon D_{n+1}^{1,f}\to D_{n+1}=\cup_{i\ge 0}A_{n+i}$:
$$B_i^+=\{x\in D_{n+1}^{1,f}\st f^{S_n-1}(x)\in A^i_{n+i}\}.$$
Then (\ref{gga}) is implied by
\beq
\frac{1}{m(D_{n+1}^{1,f})}
\sum_{i\ge 0}\,\, i \cdot m(B^+_i)
\ge \gamma(\ell)
\label{pro1d}
\eeq
Because of the real bounds and the shape estimates
on the annuli $A_j$, the diameter of $A_{n+i}^i$ is of the
same order as its distance to the nearest critical value of
$f^{S_n-1}$. Hence the distortion of the restriction of
$f^{S_n-1}\colon D_{n+1}^{1,f}\to D_{n+1}$
to the diffeomorphism
$$f^{S_n-1}\colon B_i^+\to A_{n+i}^i$$
is uniformly bounded (for all $i\ge 1$ and
all large $\ell$ and $n$).
From Theorem~\ref{miracle} we have very good estimates
for this map (on the real line) and
combined with the bounded distortion on $B_i^+$
this gives a uniform constant $C>0$ such that
\beq
|Df^{S_n-1}(x)|\le C \frac{i^{3/2}}{\ell}
\frac{|u_{n-2}-\hat u_{n-2}|}{|r_n^f-u_n^f|}
\label{esst}
\eeq
for each $x\in B_i^+$.
Since $f^{S_n-1}$ maps $B^+_i$ diffeomorphically onto $A^i_{n+i}$
and areas are distorted with the {\it square} of the
Jacobian of the map, (\ref{esst}) implies
\beq
\frac {m(A^i_{n+i})} {m(B^+_i)}
\le
\frac{i^3}{\ell^2}
\left[\frac{|u_{n-2}-\hat u_{n-2}|}{|r_n^f-u_n^f|}\right]^2.
\label{pro1e}
\eeq
Here we have used (\ref{esst}) and that the distortion of the size of areas
is measured by the square of the Jacobian of the map.
By the corollary to Theorem~\ref{complb}, the `height'
of $D_{n+1}^{1,f}$ is comparable to its `width' and also comparable
to $|r_n^f-u_n^f|$, and in particular,
\beq
\frac{m(D_{n+1})}{m(D_{n+1}^{1,f})}
\ge K \left[\frac{|u_{n-2}-\hat u_{n-2}|}{|r_n^f-u_n^f|}\right]^2
\label{pro1f}
\eeq
for some uniform constant $K$.
Combining (\ref{pro1e}) and (\ref{pro1f}), we get
\beq
\frac{m(B_i^+)}{m(D_{n+1}^{1,f})}\ge C\cdot \frac{\ell^2}{i^3}\cdot
\frac{m(A^i_{n+i})}{m(D_{n+1})}.
\label{lasteq}
\eeq
Since for $i=1,2,\dots,\ell$ the area of $m(A^i_{n+i})$ is
of the order $i/\ell^2$ times the area of $m(D_{n+1})$,
the last inequality yields
\beq
\frac{m(B_i^+)}{m(D_{n+1}^{1,f})}\ge C\cdot \frac{1}{i^2},
\mbox{ for }i=1,2,\dots,\ell .
\label{logdr}
\eeq
(For $i=\ell,\ell+1,\dots$, the area of $m(A^i_{n+i})$
is of the order $e^{-i/\ell}/\ell$ times the area of
of $m(D_{n+1})$ and so (\ref{lasteq})
also gives
\beq
\frac{m(B_i^+)}{m(D_{n+1}^{1,f})}\ge C\cdot \frac{\ell}{i^3}\cdot
e^{-i/\ell}
\mbox{ for }i=\ell+1,\ell+2,\dots .)
\label{logdr2}
\eeq
Hence, because of (\ref{logdr}), the
left hand side of (\ref{pro1d}) can be bounded from below
by
\beq
\sum_{i=0}^\ell \,\, \frac{1}{i} \ge
\text{const} \, \cdot \, \log(\ell).
\label{somdr}
\eeq
(The contribution to the expected
drift due to the `tail' terms corresponding to $i\ge \ell+1$,
see (\ref{logdr2}), is
$$
\sum_{i=\ell+1}^\infty \,\, \frac{\ell}{i^2}\cdot e^{-i/\ell}
$$
which is uniformly bounded and so does not give any essential
improvement on our previous bound.)
This concludes the proof of the first assumption of
Theorem~\ref{random}.

Next we check the second condition
from Theorem~\ref{random}.
This means that we have to find $M<\infty$ such that
$$
\frac{1}{m(A)}
\int_A(\Delta\ph)^2 \circ F^k\, dm \le M
$$
for any $A\in \A_{k+1}$ with $n=\phi(F^kA)$ large enough.
We should emphasize that $M$ does {\it not} need to be
uniform in $\ell$. Since $F^k\colon A\to A_n$
has as before bounded distortion,
it is enough to prove that there exists $M$ such that
$$
\frac{1}{m(A_n)}\int_{A_n}(\Delta\ph)^2\, dm \le M
$$
for each $n$ sufficiently large. Since the distortion of $F$ on $A_n$
is bounded by $\ell^2$ (which is bounded),
the last hand side of the previous expression is at most
$$\ell^2
\frac{1}{m(D_{n-1})}
\int_{D_{n-1}}(\ph-n)^2\, dm\,\, = \,\,
\frac{1}{m(D_{n-1})}\sum_{i=-2}^{\infty} i^2\cdot m(A_{n+i}).
$$
Since
$$\frac{m(A_{n+i})}{m(D_{n-1})}
\le C_1\cdot \frac{1}{\ell}\cdot e^{-C\cdot i/\ell},$$
the last infinite sum is, up to a multiplicative constant,
bounded from above by
$$\ell^2 \sum_{i=-2}^{\infty}
i^2 \cdot \frac{1}{\ell}\cdot e^{-C\cdot i/\ell}
\le \text{Const} \cdot \ell^5$$
which proves the second condition of Theorem~\ref{random}
and concludes the proof of the Main Theorem.
\qed
\bigskip
\medskip

To conclude this section, we would like to make
some comments on the difficulties
of giving a computer supported numerical
`estimate' for the smallest value of $\ell$ for which
the statement of the theorem holds. So take
$A\in X_{k+1}$ with $\phi(A)=n$
and let $\{W_i\}_{i\ge -2}$ be the partition of $A$ defined by
the amount of drift:
$$W_i=\{x\in A\st \Delta \phi(x)=i\}.$$
First note that
$$\frac{m(W_{-1})}{m(A)}\mbox{ and }\frac{m(W_{-2})}{m(A)}
$$
are uniformly bounded away from zero. So the chance to
go `down' is uniformly bounded away from zero.
Moveover, consider the expected drift
\beq
\frac{1}{m(A)}\sum_{i\ge -2} i\cdot \mu(W_i) =
\frac{1}{m(A)}\int_A \Delta \phi \, dm.
\label{edr}
\eeq
Suppose we would estimate this term by choosing a finite, say $k'$,
number of states neighbouring a given state and estimate
numerically the finite sum
\beq
\frac{\sum_{i\ge -2}^{k'} i\cdot \mu(W_i')}{m(A)}
\label{smu}
\eeq
where $W_i'=W_i$ for $i<k'$ and $W_{k'}'=\cup_{i\ge k'}W_i$.
By (\ref{logdr}) we see that
$$\frac{m(W_i)}{m(A)}\ge \frac{C}{i^2}$$
(in fact, one can show that the left hand side is really
of this form) and so we only can expect (\ref{smu})
to be equal to
$$\sum_{i=0}^{k'}\frac{1}{i}$$
(at least provided $k'<\ell$) and so
the estimate we would obtain in this way does not
get better for increasing $\ell$.
This means that the only way to get a good esimate for (\ref{edr})
is to take $k'$ very large!

Similarly, one has to take $k$ very large before one has
that
$$\frac{|\{x\in A\st \phi(F^k(x))-\phi(x)\ge 1\}|}{|A|}$$
gets close to one. This means that one has to iterate the induced
map a very large number of times, before `observing' the drift.

Let us finally also make a comparison
with the real one-dimensional paper \cite{BKNS},
where we did not need to take the square in
(\ref{pro1e}). This means that -- if we had known the
estimates of Section~\ref{asytex} of this paper already in that paper --
we would have obtained (up to a multiplicative constant which
is universally bounded away from zero), the following lower bound
for the expression (\ref{logdr}) in the {\it real one-dimensional}
case.
Firstly, then
$$\frac{m(W_i)}{m(A)}\sim \frac{1}{\ell}
\cdot \frac{\ell}{i^{3/2}}=\frac{1}{i^{3/2}},\mbox{ for }i=1,2,\dots,\ell,$$
where $1/\ell$ corresponds
to the size of the one-dimensional annulus $A_{n+i}$ relative to
the size of $A_n$ and $\ell/i^{3/2}$ corresponds to the
distortion of measure (when $i<\ell$).
Similarly,
$$\frac{m(W_i)}{m(A)}\sim \frac{e^{-i/\ell}}{\ell}
\cdot \frac{\ell}{i^{3/2}}=
\frac{e^{-i/\ell}}{i^{3/2}},\mbox{ for }i=\ell+1,\ell+2,\dots.$$
Hence the expected drift in the real case is
\beq
\sum_{i=0}^\ell i\cdot
\frac{1}{i^{3/2}} +
\sum_{i=\ell+1}^\infty i\cdot
\frac{e^{-i/\ell}}{i^{3/2}}
\approx C_1\sqrt{\ell} +  C_2\sqrt{\ell}
\label{onedi}
\eeq
So the drift grows in the real case much faster with $\ell$
then in the complex case!
Moreover, note that the contribution due to the `tail' (i.e., the
term of the form $\sum_{\ell+1}^\infty$ in (\ref{onedi})
which can be derived immediately from Proposition~\ref{rootl})
already cause a large drift when $\ell$ is large.
Instead, of the above estimates we used in \cite{BKNS}
a simple Koebe estimates which implies that the measure distorts by at most
$\ell/i^2$ and so we were merely able to get the weaker bound
$$\sum_{i=0}^\ell i\cdot
\frac{\ell}{i^{2}}\frac{1}{\ell}=\sum_{i=0}^\ell
\frac{1}{i}\approx \log(\ell).$$
This -- not so sharp estimate --
is sufficient in the one-dimensional case.
In our `two-dimensional' case it is not enough because
then the $i^2$ term in (\ref{logdr}) would have to be replaced by
$i^3$ and the term
$i$ would become $i^2$ in
(\ref{somdr}). This would lead to a series with uniformly bounded sum. Hence
the expected drift would be of constant magnititude. Whether it would be
positive or not would depend on the numerical value of some constants.


\end{document}